\documentclass[10pt,francais]{smfart}

\usepackage[dvips]{epsfig}

\usepackage[francais]{babel}

\usepackage{smfthm}

\newcommand{\C}{\mathbb C}
\renewcommand{\P}{\mathbb P}
\newcommand{\N}{\mathbb N}
\newcommand{\Z}{\mathbb Z}

\renewcommand{\L}{\mathbb L}

\newcommand{\Q}{\mathbb Q}

\newcommand{\G}{\mathbb G}

\newtheorem{theorem}{Th\'eor\`eme}[section]
\newtheorem{lemma}[theorem]{Lemme}
\newtheorem{proposition}[theorem]{Proposition}
\newtheorem{corollary}[theorem]{Corollaire}
\newtheorem{rappel}[theorem]{}
\newtheorem{problem}[theorem]{Probl\`eme}

\theoremstyle{definition}
\newtheorem{definition}[theorem]{D\'efinition}
\newtheorem{remark}[theorem]{Remarque}

\newtheorem{notation}[theorem]{Notation}
\newtheorem{question}[theorem]{Question}
\newtheorem{example}[theorem]{Exemple}
\newtheorem{conjecture}[theorem]{Conjecture}
\newtheorem{exercice}[theorem]{Exercice}

\newcommand{\cal}{\mathcal}
\newcommand{\wh}{\widehat}
\newcommand{\lan}{\langle}
\newcommand{\ran}{\rangle}

\newcommand{\bi}{\begin{itemize}}
\newcommand{\ei}{\end{itemize}}
\newcommand{\be}{\begin{enumerate}}
\newcommand{\ee}{\end{enumerate}}

\newcommand{\bpf}{\begin{proof}}
\newcommand{\epf}{\end{proof}}

\newcommand{\bt}{\begin{theorem}}
\newcommand{\et}{\end{theorem}}
\newcommand{\brap}{\begin{rappel}}
\newcommand{\erap}{\end{rappel}}
\newcommand{\bnt}{\begin{notation}}
\newcommand{\ent}{\end{notation}}
\newcommand{\bd}{\begin{definition}}
\newcommand{\ed}{\end{definition}}
\newcommand{\ble}{\begin{lemma}}
\newcommand{\ele}{\end{lemma}}
\newcommand{\bpr}{\begin{proposition}}
\newcommand{\epr}{\end{proposition}}
\newcommand{\bre}{\begin{remark}}
\newcommand{\ere}{\end{remark}}
\newcommand{\bco}{\begin{corollary}}
\newcommand{\eco}{\end{corollary}}
\newcommand{\beq}{\begin{equation}}
\newcommand{\eeq}{\end{equation}}
\newcommand{\bq}{\begin{question}}
\newcommand{\eq}{\end{question}}
\newcommand{\bp}{\begin{problem}}
\newcommand{\ep}{\end{problem}}
\newcommand{\beqn}{\begin{eqnarray*}}
\newcommand{\eeqn}{\end{eqnarray*}}
\newcommand{\bex}{\begin{example}}
\newcommand{\eex}{\end{example}}
\newcommand{\ber}{\begin{exercice}}
\newcommand{\eer}{\end{exercice}}
\newcommand{\sct}{\section}
\newcommand{\ssct}{\subsection}
\newcommand{\sk}{\smallskip}
\newcommand{\bk}{\bigskip}
\newcommand{\nk}{\noindent}
\newcommand{\pl}{\partial}
\newcommand{\ov}{\overline}
\newcommand{\fr}{\frac}
\newcommand{\bcj}{\begin{conjecture}}
\newcommand{\ecj}{\end{conjecture}}

\newcommand{\bck}{\backslash}

\newcommand{\Si}{\Sigma}
\newcommand{\s}{\sigma}

\renewcommand{\t}{\theta}
\renewcommand{\a}{\alpha}
\renewcommand{\b}{\beta}
\newcommand{\g}{\gamma}

\renewcommand{\l}{\lambda}
\renewcommand{\O}{\Omega}

\newcommand{\bs}{\boldsymbol}

\title[]{ Sur les vari\'et\'es $X\subset \P^N$ telles que par $n$
points passe une courbe de $X$ de degr\'e donn\'e}

\author[L. Pirio et J.-M. Tr\'epreau]{Luc Pirio et Jean-Marie Tr\'epreau}

\address{L. Pirio, IRMAR, UMR 6625 du CNRS, Universit\'e  Rennes1, Campus de beaulieu, 
35042 Rennes Cedex.}
\email{luc.pirio@univ-rennes1.fr} 

\address{J.-M.~Tr\'epreau, 
U.P.M.C., UMR 7586 du CNRS,  bureau 26-16-5-21, 4 place Jussieu, 75005 Pa\-ris.}
\email{trepreau@math.jussieu.fr}

\begin{document}

\maketitle

\begin{abstract}
Soit $\,r\geq 1$, $\,n\geq 2$, et $\,q\geq n-1$ des entiers.
On introduit la classe $\cal{X}_{r+1,n}(q)$ des sous-vari\'et\'es $X$
de dimension $r+1$ d'un espace projectif,  telles que
\bi
\item pour $(x_1,\ldots,x_n)\in X^n$ g\'en\'erique, il existe une courbe
rationnelle normale de degr\'e $q$, contenue dans $X$ et passant par
les points $x_1,\ldots,x_n$ ;
\item $X$ engendre un espace projectif dont la dimension,
pour $r$, $n$ et $q$ donn\'es, est la plus grande possible
compte tenu de la premi\`ere propri\'et\'e.
\ei
Sous l'hypoth\`ese $q\neq 2n-3$, on d\'etermine toutes les vari\'et\'es
$X$ appartenant \`a la classe $\cal{X}_{r+1,n}(q)$. 
On montre en particulier qu'il existe une vari\'et\'e  $X_0\subset \P^{r+n-1}$
de degr\'e minimal $n-1$ et une application birationnelle 
$X_0\dasharrow X$ qui envoie une section de $X_0$ par 
un $\P^{n-1}\subset \P^{r+n-1}$ g\'en\'erique sur une 
courbe rationnelle normale de degr\'e $q$.

\sk
Sans hypoth\`ese sur $q$, on d\'efinit sur l'espace des courbes 
rationnelles normales de degr\'e $q$ contenues dans la vari\'et\'e $X\in \cal{X}_{r+1,n}(q)$
une structure quasi-grassmannienne. La 
vari\'et\'e $X$ est de la forme pr\'ec\'edente si et seulement si 
cette structure est localement isomorphe \`a la structure standard, celle 
de la grassmannienne des $(n-1)$-plans de $\P^{r+n-1}$.

\sk
Le probl\`eme de la d\'etermination des vari\'et\'es $X\in \cal{X}_{r+1,n}(2n-3)$
reste ouvert. Nous donnons quelques exemples de vari\'et\'es des classes $\cal{X}_{r+1,3}(3)$
et $\cal{X}_{r+1,4}(5)$ qui ne sont pas de la forme qu'on vient de d\'ecrire.

\end{abstract}

\sct{Introduction}
\label{I}

\ssct{Les classes $\cal{X}_{r+1,n}(q)$}

On consid\`ere un probl\`eme de g\'eom\'etrie projective complexe,
auquel nous ont conduits nos travaux r\'ecents sur 
l'alg\'ebrisation des tissus de rang maximal, voir Pirio-Tr\'epreau~\cite{PT}.

\sk
Soit $\,r\geq 1$, $\,n\geq 2$ et $\,q\geq n-1$ des entiers.
\bd
\label{D1}
On note $\cal{X}_{r+1,n}(q)$ la classe des vari\'et\'es $X$,
sous-vari\'et\'es  alg\'ebriques irr\'eductibles de dimension 
$r+1$ d'un espace projectif quelconque, telles que :
\be
\item[1)]  pour $(x_1,\ldots,x_n)\in X^n$ g\'en\'erique, il existe une courbe
rationnelle normale de degr\'e $q$, contenue dans $X$ et passant par
les points $x_1,\ldots,x_n$ ;
\item[2)] $X$ engendre un espace projectif dont la dimension, not\'ee $\pi_{r,n}(q)$,
est la plus grande possible, 
compte tenu de la premi\`ere propri\'et\'e.
\ee
\ed
La dimension $\pi_{r,n}(q)$ de l'espace engendr\'e par une vari\'et\'e de la classe 
$\cal{X}_{r+1,n}(q)$ est donn\'ee par la formule (\ref{pi}) de la 
section suivante.

Notre r\'esultat principal sera la d\'etermination 
de toutes les vari\'et\'es $X$ de la classe $\cal{X}_{r+1,n}(q)$
sous l'hypoth\`ese $q\neq 2n-3$. Le cas $q=2n-3$ restera ouvert.

\bk
L'hypoth\`ese d'irr\'eductibilit\'e est en fait une cons\'equence 
des autres hypoth\`eses, la d\'emonstration 
est laiss\'ee au lecteur. On peut aussi
remplacer la Propri\'et\'e~1) par la Propri\'et\'e 1') suivante, plus faible  :
\be
\item[1')] pour $(x_1,\ldots,x_n)\in X^n$ g\'en\'erique, il existe une courbe
irr\'eductible de degr\'e $\leq q$, contenue dans $X$ et passant par
les points $x_1,\ldots,x_n$,
\ee
sans modifier la classe $\cal{X}_{r+1,n}(q)$ que l'on d\'efinit, voir la Proposition \ref{generalisation}.

Il est tentant de penser, mais nous ne le d\'emontrons pas,
qu'on peut m\^eme la remplacer 
par la Propri\'et\'e 1'') suivante, encore plus faible :
\be
\item[1'')] pour $(x_1,\ldots,x_n)\in X^n$ g\'en\'erique, il existe une courbe
irr\'eductible contenue dans $X$, passant par
$x_1,\ldots,x_n$ et engendrant un espace de dimension $\leq q$.
\ee

\ssct{Le th\'eor\`eme de la borne}

Notre premi\`ere t\^ache est de d\'eterminer la dimension $\pi_{r,n}(q)$
de l'espace engendr\'e par une vari\'et\'e $X\in \cal{X}_{r+1,n}(q)$.
Ce sera la seule d\'emonstration de l'introduction. Elle est assez
typique de cet article.

\sk
Soit $X\subset \P^N$ une vari\'et\'e alg\'ebrique irr\'eductible de dimension $r+1$. 
On note $\lan X \ran$ le sous-espace projectif qu'elle engendre, $X_{\rm reg}$ sa partie r\'eguli\`ere,
$X_{\rm sing}$ sa partie singuli\`ere. On note ${\rm CRN}_q(X)$ 
l'ensemble des courbes rationnelles normales de degr\'e $q$ contenues dans $X$.
Si $\L$ et $\L'$ sont des sous-espaces d'un espace projectif,
on note $\L\oplus \L'$ le sous-espace qu'ils engendrent s'il 
est de dimension (maximale) $\dim \,\L + \dim \,\L' + 1$ et on dit que 
c'est {\em la somme directe projective} de $\L$ et de $\L'$.

\sk
La notion d'espace osculateur ou d'osculateur sera omnipr\'esente dans la suite.
\'Etant donn\'e $x\in X_{\rm reg}$ et $k\in \N$, on note $X_x(k)$ 
{\em l'osculateur \`a l'ordre $k$ de $X$ en $x$}. C'est 
par exemple le sous-espace projectif de $\lan X \ran\subset \P^N$ passant par $x$
et dont la direction est le sous-espace vectoriel de $T_x\P^N$
engendr\'e, dans une carte affine identifi\'ee \`a $\C^N$,
mais la d\'efinition ne d\'epend pas de cette identification, 
par les d\'eriv\'ees $v'(0),\ldots,v^{(k)}(0)$
des germes de courbes param\'etr\'ees $v: (\C,0) \rightarrow (X,x)$
telles que $v(0)=x$. C'est une notion projective.

On~a toujours 
$$
\dim \, X_x(k) + 1 \leq {r+1+k \choose r+1}.
$$
On dit que $X$ est {\em  $k$-r\'eguli\`ere} en $x\in X_{\rm reg}$ si l'in\'egalit\'e 
ci-dessus est une \'egalit\'e et que $X$ est $k$-r\'eguli\`ere 
si $X$ est $k$-r\'eguli\`ere au point g\'en\'erique de $X$.

\sk
Comme la notion d'osculateur joue un grand r\^ole dans 
cet article, on rappelle dans l'Appendice les propri\'et\'es,
toutes \'el\'ementaires, qu'on utilisera sans r\'ef\'erence,
avec au moins des esquisses de d\'emonstrations. Le lecteur 
est invit\'e \`a  consulter les quelques \'enonc\'es de cet 
appendice.

\sk
La deuxi\`eme propri\'et\'e dans la D\'efinition \ref{D1} est pr\'ecis\'ee par l'\'enonc\'e suivant :
\bt
\label{Th1}
Une vari\'et\'e $X\in \cal{X}_{r+1,n}(q)$ engendre un espace de dimension 
le nombre $\pi_{r,n}(q)$ donn\'e par 
\beq
\label{pi}
\pi_{r,n}(q) + 1 =  m{r + \rho +1 \choose r+1} + (n-1-m){r + \rho \choose r+1},
\eeq
o\`u $q = \rho(n-1) + m - 1$ est la division euclidienne de $q$ par $n-1$.
\et
On d\'emontre ici la seule partie de l'\'enonc\'e qui est nouvelle : 
si l'on d\'efinit le nombre $\pi_{r,n}(q)$ 
par (\ref{pi}), une vari\'et\'e $X\in \cal{X}_{r+1,n}(q)$ engendre un espace de dimension 
$\leq \pi_{r,n}(q)$. \`A la terminologie pr\`es, le fait qu'il existe des vari\'et\'es
projectives qui v\'erifient la premi\`ere propri\'et\'e dans la D\'efinition \ref{D1}
et qui engendrent un espace de dimension $\pi_{r,n}(q)$ est d\'ej\`a connu,
voir la fin de cette introduction et le Chapitre~5.
\bpf
Soit $X\subset \P^N$ une sous-vari\'et\'e de dimension $r+1$, v\'erifiant la 
premi\`ere propri\'et\'e dans la D\'efinition \ref{D1}.
Par d\'efinition, il existe un $n$-uplet $(a_1,\ldots,a_n)\in (X_{\rm reg})^n$ 
de points deux-\`a-deux distincts tel que, pour tout $x\in X$ voisin de $a_n$,
il existe une courbe $C(x)\in \text{CRN}_q(X)$ passant par
les points $a_1,\ldots,a_{n-1}$ et $x$. 

Consid\'erons des osculateurs
$$
X_{a_i}(\rho_i) \;\; (i=1,\ldots,n-1), \;\;\;\;  \text{avec} \;\;\;\;  \sum_{i=1}^{n-1}(\rho_i+1) = q+1.
$$
Une courbe $C(x)$ est contenue
dans l'espace engendr\'e par ses osculateurs $C(x)_{a_i}(\rho_i)$, 
donc dans celui engendr\'e par les osculateurs $X_{a_i}(\rho_i)$,
$i=1,\ldots,n-1$. Comme les courbes $C(x)$ recouvrent un voisinage de $a_n$ dans $X$,
cet espace contient aussi $X$ et donc :
\beq
\label{1.1}
\dim \, \lan X\ran + 1 \leq  \sum_{i=1}^{n-1} {r+1+\rho_i \choose r+1}.
\eeq
Il suffit de prendre $\rho_i=\rho$ pour $m$ valeurs de $i$ et $\rho_i=\rho-1$ pour
les $n-m-1$ autres valeurs de $i$. On obtient le r\'esultat.
\epf
On a l'\'egalit\'e dans (\ref{1.1}) si et seulement si les 
osculateurs $X_{a_i}(\rho_i)$, $i=1,\ldots,n-1$, sont de dimensions maximales 
et en somme directe projective. En permutant les points $a_1,\ldots,a_n$, on obtient 
le r\'esultat suivant dont on d\'emontrera une sorte de  r\'eciproque au d\'ebut du Chapitre~2.
\ble
\label{L1}
Soit $X\in \cal{X}_{r+1,n}(q)$ et $a_1,\ldots,a_n$ des points deux-\`a-deux distincts 
de $X_{\rm reg}$ tels que, pour tout $(x_1,\ldots,x_n)\in X^n$ voisin de $(a_1,\ldots,a_n)$,
il existe une courbe $C\in \text{\rm CRN}_q(X)$ passant par
les points $x_1,\ldots,x_n$. Si $q=\rho(n-1)+m-1$ est la division euclidienne de $q$ par $n-1$,
$X$ est $\rho$-r\'eguli\`ere 
en chaque point $a_i$ et pour toute permutation $\s$ de $\{1,\ldots,n\}$, on a :
\beq
\label{sommedirect}
\lan X\ran  = 
(\oplus_{i=1}^m X_{a_{\s(i)}}(\rho))\oplus (\oplus_{i=m+1}^{n-1}X_{a_{\s(i)}}(\rho-1)).
\eeq
\ele

\ssct{Exemples : les classes $\cal{X}_{r+1,n}(n-1)$ et $\cal{X}_{r+1,2}(q)$} 
\label{I;2}    

Consid\'erons d'abord le cas particulier d'une vari\'et\'e $X\in \cal{X}_{r+1,n}(n-1)$.
La formule (\ref{pi}) donne $\pi_{r,n}(n-1) = r+n-1$.

\sk
La vari\'et\'e $X$, de dimension $r+1$, engendre un espace $\P^{r+n-1}$ 
et v\'erifie que, pour $(x_1,\ldots,x_n)\in X^n$ g\'en\'erique, il existe une 
courbe $C\in \text{CRN}_{n-1}(X)$ qui passe par $x_1,\ldots,x_n$. Il est bien connu 
que, pour $(x_1,\ldots,x_n)\in X^n$ g\'en\'erique, les $n$ points $x_1,\ldots,x_n$
engendrent un $\P^{n-1}$ et que $X\cap \P^{n-1}$
est une courbe irr\'eductible (et r\'eduite). C'est donc la courbe $C$. 
On en d\'eduit que $X$ est de degr\'e $n-1$, le degr\'e minimal d'une vari\'et\'e projective de dimension 
$r+1$ qui engendre un espace de dimension $r+n-1$. Une telle vari\'et\'e est 
appel\'ee une {\em vari\'et\'e minimale}. La r\'eciproque est \'evidente :
\ble
\label{L2}
Les vari\'et\'es appartenant \`a la classe $\cal{X}_{r+1,n}(n-1)$ sont les vari\'et\'es minimales
de dimension $r+1$ et de degr\'e $n-1$.
\ele 
Ces vari\'et\'es sont connues. Nous en rappellerons la classification  dans le Chapitre~5.

\bk
Une vari\'et\'e de Veronese de dimension $s\geq 1$ et d'ordre $q$ est une vari\'et\'e 
projective 
qui, dans l'espace qu'elle engendre, est l'image de $\P^s$ par un plongement
associ\'e au syst\`eme 
lin\'eaire $\left| \cal{O}_{\P^s}(q) \right|$. Par exemple, une vari\'et\'e
de Veronese de dimension $1$ et d'ordre $q$ est une courbe rationnelle 
normale de  degr\'e $q$.

\sk
Bompiani a \'etudi\'e dans \cite{Bo} un probl\`eme de g\'eom\'etrie diff\'erentielle projective
dont il a r\'eduit la solution, et cette r\'eduction occupe presque tout l'article, 
\`a la d\'emonstration d'un cas particulier de son r\'esultat principal qu'on
peut \'enoncer ainsi : 
\bt
\label{Th2}
Les vari\'et\'es appartenant \`a la classe $\cal{X}_{r+1,2}(q)$
sont les vari\'et\'es de Veronese de dimension $r+1$ et d'ordre $q$.
\et
La d\'emonstration de Bompiani, par r\'ecurrence sur $r$, est simple et synth\'etique
mais gu\`ere convaincante sauf pour $r=1$.
Dans le doute, le second auteur \cite{Tr} a d\'emontr\'e, essentiellement par un 
calcul formel, le r\'esultat plus fort suivant :
\bt
\label{TR1}
Soit $X\subset  \P^N$ un germe de vari\'et\'e lisse et $q$-r\'egulier 
en $x^\star\in \P^N$. Si pour tout $x\in X$ voisin de $x^\star$,
il existe une courbe rationnelle
normale de degr\'e $q$ localement contenue
dans $X$, qui passe par $x^\star$ et $x$, alors $X$ est une vari\'et\'e de Veronese d'ordre $q$.
\et
Pour obtenir la m\^eme conclusion, Bompiani 
suppose que l'hypoth\`ese de l'\'enonc\'e pr\'ec\'edent
reste v\'erifi\'ee si l'on remplace le point de base $x^\star$ par un point voisin.

Dans \cite{Tr} on trouvera aussi une caract\'erisation des germes 
de vari\'et\'es lisses et $q$-r\'eguliers en $x^\star\in \P^N$
tels que, pour tout $x\in X$ voisin de $x^\star$,
il existe une courbe lisse, localement contenue dans 
$X$ et passant par $x^\star$ et $x$, qui engendre un $q$-plan.
Selon l'avis de son auteur, ce r\'esultat rend plausible que, dans la D\'efinition \ref{D1},
on pourrait remplacer la premi\`ere propri\'et\'e par la Propri\'et\'e 1'')
sans modifier la classe $\cal{X}_{r+1,n}(q)$ que l'on d\'efinit.

\ssct{\'Enonc\'es des principaux r\'esulats}
\label{I;3}

Notre r\'esultat principal est le suivant, voir aussi le Th\'eor\`eme \ref{final} 
du Chapitre 3 pour un \'enonc\'e plus pr\'ecis.  
\bt
\label{Th4}
Soit $X\in \cal{X}_{r+1,n}(q)$. Si $q\neq 2n-3$, ou si $q=2n-3$ et $r=1$ ou $n=2$,
il existe une vari\'et\'e $X_0\subset \P^{r+n-1}$ de degr\'e minimal $n-1$,
et une application birationnelle $\phi : \; X_0\dasharrow X$ telle que 
l'image d'une section de $X_0$ par un $\P^{n-1}\subset \P^{r+n-1}$
g\'en\'erique 
appartienne \`a $\text{\rm CRN}_q(X)$.
\et
 
Pour autant que nous le sachions, ce r\'esultat est nouveau sauf dans le cas $q=n-1$, 
dans le cas $n=2$ et dans le cas $r=1$. 

\sk
Dans le cas $q=n-1$, on retrouve le Lemme \ref{L2},
qui est standard. 
Dans le cas $n=2$, on retrouve 
le th\'eor\`eme de Bompiani. 
Nous ne le red\'emontrons pas ici.
Il joue un r\^ole tr\`es important dans la d\'emonstration 
du r\'esultat g\'en\'eral.

\sk
Le cas $r=1$ des surfaces est particulier. Il est beaucoup
plus simple pour la raison que, dans ce cas, une courbe 
de $X$ est un diviseur. La th\'eorie des syst\`emes lin\'eaires
permet alors de d\'emontrer sans trop de difficult\'es le Th\'eor\`eme \ref{Th4}.
M\^eme si nous n'avons pas trouv\'e l'\'enonc\'e
correspondant dans la litt\'erature, il est bien possible qu'il soit 
folklorique.

\sk
La d\'emonstration du Th\'eor\`eme \ref{Th4} occupe le Chapitre 2 et le Chapitre 3.
Le lecteur d\'esireux d'en  avoir une id\'ee peut survoler 
les introductions \`a ces chapitres.

\sk
Le Th\'eor\`eme \ref{Th4} sugg\`ere la d\'efinition suivante :
\bd
\label{D2}
La vari\'et\'e $X\in \cal{X}_{r+1,n}(q)$ est {\em standard},
sous-entendu comme \'el\'e\-ment de $\cal{X}_{r+1,n}(q)$,
s'il existe une vari\'et\'e 
minimale $X_0\in \cal{X}_{r+1,n}(n-1)$ et une application 
birationnelle $\phi : \; X_0\dasharrow X$ telle que l'image 
d'une section de $X_0$ par un $\P^{n-1}\subset \P^{r+n-1}$ 
g\'en\'erique appartienne \`a $\text{\rm CRN}_q(X)$.
On dit alors que {\em $X$ est associ\'ee \`a $X_0$}
par l'application $\phi$.
Elle est {\em sp\'eciale} dans le cas contraire.
\ed
Avec cette terminologie, on peut r\'e\'ecrire le Th\'eor\`eme \ref{Th4} 
sous la forme suivante.

\bk
{\em Si $q\neq 2n-3$, ou si $r=1$ ou $n=2$, toute vari\'et\'e de 
la classe $\cal{X}_{r+1,n}(q)$ est standard.}

\bk
Une m\^eme vari\'et\'e $X$ peut appartenir \`a plusieurs
classes $\cal{X}_{r+1,n}(q)$ distinctes. 
L'exemple de la vari\'et\'e de Veronese 
de dimension $3$ et d'ordre $3$ est curieux : cette vari\'et\'e
est un \'el\'ement standard de la classe $\cal{X}_{3,2}(3)$
et un \'el\'ement sp\'ecial de la classe $\cal{X}_{3,6}(9)$,
comme on le verra dans le Chapitre 6.

\bk
Si $X\in \cal{X}_{r+1,n}(q)$, nous montrerons dans le Chapitre 4
qu'un ouvert $\Si_q(X)$  de l'ensemble ${\rm CRN}_q(X)$, 
 qu'on d\'efinira pr\'ecis\'ement, admet une ${\rm G}_{r,n}$-structure naturelle.
On appelle ainsi le type de $G$-structure 
model\'ee sur la structure infinit\'esimale de 
la grassmannienne $\G_{r,n}$ des $\P^{n-1}$ de $\P^{r+n-1}$.
Si par exemple $X_0$  est une vari\'et\'e minimale de degr\'e $n-1$,
c'est la structure induite par l'application qui,
\`a un $\P^{n-1}\in \G_{r,n}$ g\'en\'erique, associe 
la courbe $X_0\cap \P^{n-1}\in \Si_{n-1}(X_0)$. Nous d\'emontrerons  : 
\bt
\label{Th6}
La vari\'et\'e $X\in \cal{X}_{r+1,n}(q)$ est standard 
si et seulement si la  ${\rm G}_{r,n}$-structure naturelle
de $\Si_q(X)$ est int\'egrable.
\et

Le Chapitre 5 peut \^etre lu ind\'ependamment du reste de l'article. 
Il contient la classification des 
vari\'et\'es standards, qu'on obtiendra \`a partir de celle des vari\'et\'es minimales
et de la structure de leurs groupes de Picard. Cette classification est donn\'ee 
par le Th\'eor\`eme~\ref{class}, dont l'\'enonc\'e est trop long pour qu'on le
reproduise ici.

\sk
Le premier exemple que nous avons connu d'une vari\'et\'e sp\'eciale nous a \'et\'e 
communiqu\'e par F. Russo \cite{Ru}. Il s'agit de l'image de $\P^1\times \P^1\times \P^1$
par le plongement de Segre de type $(1,1,1)$. C'est une vari\'et\'e sp\'eciale
de la classe $\cal{X}_{3,3}(3)$.

Dans le Chapitre 5, nous pr\'esenterons cet exemple ainsi que quelques autres,
obtenant en particulier le r\'esultat suivant. 
\bt
\label{Th5}
Pour tout $r\geq 2$, les classes $\cal{X}_{r+1,3}(3)$ et $\cal{X}_{r+1,4}(5)$
contiennent des vari\'et\'es sp\'eciales.
\et
Cet \'enonc\'e ne fait qu'effleurer le probl\`eme tr\`es int\'eressant 
de la classification compl\`ete des vari\'et\'es des classes $\cal{X}_{r+1,n}(2n-3)$ avec $r\geq 2$
et $n\geq 3$.
Nous connaissons quelques exemples encore de vari\'et\'es sp\'eciales en plus de ceux que
nous mentionnerons mais, \`a l'heure qu'il est, nous savons tr\`es peu de choses.
Nous ne  savons m\^eme pas s'il existe des vari\'et\'es sp\'eciales  
dans toutes ces classes.

\ssct{Une relation avec les vari\'et\'es de genre maximal}
\label{I;4}

Pour finir, nous  pr\'e\-sentons sans d\'emonstration la relation int\'eressante 
qui existe entre les classes $\cal{X}_{r+1,n}(q)$ et les vari\'et\'es de genre g\'eom\'etrique maximal.

\sk
Soit $Z\subset \P^{r+n-1}$ une vari\'et\'e projective de dimension $r$ et de degr\'e $d$, qui engendre $\P^{r+n-1}$.
Si $Z$ est irr\'eductible et lisse, le {\em genre g\'eom\'etrique} de $Z$ est la dimension 
$g(Z)$ de l'espace $H^0(Z,\O^r_Z)$ des $r$-formes holomorphes sur $Z$. 
 \'Ecrivons 
$$
d - 1 = \s (n-1) + m, \qquad m\in \{1,\ldots,n-1\}.
$$
Pour une dimension $r$ et un degr\'e $d$ donn\'es, le genre g\'eom\'etrique de $Z$
est au plus \'egal, et le r\'esultat est optimal, au nombre 
\beq
\label{castel1}
g_{r,n}(d) = m \, {\s + 1 \choose r+1 } + (n-1-m) { \s \choose r+1 }.
\eeq
C'est {\em la borne de Castelnuovo-Harris}. Deux d\'emonstrations au moins de ce r\'esultat 
sont connues. La premi\`ere est due \`a Chern et Griffiths \cite{CG} et repose 
sur le th\'eor\`eme d'addition d'Abel et un argument combinatoire inspir\'e de la 
th\'eorie des tissus. La deuxi\`eme est due \`a Harris~\cite{Ha} et appartient
\`a la g\'eom\'etrie alg\'ebrique. La formule (\ref{castel1}) est 
une \'ecriture l\'eg\`erement modifi\'ee de la formule de \cite{Ha}, page 65.

\sk
En fait, pour le probl\`eme qui nous int\'eresse, c'est trop demander que 
$Z$ soit lisse ou m\^eme que $Z$ soit irr\'eductible. 

\sk
Dans \cite{Gr}, en supposant seulement que $Z$ est r\'eduit,
Griffiths associe \`a toute $r$-forme rationnelle $\omega$ sur $Z$
sa trace $\text{Tr}\, (\omega)$, une $r$-forme rationnelle
sur la grassmannienne $\G_{r,n}$.
Il introduit l'espace des $r$-formes rationnelle sur $Z$
de trace nulle. (Si $Z$ est lisse, ce sont les formes holomorphes
sur $Z$.) Faute d'une terminologie traditionnelle, appelons 
{\em genre g\'eom\'etrique corrig\'e} de $Z$ la dimension
$g(Z)$ de cet espace. On suppose de plus qu'un $\P^{n-1}\subset \P^{r+n-1}$ g\'en\'erique coupe $Z$ en $d$
points en position g\'en\'erale dans ce $\P^{n-1}$.
Sous cette hypoth\`ese,
on a encore l'in\'egalit\'e $g(Z)\leq g_{r,n}(d)$.

\sk
Supposons $d \geq (r+1)(n-1)+2$ et posons $q=d- r(n-1)- 2$. On a 
donc $q \geq n-1$ et $q = (\s - r)(n-1) + m - 1$ est la division euclidienne 
$q = \rho(n-1)+m-1$ de $q$ par $(n-1)$. La comparaison 
des formules (\ref{castel1}) et (\ref{pi}) donne $g_{r,n}(d) = \pi_{r,n}(q) + 1$.

\bk
{\em Le nombre $\pi_{r,n}(q)+1$ est le genre g\'eom\'etrique corrig\'e 
maximal d'une vari\'et\'e de dimension~$r$ et de degr\'e $d=q+r(n-1)+2$,
qui engendre un  espace de dimension $r+n-1$.}

\bk
Rappelons maintenant  quelques r\'esultats de Harris \cite{Ha} dans le cas 
o\`u $Z$ est irr\'eductible. Si $g(Z)=g_{r,n}(d)$,
la vari\'et\'e $Z$
est contenue dans une vari\'et\'e minimale $X_0\in \cal{X}_{r+1,n}(n-1)$.
L'espace des $r$-formes rationnelles sur $Z$ de trace 
nulle d\'efinit une {\em application canonique} $Z \dasharrow \P^{\pi_{r,n}(q)}$
qui, c'est une propri\'et\'e importante, se prolonge \og canoniquement \fg\,
en une application rationnelle $c_Z : X_0\dasharrow \P^{\pi_{r,n}(q)}$,
dont l'image est de dimension $r+1$.
(Celle-ci ne d\'epend que de la classe du
diviseur $Z$ dans le groupe de Picard de $X_0$.)
L'image d'une section de $X_0$ par un $\P^{n-1}\subset \P^{r+n-1}$
g\'en\'erique est une courbe rationnelle normale de degr\'e $q$.
Il en r\'esulte que l'image $c_Z(X_0)$ de $X_0$ est une 
vari\'et\'e de la classe $\cal{X}_{r+1,n}(q)$.

\sk
Nous avons :
\bt
\label{Th7}
Toute vari\'et\'e standard $X\in \cal{X}_{r+1,n}(q)$ associ\'ee,
au sens de la D\'efinition \ref{D2},  
\`a une vari\'et\'e minimale
$X_0\in \cal{X}_{r+1,n}(n-1)$, est l'image de $X_0$
par l'application $c_Z: \, X_0\dasharrow \P^{\pi_{r,n}(q)}$ 
associ\'ee \`a une vari\'et\'e $Z\subset X_0$, de dimension $r$,
de degr\'e  $d=q+r(n-1)+2$ et de genre g\'eom\'etrique corrig\'e maximal, i.e. \'egal \`a $g_{r,n}(d)$.
\et
C'est essentiellement une cons\'equence du Th\'eor\`eme \ref{Th1}, \`a ceci pr\`es
qu'il faut v\'erifier la compatibilit\'e entre les points de vue de \cite{Ha}
et de \cite{Gr}. Comme cette v\'erification fait appel \`a des notions 
qui ne sont pas introduites dans cet article, nous la reportons \`a \cite{PT}.

\ssct{Remarques sur le style}

Le style de cet article est relativement \'el\'ementaire. 
Nous n'utilisons pas le langage de la g\'eom\'etrie alg\'ebrique 
moderne.  Une certaine familiarit\'e 
avec la g\'eom\'etrie projective complexe, par exemple 
telle qu'elle est expos\'ee dans le \og cours pr\'eparatoire \fg\, 
de Mumford \cite{Mu} est suffisante pour sa lecture.
Le r\'edacteur a souvent pr\'ef\'er\'e donner une d\'emonstration,
surtout si elle est courte, plut\^ot qu'une r\'ef\'erence \`a la litt\'erature.

Nous esp\'erons que nos r\'esultats int\'eresseront quelques
g\'eom\`etres et qu'ils ne seront pas rebut\'es par le style, qu'ils
trouveront peut-\^etre d\'emod\'e. L'expert pourra toujours sauter les passages qu'il juge
exag\'er\'ement pointilleux.

Enfin, certains arguments sont inspir\'es par 
des th\'eor\`emes classiques, par exemple le th\'eor\`eme de d\'eformation de Kodaira~\cite{Ko}
ou la caract\'erisation des syst\`emes lin\'eaires due 
\`a Enriques~\cite{En}. Compte tenu du cadre tr\`es particulier 
de notre \'etude, il s'est toutefois av\'er\'e plus simple de s'en passer que 
de chercher \`a savoir si d'\'eventuelles  
versions de ces r\'esultats pouvaient \^etre appliqu\'ees.

\sct{Projections sur des vari\'et\'es de Veronese et applications}

\ssct{Introduction}

Dans tout ce chapitre, $\,r\geq 1$, $\,n\geq 2$, et $\,q\geq n-1$ sont
des entiers fix\'es. On  note 
$$
q = \rho(n-1) + m - 1
$$
la division euclidienne de $q$ par $n-1$.

On note  $N$ l'entier $\pi_{r,n}(q)$ d\'efini
par la formule~(\ref{pi}) et on se donne une vari\'et\'e $X\subset \P^N$ de la classe $\cal{X}_{r+1,n}(q)$.

\bk
On va obtenir des propri\'et\'es importantes de 
$X$ gr\^ace aux projections 
$\P^N\dasharrow X_{a_i}(\rho_i)$, $\, \rho_i\in \{\rho,\rho-1\}$,
 associ\'ees aux d\'ecompositions 
$\P^N = \oplus_{i=1}^{n-1} X_{a_i}(\rho_i)$ de $\P^N$ 
en somme directe d'osculateurs qui apparaissent dans le Lemme \ref{L1}.

\sk
On montrera en effet que l'image de la vari\'et\'e $X$ par 
une projection de ce type est une vari\'et\'e
de la classe $\cal{X}_{r+1,2}(\rho_i)$. Le point crucial 
est alors que ces vari\'et\'es sont connues. D'apr\`es le 
Th\'eor\`eme \ref{Th2}, ce sont des vari\'et\'es de Veronese
d'ordre $\rho_i$. 

\sk
Gr\^ace \`a cela, on sera en mesure de montrer que $X$
est lisse au voisinage de toute courbe $C$ \'el\'ement d'un 
ouvert dense $\Si_q(X)$, qu'on d\'efinira pr\'ecis\'ement,
de la seule composante irr\'eductible int\'eressante 
de ${\rm CRN}_q(X)$. Ce r\'esultat 
sera utile dans le Chapitre 3. 

On montrera aussi
que, pour toute courbe $C\in \Si_q(X)$, le fibr\'e
normal $N_CX$ est une somme directe de $r$ fibr\'es en droites 
de degr\'e $n-1$. C'est cette propri\'et\'e qui 
fait que $\Si_q(X)$ admet la ${\rm G}_{r,n}$-structure
naturelle qu'on consid\'erera dans le Chapitre~4. 

Les projections sur des vari\'et\'es de Veronese 
joueront encore un r\^ole majeur dans le Chapitre~3
pour construire, dans les cas favorables, le syst\`eme lin\'eaire 
qui permet d'associer $X$ \`a une vari\'et\'e minimale de 
degr\'e $n-1$, au sens de la D\'efinition \ref{D2}.

\ssct{Une r\'eciproque au Lemme \ref{L1}}

On note $\text{CR}_q(\P^N)$ la composante irr\'educ\-tible 
de la vari\'et\'e de Chow des $1$-cycles effectifs 
de degr\'e $q$ de $\P^N$ qui contient, comme ouvert dense,
l'ensemble ${\rm CRN}_q(\P^N)$ des courbes rationnelles
normales de degr\'e $q$.

\sk
Un \'el\'ement $C$ de ${\rm CR}_q(\P^N)$ s'\'ecrit  
$$
C = m_1C_1 + \cdots + m_dC_d, \qquad m_1,\ldots,m_d\in \N^*,
$$
o\`u $C_i$, $\,i=1,\ldots,d$,  est une courbe irr\'eductible
de degr\'e $q_i$ et $m_1q_1+\cdots+m_dq_d=q$.
De plus, {\em le support $|C|=\cup_{i=1}^dC_i$ du cycle $C$ est connexe}.
D'autre part, les composantes $C_i$ sont en fait des courbes rationnelles, mais 
nous n'utiliserons pas cette propri\'et\'e.

\sk
Si $X$ est une sous-vari\'et\'e alg\'ebrique irr\'eductible de $\P^N$,
on note ${\rm CR}_q(X)$ la sous-vari\'et\'e alg\'ebrique  des 
\'el\'ements de ${\rm CR}_q(\P^N)$ dont le support est contenu dans $X$.
La vari\'et\'e d'incidence 
$$
I = \{(C, {\boldsymbol x}) \in \text{CR}_q(X)\times X^n, \; {\boldsymbol x}\subset |C|\}
$$
est une vari\'et\'e projective compacte. 
Si $X\in \cal{X}_{r+1,n}(q)$,
l'image de la projection canonique $I\rightarrow X^n$ contient un ouvert non vide de $X^n$,
donc est \'egale \`a $X^n$. Le m\^eme argument montre que si $X$ engendre un espace de dimension 
$\pi_{r,n}(q)$ et si l'on suppose seulement que la premi\`ere propri\'et\'e de la D\'efinition \ref{D1}
est v\'erifi\'ee au voisinage d'{\em un} $n$-uplet ${\boldsymbol a}\in X^n$, alors 
la vari\'et\'e $X$ appartient \`a la classe $\cal{X}_{r+1,n}(q)$. Nous utiliserons 
cette remarque \`a plusieurs reprises.

\bk
On a la r\'eciproque suivante au Lemme \ref{L1} :
\ble
\label{rec-L1}
Soit $X\in \cal{X}_{r+1,n}(q)$ et $a_1,\ldots,a_n$ des points deux-\`a-deux distincts 
de $X_{\rm reg}$ tels que, pour toute permutation $\s$ de $\{1,\ldots,n\}$, on ait :
\beq
\label{sommedirect}
\lan X\ran  = 
(\oplus_{i=1}^{m} X_{a_{\s(i)}}(\rho_i)) \oplus (\oplus_{i=m+1}^{n-1}  X_{a_{\s(i)}}(\rho-1)).
\eeq
Toute courbe alg\'ebrique connexe $C\subset X$ 
de degr\'e $\leq q$ qui passe par les points 
$a_1,\ldots,a_n$ est une courbe rationnelle normale de degr\'e $q$.
\ele
Il en existe au moins une d'apr\`es les remarques pr\'ec\'edentes.
\bpf
Soit d'abord $p\leq n-1$ et $C'\subset X$ une courbe alg\'ebrique {\em connexe} passant par au moins $p$
points parmi $a_1,\ldots,a_n$, par exemple par $a_1,\ldots,a_p$.
Le degr\'e de la courbe $C'$ est au moins \'egal \`a la dimension de l'espace qu'elle engendre.
D'autre part, il existe un $p$-uplet $(b_1,\ldots,b_p)$,
aussi voisin qu'on veut de $(a_1,\ldots,a_p)$, compos\'e de points de $C'_{\rm reg}$.
Selon que $p$ est $\leq m$ ou est $>m$, on peut supposer que les osculateurs 
$$
X_{b_1}(\rho), \ldots, X_{b_p}(\rho) 
\, ; \;\; \text{ou} \;\; 
X_{b_1}(\rho), \ldots, X_{b_m}(\rho), X_{b_{m+1}}(\rho-1),\ldots,X_{b_p}(\rho-1),
$$
sont en somme directe projective. Dans chaque cas, la courbe $C'$ est lisse et donc $\rho_i$-r\'eguli\`ere 
en $b_i$, $\,i=1,\ldots,p$, et les osculateurs correspondants $C'_{b_i}(\rho_i)\subset X_{b_i}(\rho_i)$
sont en somme directe projective. On a donc, suivant les cas  :
$$
\dim\, \lan C'\ran + 1 \geq  p(\rho+1)
\, ; \;\;  \text{ou} \; \dim\, \lan C'\ran + 1 \geq  m(\rho+1)+ (p-m)\rho.
$$
En r\'esum\'e, {\em  une courbe alg\'ebrique connexe qui passe par $p\leq n-1$ points 
parmi $a_1,\ldots,a_n$ engendre un espace de dimension $\geq \rho p + \min(p,m)-1$.}

\sk
Soit maintenant $C$ une courbe {\em connexe} de degr\'e $\leq q$ qui 
passe par $a_1,\ldots,a_{n-1}$. Elle engendre un espace de dimension
$\geq q$
d'apr\`es ce qui pr\'ec\`ede. Si elle 
est irr\'eductible, c'est donc une courbe rationnelle normale de 
degr\'e $q$. Il reste \`a montrer qu'elle est irr\'eductible. 

\sk
Si ce n'est pas le cas, on peut \'ecrire 
$C= C' \cup C''$ o\`u $C'$ et $C''$ sont des courbes connexes 
de degr\'e $\geq 1$ et $q\geq \deg C = \deg C' + \deg C''$.
Si $C'$ passe par tous les points $a_1,\ldots,a_n$, o\`u m\^eme 
par $(n-1)$ d'entre eux, elle 
est de degr\'e $\geq q$, ce qui est impossible.

On peut donc supposer que $C'$ passe par exactement $p\in \{1,\ldots,n-1\}$
points parmi $a_1,\ldots,a_n$ et $C''$ par les autres. On obtient
$q\geq \rho n + \min \,(p,m) + \min \,(n-p,m) - 2$.
Comme $q=\rho(n-1)+m-1$ et $\rho\geq 1$, il vient $m \geq \min \,(p,m) + \min \,(n-p,m)$
et enfin, puisque  $m$, $p$ et $n-p$ sont strictement positifs, on obtient 
$m \geq n$, une contradiction.
\epf

\sk
On a mentionn\'e le r\'esultat suivant dans la Section 1.1. S'il est int\'eressant en tant que tel,
il ne servira pas dans la suite.
\bco
\label{generalisation}
Soit $X$ une sous-vari\'et\'e  alg\'ebrique irr\'eductible de dimension 
$r+1$ d'un espace projectif telle que, 
pour $(x_1,\ldots,x_n)\in X^n$ g\'en\'erique, il existe une courbe
irr\'eductible de degr\'e $\leq q$, contenue dans $X$ et passant par
les points $x_1,\ldots,x_n$. La vari\'et\'e $X$ engendre un espace de dimension $\leq \pi_{r,n}(q)$.
Si elle  engendre un espace de dimension $\pi_{r,n}(q)$, elle 
appartient \`a la classe $\cal{X}_{r+1,n}(q)$.
\eco
\bpf
Soit $q=\rho(n-1)+m-1$ la division euclidienne de $q$ par $n-1$.

Dans la d\'emonstration pr\'ec\'edente, on n'a pas utilis\'e 
le fait que $X$ appartient \`a la classe $\cal{X}_{r+1,n}(q)$
mais seulement le fait que $X$ engendre un espace de dimension 
$\pi_{r,n}(q)$ et v\'erifie la propri\'et\'e (\ref{sommedirect}).
Il suffit donc, pour d\'emontrer le corollaire, de montrer 
que $X$ engendre un espace de dimension $\leq \pi_{r,n}(q)$ et que,
si $X$  engendre un espace de dimension $\pi_{r,n}(q)$,
la propri\'et\'e (\ref{sommedirect}) est v\'erifi\'ee pour 
$(a_1,\ldots,a_n)\in X^n$ g\'en\'erique.

\sk
On fait l'hypoth\`ese, plus faible que celle de l'\'enonc\'e, que pour $(x_1,\ldots,x_n)\in X^n$ g\'en\'erique, 
il existe une courbe
irr\'eductible, engendrant un espace de dimension $\leq q$, contenue dans $X$ et passant par
$x_1,\ldots,x_n$. On adapte la d\'emonstration du Th\'eor\`eme~\ref{Th1}.

\sk
Soit ${\boldsymbol a}=(a_1,\ldots,a_n)$ un $n$-uplet de points deux-\`a-deux distincts de $X_{\rm reg}$
tels que l'hypoth\`ese pr\'ec\'edente soit v\'erifi\'ee pour tout ${\boldsymbol x}\in X^n$ voisin de $\boldsymbol a$.
On peut imposer au $n$-uplet $\boldsymbol a$ de v\'erifier la propri\'et\'e (g\'en\'erique)
suivante : pour toute permutation $\s$ de $\{1,\ldots,n\}$, la dimension de l'espace engendr\'e par les osculateurs 
$$
X_{a_{\s(1)}}(\rho), \ldots, X_{a_{\s(m)}}(\rho),  X_{a_{\s{(m+1)}}}(\rho-1), \ldots X_{a_{\s(n-1)}}(\rho-1),
$$
est maximale, parmi celles qu'on peut obtenir quand $\boldsymbol a$ d\'ecrit $X_{\rm reg}^n$.

\sk
Il suffit maintenant de consid\'erer le cas o\`u la permutation $\s$ est l'identit\'e.
Notons $X_{a_i}(\rho_i)$ les osculateurs qui interviennent, $i=1,\ldots,n-1$.
Les courbes $C$, dont on suppose l'existence, qui passent par $a_1,\ldots,a_{n-1}$ et un point $x$
voisin de $a_n$ recouvrent un voisinage de $a_n$ dans $X$ donc 
engendrent $\lan X\ran$. D'autre part, une telle courbe $C$ \'etant fix\'ee,
l'espace de dimension $\leq q$ qu'elle engendre est aussi engendr\'e par ses osculateurs $C_{b_i}(\rho_i)$, 
pour ${\boldsymbol b}=(b_1,\ldots,b_{n-1})\in C_{\rm reg}^{n-1}$ g\'en\'erique.
Elle est donc contenue dans l'espace engendr\'e par les osculateurs $X_{b_i}(\rho_i)$.
En faisant tendre $\boldsymbol b$ vers $\boldsymbol a$ et compte tenu des conditions impos\'ees
au $n$-uplet $\boldsymbol a$, on obtient qu'elle est contenue dans l'espace 
engendr\'e par les osculateurs $X_{a_i}(\rho_i)$. On conclut comme \`a la fin de la d\'emonstration
du Th\'eor\`eme \ref{Th1}.
\epf

\ssct{Objets admissibles ; notations}

Plut\^ot que d'abuser du mot {\em g\'en\'erique}, on introduit une terminologie 
et des notations qui seront d'un usage courant dans la suite.

\sk
{\em Quand on parle d'un {\em $p$-uplet}
${\boldsymbol a}=(a_1,\ldots,a_p)$ de points,
il est sous-entendu que ces
points sont deux-\`a-deux distincts\footnote{Les $p$-uplets avec $p \geq 2$ seront toujours not\'es par des minuscules latines
grasses ${\boldsymbol a}, \,\wh{\boldsymbol a}, \,{\boldsymbol  x}, \,\wh{\boldsymbol x} \ldots$.}.} 

Abusivement, on confondra souvent un $p$-uplet $\boldsymbol a$ et l'ensemble $\{a_1,\ldots,a_p\}$ sous-jacent
en s'autorisant des notations telles que ${\boldsymbol a}\subset C$ ou ${\boldsymbol b} \subset {\boldsymbol a}$...
Si $\boldsymbol a$ est un  $p$-uplet et $\boldsymbol b$ un $q$-uplet disjoint de $\boldsymbol a$, 
$({\boldsymbol a},{\boldsymbol b})$ d\'esigne le $(p+q)$-uplet obtenu par juxtaposition.
\bd
\label{D3.1}
Un $n$-uplet ${\boldsymbol a}=(a_1,\ldots,a_n)\in X^n$ est {\em admissible} si ${\boldsymbol a}\in X_\text{reg}^n$
et si, pour toute permutation $\s$ de $\{1,\ldots,n\}$, on a :
$$
\P^N  = 
(\oplus_{i=1}^m X_{a_{\s(i)}}(\rho))\oplus (\oplus_{i=m+1}^{n-1}  X_{a_{\s(i)}}(\rho-1)).
$$
De fa\c{c}on \'equivalente, $a_1,\ldots,a_n$ sont des points deux-\`a-deux distincts de $X_{\rm reg}$ et, pour toute permutation
$\s$ de $\{1,\ldots,n\}$, les courbes $C\in {\rm CRN}_q(X)$ qui passent par les points 
$a_{\s(1)},\ldots,a_{\s(n-1)}$ recouvrent un voisinage de $a_{\s(n)}$.
\ed
Ce sont en fait des propri\'et\'es  de l'ensemble sous-jacent \`a ${\boldsymbol a}$.
Le fait qu'elles sont \'equivalentes r\'esulte de la d\'emonstration du Th\'eor\`eme \ref{Th1}
 et du Lemme \ref{rec-L1}.

\sk
On \'etend la notion d'\^etre admissible \`a d'autres objets :
\bd
\label{D3.2}
Si $1\leq p\leq n-1$, un $p$-uplet ${\boldsymbol a}\in X^p$ est {\em admissible} si 
on peut le compl\'eter en un $n$-uplet admissible $({\boldsymbol a},{\boldsymbol b})$. 
Une courbe $C$ de $X$ est {\em admissible} si $C$ appartient 
\`a ${\rm CRN}_q(X)$ et contient un $n$-uplet admissible.
\ed
On utilisera syst\'ematiquement les notations suivantes.

\bk
\nk
\; -- \;  $\,X^{(p)}_{\rm adm}$ est l'ouvert dense de $X^p$
des $p$-uplets admissibles de $X$. 
On note $X_{\rm adm}$ au lieu de
$X^{(1)}_{\rm adm}$. C'est l'ensemble des {\em points admissibles} de $X$.

\bk
\nk
\; -- \;  Si ${\boldsymbol a}\in X^{(p)}_{\rm adm}$ et $\,p+p'\leq n$, on note  
$X^{(p')}_{\rm adm}({\boldsymbol a})$ l'ouvert dense de $X^{p'}$ d\'efini par 
la formule $X^{(p')}_{\rm adm}({\boldsymbol a})  = \{{\boldsymbol b}\in X^{p'}, \;\; ({\boldsymbol a},{\boldsymbol b}) \;\; \text{est admissible}\}$.
On note $X_{\rm adm}({\boldsymbol a})$ au lieu de $X^{(1)}_{\rm adm}({\boldsymbol a})$.

\bk
\nk
\; -- \; $\Si_q(X)$ est l'ensemble des courbes admissibles de $X$.

\bk
\nk
\; -- \; $\Si_q(X;{\boldsymbol a})$ est l'ensemble des courbes $C\in \Si_q(X)$
qui contiennent le $p$-uplet $\boldsymbol a$.

\bk
Le fait que $X^{(p')}_{\rm adm}({\boldsymbol a})$ est un ouvert dense de $X^{p'}$
se v\'erifie sans difficult\'e.

\bk
Consid\'erons l'exemple fondamental d'une vari\'et\'e  $X\in \cal{X}_{r+1,2}(q)$,
autrement dit d'une vari\'et\'e de Veronese d'ordre $q$.
Dans ce cas, $X_{\rm adm}=X$ et $X_{\rm adm}(a)=X\bck \{a\}$
pour tout $a\in X$.
Tout $2$-uplet de $X$ est admissible et $\Si_q(X) = {\rm CRN}_q(X) = {\rm CR}_q(X)$
est une vari\'et\'e compacte lisse de dimension $2r$.
Un isomorphisme de Veronese $\P^{r+1}\rightarrow X$
induit un isomorphisme de la paire $(\P^{r+1},\Si_1(\P^{r+1}))$,
o\`u  $\Si_1(\P^{r+1})$  est la grassmannienne $\G_{r,2}$
des droites de $\P^{r+1}$, sur la paire $(X,\Si_q(X))$.
Les propri\'et\'es de la seconde paire qu'on utilisera sont simplement des traductions 
de propri\'et\'es bien connues de la premi\`ere.

\ssct{Projections sur une vari\'et\'e de Veronese}

On sait que, si $n=2$, $X$ est une vari\'et\'e de Veronese d'ordre $q$.
On suppose maintenant $n\geq 3$. 
\bd
Une {\em pond\'eration} est un \'el\'ement $(\rho_1,\ldots,\rho_{n-1})$ de $\N^{n-1}$
tel qu'on ait  $\rho_i=\rho$ pour $m$ et $\rho_i=\rho-1$ pour $n-1-m$ valeurs de $i\in \{1,\ldots,n-1\}$.

Une pond\'eration \'etant fix\'ee, si $p\in \{1,\ldots,n-2\}$
et si $\max (\rho_{p+1},\ldots,\rho_{n-1})\geq 1$, on associe 
\`a tout $p$-uplet admissible ${\boldsymbol a}=(a_1,\ldots,a_p)$  
les projections $\t_{{\boldsymbol a}}$ de centre la somme directe 
$\oplus_{i=1}^{p} \, X_{a_i}(\rho_i)$.
On dit que $\t_{\boldsymbol a}$ est {\em une projection osculatrice associ\'ee \`a 
$\boldsymbol a$.}
\ed
La condition $\max (\rho_{p+1},\ldots,\rho_{n-1})\geq 1$ est l\`a pour \'eviter 
des projections dont la restriction \`a $X$ d\'eg\'en\`ere. On 
pourra toujours supposer sans dommage qu'on a  $\rho_{n-1}\geq 1$.

\sk
On consid\`ere d'abord des projections osculatrices associ\'ees \`a des $(n-2)$-uplets
admissibles. Ce sont de loin les plus importantes.

Pour fixer les id\'ees, on choisit dans ce chapitre, sauf la derni\`ere section,
la pond\'eration 
\beq
\label{ponderation}
(\rho_1,\ldots,\rho_{n-1}) = (\rho-1,\ldots,\rho-1,\rho,\ldots,\rho),
\eeq
(on utilisera les deux notations) mais on aurait un r\'esultat analogue au suivant pour une pond\'eration 
g\'en\'erale avec $\rho_{n-1}\geq 1$, la projection de $X$ \'etant alors une vari\'et\'e de 
Veronese d'ordre $\rho_{n-1}$.  

\sk
Le r\'esultat suivant est l'outil-cl\'e pour toute la suite.
\bpr
\label{proj}
Soit $\boldsymbol a$ un $(n-2)$-uplet admissible de $X$ et $\t_{\bs a}$ une projection 
osculatrice associ\'ee.
Elle induit une application birationnelle $\t_{\bs a}: X\dasharrow X'$, o\`u $X' = \t_{\bs a}(X)$ 
est une vari\'et\'e de Veronese d'ordre $\rho$, et un diff\'eomorphisme de $X_{\rm adm}({{\boldsymbol a}})$ sur son image.

Si $\,{\boldsymbol x}\in X_{\rm adm}^{(2)}({\boldsymbol a})$, il existe un seul
cycle $C\in {\rm CR}_q(X)$ qui contient $({\boldsymbol a},{\boldsymbol x})$.
Il appartient \`	a $\Si_q(X;{{\boldsymbol a}})$ et $\,\t_{\bs a}(C)$ est l'unique \'el\'ement 
de $\Si_\rho(X')$ qui contient $\t_{\bs a}({\boldsymbol x})$. De plus, 
la projection $\t_{\bs a}$ est d\'efinie comme morphisme au voisinage de $C\bck {\boldsymbol a}$ 
et $\t_{\bs a} :C\bck {\boldsymbol a} \rightarrow \t_{\bs a}(C)$
est la restriction d'un isomorphisme de $C$ sur $\t_{\bs a}(C)$.
\epr
\bpf
Comme ${\boldsymbol a}=(a_1,\ldots,a_{n-2})$ est fix\'e, on ne note pas ${\boldsymbol a}$ en indice.
On consid\`ere une projection de centre $\Q=\oplus_{i=1}^{n-2} \, X_{a_i}(\rho_i)$, soit 
\beq
\label{projP}
\t: \; \P^N \dasharrow X_c(\rho),
\eeq
o\`u $c\in X_{\rm adm}({\boldsymbol a})$. Bien s\^ur le choix de la cible est sans importance.

\sk
On note $X'$ l'image $\t(X)$ de $X$ et $\t: X\dasharrow X'$
l'application induite\footnote{
Soit $\wh{X}$ une vari\'et\'e irr\'eductible et $\t: \wh{X}\dasharrow \P^m$ une application rationnelle.
Celle-ci est d\'efinie comme morphisme sur un ouvert dense $\wh{X}'$ de $\wh{X}$.
Suivant l'usage, si $X$ est une sous-vari\'et\'e irr\'eductible de $\wh{X}$ 
{\em et si $X\cap \wh{X}'$ est non vide}, on note $\t(X)$ 
l'adh\'erence de Zariski de $\t(X\cap \wh{X}')$ 
et $\t : X\dasharrow \t(X)$ l' application rationnelle 
induite par  la restriction de $\t$ \`a $X\cap\wh{X}'$.}.
Comme $\lan X\ran =\P^N$, on a $\lan X'\ran \, = X_c(\rho)$, donc  
$\dim \, \lan X'\ran = \pi_{r,2}(\rho)$.

\sk
Soit ${\boldsymbol x}=(x_1,x_2)\in X_{\rm adm}^{(2)}({\boldsymbol a})$.
Comme $\P^N=\Q\oplus X_{x_1}(\rho)$, la projection $\t$ induit un diff\'eomorphisme local 
du germe de $X$ en $x_1$ sur un germe lisse et $\rho$-r\'egulier $\tilde{X}'\subset X'$ en $\t(x_1)$,
de dimension $r+1$. On ne sait pas si $X'$ est lisse aux  points $\t(x_1)$ et $\t(x_2)$
mais c'est bien s\^ur vrai
pour un choix g\'en\'erique de ${\boldsymbol x}\in X_{\rm adm}^{(2)}({\boldsymbol a})$.  

\bk
Soit  $C$ {\em un} \'el\'ement de ${\rm CR}_q(X)$ qui contient $({\boldsymbol a},{\boldsymbol x})$.
Compte tenu du Lemme \ref{rec-L1}, $C$ est une courbe admissible.
Elle n'est pas contenue dans $\Q$ et d'apr\`es ce qui pr\'ec\`ede, $\t(C)$
est une courbe irr\'eductible $C'\subset X'$.

Si $H'$ est un hyperplan de $X_c(\rho)$ et $H=\Q\oplus H'$,
la courbe $C$ coupe $H$ en $a_i$ avec au moins la multiplicit\'e $\rho_i+1$,
$\,i=1,\ldots,n-2$. Le  degr\'e de son image $C'$ est donc major\'e par
$q - \sum_{i=1}^{n-2}(\rho_i+1)$,
c'est-\`a-dire par $\rho$. 
D'autre part, l'image du germe de $C$ en $x_1$ 
est un germe de courbe lisse contenu dans le germe $\rho$-r\'egulier 
$\tilde{X}'$ et donc $\rho$-r\'egulier. Ainsi $C'$ est une courbe 
$\rho$-r\'eguli\`ere de degr\'e $\leq \rho$, une courbe rationnelle normale 
de degr\'e $\rho$. 

\sk
Comme les courbes $C$ et $C'$ sont lisses, l'application rationnelle
induite 
$\t: C\dasharrow C'$ est en fait un morphisme qu'on note $\t_C: C\rightarrow C'$,
pour \'eviter les confusions ($\t_C({\boldsymbol a})$ d\'epend de $C$).
Reprenons le d\'ecompte de $C\cap H$. On a $C\cap \Q={\boldsymbol a}$, \,
sinon on trouverait que le degr\'e de $C'$
est au plus $\rho-1$, et comme le cardinal de $C\cap (H\bck \Q)$ est $\leq \rho$,
$\, \t_C$ est un isomorphisme. En r\'esum\'e, 

{\em la projection $\t$ est d\'efinie comme morphisme 
au voisinage de $C\bck {\bs a}$ et sa restriction \`a $C\bck {\bs a}$
se prolonge en un isomorphisme de $C$ sur $\t(C)$, un \'el\'ement de ${\rm CRN}_\rho(X')$.}

\sk
Dans le raisonnement qu'on vient de faire, ${\boldsymbol x}\in X_{\rm adm}^{(2)}({\boldsymbol a})$
\'etait quelconque et, comme on a dit, pour ${\boldsymbol x}$ voisin 
d'un $2$-uplet ${\boldsymbol x}_0$ bien choisi, les deux points distincts $\t(x_1)$
et $\t(x_2)$ appartiennent \`a $X'_{\rm reg}$. 
Ainsi, pour toute paire $(x'_1,x'_2)$ d'un ouvert non vide de $X'^2$,
il existe une courbe $C'\in {\rm CRN}_\rho(X')$ qui passe par $x'_1$  et $x'_2$.
Comme $X'$ engendre un espace de dimension $\pi_{r,2}(\rho)$,
$X'$ est une vari\'et\'e de la classe $\cal{X}_{r+1,2}(\rho)$.
Compte tenu du Th\'eor\`eme~\ref{Th2}, 

{\em $X'$ est une vari\'et\'e de Veronese 
d'ordre $\rho$, en particulier $X'$ est lisse}.

\sk
Soit $x_1,x_2\in X_{\rm adm}({\boldsymbol a})$. Comme $X_{\rm adm}({\boldsymbol a},x_1)$ et $X_{\rm adm}({\boldsymbol a},x_2)$
sont des ouverts denses de $X$, on peut consid\'erer un point $x_0$
de $X$ tel que les $n$-uplets ${\boldsymbol x}_1=({\boldsymbol a},x_0,x_1)$
et ${\boldsymbol x}_2=({\boldsymbol a},x_0,x_2)$ soient admissibles. Soit
$C_1,C_2\in {\rm CR}_q(X)$ des cycles qui contiennent respectivement ${\boldsymbol x}_1$ et ${\boldsymbol x}_2$.
Compte tenu du Lemme \ref{rec-L1}, ce sont des courbes admissibles. 
Si $\t(x_1)=\t(x_2)$, les courbes admissibles $C_1$ et $C_2$
ont la m\^eme image par $\t$, {\em la} courbe $C'\in {\rm CRN}_\rho(X')$
qui passe par $\t(x_0)$ et $\t(x_1)=\t(x_2)$. Ainsi $C_1$ et $C_2$ ont le m\^eme germe en $x_0$,
donc $C_1=C_2$. Comme la restriction de $\t$ \`a $C_1\bck {\bs a}$ est injective, on obtient $x_1=x_2$.

\sk
Ceci montre que la restriction de $\t$
\`a $X_{\rm adm}({{\boldsymbol a}})$ est injective. Comme elle est 
de rang constant $r+1$, c'est un diff\'eo\-morphisme de $X_{\rm adm}({\boldsymbol a})$ sur son image dans $X'$.
En particulier, $\, \t: X\dasharrow X'$ est birationnelle. 

Ceci montre aussi que, pour tout ${\bs x}\in X^{(n)}_{\rm adm}$, il 
existe un et un seul cycle $C\in {\rm CR}_q(X)$, n\'ecessairement un \'el\'ement
de ${\rm CRN}_q(X)$,  dont le support contient $\bs x$.
\epf

\ssct{Une traduction et une premi\`ere application}

La Proposition \ref{proj} a une jumelle qu'on obtient en composant 
l'application osculatrice $\t_{\bs a}: X\dasharrow X'$
avec l'inverse d'un isomorphisme de Veronese de $\P^{r+1}$ sur $X'$.
\bpr
\label{proj-jumelle}
Soit $n\geq 3$, $\,X$ une vari\'et\'e de la classe $\cal{X}_{r+1,n}(q)$
et $\boldsymbol a$ un $(n-2)$-uplet admissible de $X$. Il existe 
une application birationnelle $\pi_{\boldsymbol a} : X\dasharrow \P^{r+1}$,
telle $\pi_{\bs a}(C)$ est une droite de $\P^{r+1}$ pour toute courbe $C\in \Si_q(X;{\boldsymbol a})$.

L'application $\pi_{\boldsymbol a}$ induit un diff\'eomorphisme de $X_{\rm adm}({{\boldsymbol a}})$
sur son image dans $\P^{r+1}$ et si $C\in \Si_q(X;{\boldsymbol a})$,
$\,\pi_{\boldsymbol a}$ est un morphisme au voisinage de $C\bck {\boldsymbol a}$ 
dont la restriction \`a $C\bck {\boldsymbol a}$ se prolonge en un isomorphisme 
de la courbe $C$ sur la droite $\pi_{\boldsymbol a}(C)$.
\epr
On dira aussi que l'application birationnelle $\pi_{\boldsymbol a}: X\dasharrow \P^{r+1}$
est {\em associ\'ee \`a $\boldsymbol a$}. Elle  est d\'etermin\'ee \`a la composition 
\`a gauche pr\`es par un automorphisme de $\P^{r+1}$. 
Selon la situation on utilisera la Proposition \ref{proj}
ou la Proposition \ref{proj-jumelle}.

\sk
L'un des points de la Proposition \ref{proj} s'\'enonce sans r\'ef\'erence \`a une 
projection osculatrice. {\em Si ${\boldsymbol x}\in X^{(n)}_{\rm adm}$, 
un seul \'el\'ement de ${\rm CR}_q(X)$ contient $\boldsymbol x$, c'est une courbe admissible.}
L'adh\'erence de Zariski $\overline{\Si}_q(X)$ de $\Si_q(X)$
dans ${\rm CR}_q(X)$ est donc la seule composante $Z$ de $\text{\rm CR}_q(X)$, telle
que la projection $\{(C,{\boldsymbol x})\in Z\times X^n, \;\; {\boldsymbol x}\subset |C|^n\} \rightarrow X^n$
soit surjective.

Consid\'erons la vari\'et\'e d'incidence $I = \{(C, {\boldsymbol x})\in  \overline{\Si}_q(X)\times X^n, \; {\boldsymbol x}\subset |C|\}$
et la projection $p: I\rightarrow X^n$. L'image r\'eciproque d'un \'el\'ement ${\bs x}\in X^{(n)}_{\rm adm}$
est r\'eduite \`a un point. Il en r\'esulte qu'il existe une et une seule
composante irr\'eductible $I_0$ de $I$ telle que le morphisme induit $p: I_0\rightarrow X^n$
soit surjectif. Pour la m\^eme raison, c'est une application birationnelle,
dont l'inverse $p^{-1} : X^n \dasharrow I_0$ 
n'a pas de point d'ind\'etermination sur l'ouvert lisse $X^{(n)}_{\rm adm}$.
On a donc le r\'esultat suivant.
\ble
\label{anal}
L'adh\'erence de Zariski $\overline{\Si}_q(X)$ de l'ensemble $\Si_q(X)$
des courbes admissibles de $X$ est la seule composante irr\'eductible de $\text{\rm CR}_q(X)$
qui v\'erifie que, pour ${\boldsymbol x}\in X^n$, il existe  
un \'el\'ement de cette composante dont le support contient $\boldsymbol x$.
L'application 
\beq
\label{C}
\g:  \; X_{\rm adm}^{(n)} \rightarrow \Si_q(X), 
\eeq
qui \`a tout ${\boldsymbol x}\in X^{(n)}_{\rm adm}$ 
associe l'unique cycle $C\in {\rm CR}_q(X)$ dont le support contient $\boldsymbol x$,
en fait une courbe admissible,
est analytique et ouverte.
\ele

\ssct{Applications}

Comme on ne fait pas d'hypoth\`ese de r\'egularit\'e sur $X$,
on pr\'ef\`ere dans un premier temps voir $X$ comme une sous-vari\'et\'e de $\P^N$
et $\Si_q(X)$ comme une sous-vari\'et\'e de 
la vari\'et\'e lisse (non ferm\'ee) ${\rm CRN}_q(\P^N)$.

\sk
On s'int\'eresse \`a l'application analytique et ouverte $\g:  \; X_{\rm adm}^{(n)} \rightarrow \Si_q(X)$,
d\'efinie  dans le lemme pr\'ec\'edent. 

\sk
Si $\wh{\bs a}\in X^{(n)}_{\rm adm}$,
l'espace tangent \`a la vari\'et\'e lisse ${\rm CRN}_q(\P^N)$
en $\g(\wh{\boldsymbol a})$ s'identifie \`a un sous-espace,
en fait \`a tout l'espace mais ce n'est pas important ici, 
de l'espace 
des sections analytiques globales du fibr\'e normal $N_{\g(\wh{\boldsymbol a})}\P^N$
\`a $\g(\wh{\boldsymbol a})$ dans $\P^N$, qu'on note $H^0(\g(\wh{\boldsymbol a}),N_{\g(\wh{\boldsymbol a})}\P^N)$.
{\em Dans la suite, on fait cette identification :}
$$
T_{\g(\wh{\bs a})}{\rm CRN}_q(\P^N) \subset H^0(\g(\wh{\boldsymbol a}),N_{\g(\wh{\boldsymbol a})}\P^N).
$$
Rappelons 
comment elle s'interpr\`ete. 

\sk
Soit $\Lambda$ un germe de vari\'et\'e lisse et $v: \Lambda \rightarrow {\rm CRN}_q(\P^N)$
une application analytique. Quitte \`a choisir le repr\'esentant $\Lambda$ assez petit, 
il existe un rel\`evement $V: \Lambda\times \P^1 \rightarrow \P^N$ de $v$,
c'est-\`a-dire une application analytique telle que, pour $\l \in \Lambda$ fix\'e, l'application 
$V_\l: t\mapsto V(\l,t)$ est  un isomorphisme de $\P^1$ sur la courbe $v(\l)$. 

La diff\'erentielle 
$dv(\l): T_\l \Lambda \rightarrow T_{v(\l)}{\rm CRN}_q(\P^N)$
de $v$ en $\l\in \Lambda$, associe \`a $\nu \in T_{\l}\Lambda$ la section  de $N_{v(\l)}\P^N$
donn\'ee par :
\beq
\label{diff}
p\in v(\l), 
\qquad 
(dv(\l)\cdot \nu)(p)  = \fr{\pl V}{\pl \l}(\l,V_\l^{-1}(p))\cdot \nu \;\; {\rm modulo} \;\; T_p v(\l).
\eeq
On v\'erifie sans peine, et c'est bien connu, que le r\'esultat ne d\'epend pas 
du rel\`evement choisi et qu'il est local.
On entend par l\`a que si $p_0\in v(\l_0)$, pour calculer la section $dv(\l_0)\cdot \nu$
au voisinage de $p_0$ dans $v(\l_0)$, un rel\`evement local de $v$ suffit,
{\em i.e.} une application analytique $V : \Lambda \times \omega \rightarrow \P^N$, o\`u $\omega$
est un  voisinage d'un point $t_0\in \P^1$, telle qu'on ait  $\;V(\l_0,t_0)=p_0$ et que, 
pour tout $\l \in \Lambda$, $\;t\mapsto V(\l,t)$ soit un plongement $V_\l : \omega \rightarrow v(\l)$. 

\bk
Dans cette section, on d\'emontre un r\'esultat sur la r\'egularit\'e de la vari\'et\'e 
$X$ et un autre sur la r\'egularit\'e et la structure de la vari\'et\'e $\Si_q(X)$
des courbes admissibles. Les d\'emonstrations reposent  sur la consid\'eration 
de certaines applications partielles induites par 
l'application $\g$ d\'efinie par (\ref{C}) et sur la Proposition \ref{proj-jumelle}.

\sk
On d\'ecrit la situation mod\`ele. On se donne un $n$-uplet admissible 
$$
\wh{\bs a} = (a_1,\ldots,a_n),
$$
qu'on note aussi 
$$
\wh{\bs a}=({\bs a},{\bs b}), \;\;\; {\bs a} = (a_1,\ldots,a_{n-2}), \;\; {\bs b}=(b_1,b_2).
$$
Soit $\pi: X\dasharrow \P^{r+1}$ une application birationnelle
associ\'ee \`a ${\bs a}$. On note ${\bs b}'= \pi({\bs b})$,
c'est-\`a-dire $b'_1=\pi(b_1)$, $\,b'_2=\pi(b_2)$.
On introduit les {\em germes} d'applications  
\beq
\label{partiel}
v : \, (X^2,{\bs b}) \rightarrow {\rm CRN}_q(\P^N), 
\qquad 
l:  ((\P^{r+1})^2,{\bs b}')  \rightarrow {\rm CRN}_1(\P^{r+1}),
\eeq
o\`u $v({\bs x})= \g({\boldsymbol a},{\bs x})$ pour ${\bs x}$  voisin de $\bs b$ 
et $l({\bs x}')$ est la droite qui contient ${\bs x}'$, pour ${\bs x}'$ voisin de ${\bs b}'$.
On introduit des {\em germes} de rel\`evements de $v$ et de $l$ 
\beq
\label{Partiel}
V: \; (X^2,{\bs b}) \times \P^1 \rightarrow \P^N, 
\qquad
L: \; ((\P^{r+1})^2,{\bs b}') \times \P^1 \rightarrow \P^{r+1}.
\eeq
On les choisit {\em normalis\'es de fa\c{c}on concordante}.
Par exemple, on se donne un troisi\`eme point $b_3\in X_{\rm adm}({\bs  a})\cap \g({\wh{\bs a}})$
et on impose les conditions de normalisation 
\beq
\label{choixY}
V({\bs x},t_k) \in Y_{b_k}, \;\;\;  L({\bs x},t_k) \in \pi(Y_{b_k}), \qquad k=1,2,3,
\eeq
o\`u $Y_{b_k}$ est un germe d'hypersurface lisse donn\'e transverse \`a la courbe $\g(\wh{\bs a})$
en $b_k$ et par exemple $(t_1,t_2,t_3)=(1,0,\infty)$.

\sk
Comme le morphisme $v({\bs x})\bck {\boldsymbol a}\rightarrow l(\pi({\bs x}))$
induit par $\pi$ est la restriction d'un isomorphisme de $v({\bs x})$ sur $l(\pi({\bs x}))$
et compte tenu de la concordance des rel\`evements, on a :
\beq
\label{partiel-proj}
V({\bs x},t)\notin {\boldsymbol a} \;\; \Rightarrow \;\; \pi(V({\bs x},t)) = L(\pi({\bs x}),t).
\eeq

\bk
Cette remarque est  suffisante pour d\'emontrer le 
r\'esultat suivant, qui sera important dans le Chapitre~3. Il implique en particulier
que le nombre d'intersection $Y\cdot C$ est bien d\'efini, si $Y$
est une hypersurface de $X$ et $C\not\subset Y$ une courbe admissible.
\bt
\label{Xlisse}
Une vari\'et\'e $X\in \cal{X}_{r+1,n}(q)$
est lisse au voisinage de toute courbe admissible de $X$.
Tout $n$-uplet de points deux-\`a-deux distincts d'une courbe admissible est admissible.
L'ouvert  $X_{\rm adm}$ des points admissibles de $X$ est donc aussi la 
r\'eunion des courbes admissibles de $X$.
\et
\bpf
On  consid\`ere la situation qu'on vient de d\'ecrire en pr\'eambule 
\`a l'\'enonc\'e et on conserve ses notations, en particulier $(\ref{Partiel})$. 
Il est commode de choisir la normalisation (\ref{choixY}) 
des rel\`evements normalis\'es concordants  $V$ et $L$.

\sk
Comme un point admissible appartient \`a $X_{\rm reg}$, 
il suffit, pour d\'emontrer la premi\`ere partie de l'\'enonc\'e,  de montrer
que $X$ est lisse au voisinage de tout point  $p_0\in  \g(\wh{\bs a}) \bck \wh{{\bs a}}$.
Notons $Y$ pour $Y_{b_2}$, $\, Y'$ pour $\pi(Y_{b_2})$ et : 
$$
V_\star(y,t) = V(b_1,y,t), \;\;\; L_\star(y',t)  = L(b'_1,y',t),  \qquad y\in Y, \;\; y'\in Y',\;\; t\in \P^1\bck \{0\}.
$$
Soit $p_0=V({\bs b},t_0)\in \g(\wh{\bs a})\bck \wh{{\bs a}}$ donc avec $t_0\neq 0$. 
L'application $L_\star$ est de rang constant $r+1$.
Compte tenu de l'\'egalit\'e (\ref{partiel-proj}), 
il en d\'ecoule 
que $(y,t)\mapsto \pi(V_\star(y,t))$ est de rang maximal $r+1$ en $(b_2,t_0)$. 
On en d\'eduit d'abord que $(y,t)\mapsto V_\star(y,t)$
est de rang maximal $r+1$ au m\^eme point donc param\`etre un germe 
de vari\'et\'e lisse $\tilde{X}\subset X$ en $p_0$, de dimension $r+1$,
ensuite que la restriction de $\pi$ \`a $\tilde{X}$ est de rang $r+1$
en $p_0$, donc induit un diff\'eomorphisme local de $\tilde{X}$ sur son image
dans $\P^{r+1}$.

\sk
En particulier, l'image de $\tilde{X}$ contient un voisinage $U$
de $\pi(p_0)$ dans $\P^{r+1}$. 
Si $\tilde{X}$ n'\'etait pas le germe de $X$ en $p_0$, il existerait en $p_0$
un {\em autre} germe $\breve{X}\subset X$ de vari\'et\'e de dimension $r+1$
et, dans la situation qu'on consid\`ere, l'image de $\breve{X}$
contiendrait un ouvert non vide de $U$ puisque 
$\pi: X \dasharrow \P^{r+1}$ est dominante, en contradiction
avec le fait que $\pi: X\dasharrow \P^{r+1}$ est birationnelle.

\sk
Pour d\'emontrer la deuxi\`eme partie de l'\'enonc\'e il suffit,
dans la m\^eme situation, de montrer que le $n$-uplet
$({\bs a},p_0,b_2)$ est admissible. Le r\'esultat annonc\'e, que tout 
$n$-uplet de $\g(\wh{\bs a})$ est admissible, en d\'ecoule 
par it\'eration et en changeant de $(n-2)$-uplet $\bs a$.

Comme $p_0\in X_{\rm reg}$, compte tenu de la D\'efinition \ref{D3.1}
et par sym\'etrie, il suffit de montrer que les courbes admissibles 
qui passent par $\bs a$ et $b_2$ recouvrent un voisinage de $p_0$,
ce qu'on vient de faire, et que les courbes admissibles qui passent
par $\bs a$ et $p_0$ recouvrent un voisinage de $b_2$, ce qui 
maintenant est  clair, puisque les droites qui
passent par $\pi(p_0)$ et un point voisin de $b'_1$ recouvrent un voisinage de $b'_2$.
\epf

\sk
On revient \`a la situation mod\`ele 
d\'ecrite avant l'\'enonc\'e pr\'ec\'edent\footnote{Le r\'esultat qu'on 
a en vue est important seulement dans le Chapitre 4. Les consid\'erations qui suivent
ne servent pas dans la d\'emonstration du Th\'eor\`eme \ref{Th4}.}.
On s'int\'eresse maintenant aux d\'eriv\'ees des applications (\ref{partiel}),
$$
dv({\bs b})  : \, T_{\bs b}X^2  \rightarrow T_{v({\bs b})}{\rm CRN}_q(\P^N),
\qquad 
dl({\bs b}') : \, T_{{\bs b}'}(\P^{r+1})^2  \rightarrow T_{l({\bs b'})}{\rm CRN}_1(\P^{r+1}),
$$
en $\bs b$ et en ${\bs b}'=\pi({\bs b})$.
Les propri\'et\'es
de $dl({\bs b}')$ sont connues, elles concernent la grassmannienne des 
droites de $\P^{r+1}$. D'autre part, en d\'erivant l'identit\'e de (\ref{partiel-proj}) en $\bs b$, dans la direction 
${\bs \nu}\in T_{\bs b}X^2$, voir aussi la formule~(\ref{diff}), on obtient la relation 
\beq
\label{DIFF}
p\in v({\bs b})\bck {\bs a}, 
\qquad 
d\pi(p)\cdot ((dv({\bs b})\cdot {\bs \nu})(p)) = (dl({\bs b}')\cdot {\bs \nu}'))(\pi(p)),
\eeq
o\`u ${\bs \nu}=(\nu_1,\nu_2) \in T_{(b_1,b_2)}(X \times X)$
et ${\bs \nu}'= (d\pi(b_1)\cdot \nu_1,d\pi(b_2)\cdot \nu_2)$.

\sk
Notons $l$ la droite $l({\bs b}')$. Le fibr\'e normal $N_l\P^{r+1}$ est une somme directe 
de $r$ fibr\'es en droites de degr\'e $1$. On a :
$$
\ker dl({\bs b}') = T_{(b'_1,b'_2)}(l\times l), 
\qquad 
{\rm im} \, dl({\bs b}') = H^0(l,N_l\P^{r+1}).
$$
Une section globale non nulle de $N_l\P^{r+1}$ s'annule au plus en un point, simplement.

\sk
Notons $C$ la courbe $\g(\wh{\bs a})=v({\bs b})$. Pour tout $p\in C\bck {\bs a}$,
l'application birationnelle $\pi : X\dasharrow \P^{r+1}$
induit un diff\'eomorphisme local au voisinage de $p$.
La d\'eriv\'ee $d\pi(p)$ induit des isomorphismes de $T_pX$ sur 
$T_{\pi(p)}\P^{r+1}$ et de $T_pC$ sur $T_{\pi(p)}l$, donc aussi
de la fibre $(N_CX)_p$ de $N_CX$ en $p$ sur la fibre $(N_l\P^{r+1})_{\pi(p)}$
de $N_l\P^{r+1}$ en $\pi(p)$. Cette remarque s'applique en particulier 
\`a la relation $(\nu'_1,\nu'_2)=(d\pi(b_1)\cdot \nu_1,d\pi(b_2)\cdot \nu_2)$.

\sk
Du c\^ot\'e de $dv({\bs b})$, on a :
$$
\ker dv({\bs b}) = T_{(b_1,b_2)}(C\times C), 
\qquad 
{\rm im} \, dv({\bs b}) \subset H^0(C,N_CX) \subset H^0(C,N_C\P^N).
$$
La premi\`ere relation r\'esulte des propri\'et\'es de  $dl({\bs b}')$  et de (\ref{DIFF}),
la seconde rappelle qu'on a identifi\'e $T_C\,{\rm CRN}(\P^N)$
\`a un sous-espace de $H^0(C,N_C\P^N)$
et que $X$ est lisse au voisinage de $C$.
Il r\'esulte encore de (\ref{DIFF}) qu'on a un isomorphisme $\phi : {\rm im} \,  dv({\bs b}) \rightarrow {\rm im} \, dl({\bs b}')$, 
compatible 
avec l'action naturelle de $d\pi(p) :  (N_CX)_p \rightarrow (N_l\P^{r+1})_{\pi(p)}$ 
en dehors de $\bs a$. On a donc un isomorphisme  
\beq
\label{phi1}
\phi : \; {\rm im} \,  dv({\bs b}) \rightarrow H^0(l,N_l\P^{r+1}) \, ;
\qquad 
\phi(\xi)(\pi(p)) = d\pi(p)\cdot \xi(p) \;\; \text{si $p\in C\bck {\bs a}$}.
\eeq
On a presque obtenu le lemme suivant, dont l'\'enonc\'e ne fait plus r\'ef\'erence 
aux applications osculatrices.
Rappelons qu'une section non nulle de $N_CX$ engendre un sous-fibr\'e en droites 
$\L(\xi)$ de $N_CX$ dont elle est une section\footnote{
Rappelons pourquoi. Au voisinage d'un z\'ero d'ordre $p\geq 1$
de $\xi$, on identifie $C$ \`a un voisinage de $0\in\C$ et on \'ecrit $\xi(t)=t^p e(t)$ ;
alors $e$ d\'efinit un fibr\'e en droites sur un voisinage $U$ de $0$,
qui prolonge celui que $\xi$ d\'efinit sur $U\bck\{0\}$.}.
\ble
\label{modele}
Soit $C$ une courbe admissible de $X$ et $\wh{\bs a}=(a_1,\ldots,a_n)$ un $n$-uplet de~$C$.
Pour $i=1,\ldots,n$, soit $E_i \subset H^0(C,N_CX)$  l'image de la restriction de la d\'eriv\'ee $d\g(\wh{\bs a})$
au sous-espace $\{0\}\times \cdots \times \{0\}\times T_{a_i}X \times \{0\} \times \cdots \times \{0\}$
de $T_{\wh{\bs a}}X^n$. On a :
\be
\item pour tout $i\in \{1,\ldots,n\}$ l'espace $E_i$ est de dimension $r$
et une section non nulle $\xi\in E_i$ de $N_CX$ s'annule simplement aux points de $\wh{\bs a}\bck \{a_i\}$
et nulle part ailleurs ;
\item pour tout $i,j\in \{1,\ldots,n\}$ et toute section non nulle $\xi\in E_i$,
il existe une section non nulle $\eta\in E_j$ qui d\'efinit le m\^eme sous-fibr\'e en droites 
de $N_CX$ : $\L(\eta)=\L(\xi)$.
\ee
\ele 
\bpf
Par sym\'etrie, on peut supposer $i=n$. Notons $a_n=b_2$. Si $a_{n-1}=b_1$, on retrouve la
situation mod\`ele et $E_n\subset H^0(C,N_CX)$ est aussi l'image de la restriction 
de $dv({\bs b})$ au sous-espace $\{0\}\times T_{b_2}X$ de $T_{\bs b}X^2$.
Compte tenu de (\ref{phi1}), $E_n$ est de dimension $r$
et une section non nulle $\xi\in E_n$ de $N_CX$ s'annule en $a_1,\ldots,a_{n-2}$,
simplement en $a_{n-1}$ et nulle part ailleurs. 
Par sym\'etrie, on obtient la premi\`ere partie de l'\'enonc\'e.

\sk
La seconde partie est \'evidente si $i=j$. Par sym\'etrie, on peut supposer $i=n$,
$j=n-1$ et retrouver la situation mod\`ele. Soit $\xi\in E_n$ une section 
non nulle de $N_CX$. Son image $\xi'=\phi(\xi)$ engendre un sous-fibr\'e en droite 
$\L(\xi')$ de $N_l\P^{r+1}$, de degr\'e $1$, et l'image de $dl({\bs b}')$ contient les sections 
de ce fibr\'e. Soit $\eta'$ une section non nulle de $\L(\xi')$ qui
s'annule en $b'_1=\pi(a_{n-1})$. Comme $\phi$ est un isomorphisme et compte tenu de (\ref{phi1}), $\eta'=\phi(\eta)$
pour un $\eta \in E_{n-1}$. D'autre part $\eta$ est une section de $\L(\xi)$
d'apr\`es (\ref{phi1}). D'o\`u la seconde partie du lemme.
\epf

\sk
La cons\'equence suivante est importante.
\bt
\label{structure-Si}
La vari\'et\'e alg\'ebrique $\ov{\Si}_q(X)$ est de dimension $rn$
et l'ouvert dense $\Si_q(X)$ des courbes admissibles de $X$ est contenu dans 
sa partie lisse.

Pour toute courbe admissible $C$ de $X$, son fibr\'e normal $N_CX$ 
est une  somme  directe de $r$ sous-fibr\'es en droites de degr\'e $n-1$ et 
$T_C\Si_q(X)=H^0(C,N_CX)$.

Si ${\bs a}$ est un $p$-uplet de $C$,
l'ensemble
$\Si_q(X;{\boldsymbol a})$ des courbes admissibles qui contien\-nent $\boldsymbol a$
est une sous-vari\'et\'e lisse de dimension $(n-p)r$ de $\Si_q(X)$
et  $T_C\Si_q(X;{\boldsymbol a})$ est l'espace des sections de $N_CX$ 
qui s'annulent sur $\bs a$.
\et
\bpf
L'application $\g: X_{\rm adm}^{(n)} \rightarrow \Si_q(X)$ est ouverte donc 
l'image du germe de $X^{(n)}_{\rm adm}$ en un point est 
le germe de $\Si_q(X)$ au point image.
Compte tenu du th\'eor\`eme du rang, il suffit, pour 
montrer que $\Si_q(X)$ est lisse de dimension $rn$, de montrer que l'application 
$\g:  X_{\rm adm}^{(n)} \rightarrow {\rm CRN}_q(\P^N)$
est de rang constant~$rn$. 

\sk
Soit $\wh{\bs a}\in X^{(n)}_{\rm adm}$. Notons $C=\g(\wh{\bs a})$ et 
$$
E = {\rm im}\, d\g(\wh{\bs a}) \subset H^0(C,N_CX).
$$
On d\'efinit les sous-espaces $E_1,\ldots,E_n$ de $E$
comme dans le lemme pr\'ec\'edent. Ils sont de dimension $r$ et  
$E = \sum_{i=1}^n E_i$. Cette somme est directe : si $\xi_i\in E_i$,
et $\sum_{i=1}^n \xi_i=0$, on v\'erifie que chaque terme s'annule en $a_1,\ldots,a_n$,
donc est nul d'apr\`es la premi\`ere partie du Lemme \ref{modele}.
On a donc 
$$
E = \oplus_{i=1}^n E_i.
$$
Ceci donne en particulier la premi\`ere partie de l'\'enonc\'e. 

\sk
La premi\`ere partie du Lemme \ref{modele} montre aussi que,
pour tout $i\in \{1,\ldots,n\}$,  
la {\em fibre}
$E_i(p)=\{\xi(p), \;\; \xi\in E_i\}$ de $E_i$ v\'erifie $E_i(p)=(N_CX)_p$ si 
$p\in C\bck \wh{\bs a}$ ou si $p=a_i$.
Il en r\'esulte que la 
fibre  $E(p)$ de $E$ est de dimension $r$ pour tout $p\in C$.

\sk
Pour d\'eterminer la structure de $N_CX$, 
partons par exemple d'une base $(\xi^\a)_{\a=1}^r$
de $E_n$. Pour chaque $\a$, le fibr\'e $\L(\xi^\a)$ 
engendr\'e par $\xi^\a$ est de degr\'e $n-1$. Les espaces 
$$
F^\a = H^0(C,\L(\xi^\a)) \subset H^0(C,N_CX), \qquad \a=1,\ldots,r,
$$
sont de dimension $n$ et \'evidemment en somme directe.
Pour $\a$ fix\'e et compte tenu de la seconde partie du Lemme \ref{modele},
il existe, pour tout $i\in \{1,\ldots,n\}$ une section non nulle $\eta^\a_i\in E_i$
de $\L(\xi^\a)$. Les sections $\eta^\a_1,\ldots,\eta^\a_n \in E$ 
sont lin\'eairement ind\'ependantes donc engendrent $F^\a$.
En particulier $F^\a\subset E$. 

\sk
Ceci montre que $E=\oplus_{\a=1}^r F^\a$ puisque le premier membre contient le second 
et que les deux membres ont la m\^eme dimension.
Comme les fibres de $E$ sont de dimension~$r$, les 
fibr\'es $\,\L(\xi^1),\ldots,\L(\xi^r)$ sont en somme directe comme 
sous-fibr\'es de $N_CX$. On obtient que $N_CX = \oplus_{k=1}^r \L(\xi^k)$ est une
 somme directe de fibr\'es en droites de degr\'e $n-1$ et que $E=H^0(C,N_CX)$.
Ceci d\'emontre la deuxi\`eme partie de l'\'enonc\'e.

\sk
Si $1\leq p\leq n-1$, exactement le m\^eme argument qu'au d\'ebut de la d\'emonstration, appliqu\'e 
\`a la restriction de $\g$ \`a $\{(a_1,\ldots,a_p)\}\times X^{(n-p)}_{\rm adm}(a_1,\ldots,a_p)$,
montre que la vari\'et\'e $\Si_q(X;a_1,\ldots,a_p)$ est lisse et que son espace tangent
en $C$ est donn\'e par 
$$
T_C \Si_q(X;a_1,\ldots,a_p) = \oplus_{i=p+1}^n E_i.
$$
Il est de dimension $(n-p)r$. Comme $N_CX$ est une somme directe de fibr\'es en droites 
de degr\'e $n-1$, l'espace de ses sections qui s'annulent en $a_1,\ldots,a_p$
est de dimension $(n-p)r$ et donc $T_C \Si_q(X;a_1,\ldots,a_p)$ est \'egal
\`a cet espace.
\epf

Finissons par une remarque. 
On obtient des \og cartes locales \fg\,  de $\Si_q(X)$ 
au  voisinage de la courbe $C=\g(\wh{\bs a})$ en choisissant,
pour $i=1,\ldots,n,$ un germe d'hypersurface lisse $Y_i\subset X$, transverse \`a $C$ en $a_i$.
La restriction de $\g$,
$$
\g_{|Y_1\times \cdots \times Y_n}  : \; \sqcap_{i=1}^n Y_i \rightarrow \Si_q(X) 
$$
est, compte tenu de ce qui pr\'ec\`ede,  un diff\'eomorphisme local. De m\^eme,
si $1\leq p\leq n-1$, sa restriction $\,\sqcap_{i=1}^p\{a_i\} \times \sqcap_{i=p+1}^n Y_i \rightarrow \Si_q(X;a_1,\ldots,a_p)$
est un diff\'eomorphisme local.

\ssct{Projections osculatrices g\'en\'erales}

On d\'emontre l'analogue de la Proposition \ref{proj} pour une projection
osculatrice associ\'ee \`a un point admissible pond\'er\'e de $X$.
Par composition de telles projections, le r\'esultat couvre en fait
le cas d'une projection osculatrice g\'en\'erale. On se donne une pond\'eration
$$
(\rho_1,\ldots,\rho_{n-1}), \;\;\; \text{avec} \;\; \rho_{n-1}\geq 1.
$$
\bpr
\label{proj-gene}
On suppose $n\geq 3$. Soit $a\in X_{\rm adm}$
et $\t_a$ une projection osculatrice de centre $X_a(\rho_1)$. La 
vari\'et\'e 
$\t_a(X)$ appartient \`a la classe $\cal{X}_{r+1,n-1}(q-(\rho_1+1))$
et $\t_a : X\dasharrow \t_a(X)$ est birationnelle et induit un diff\'eomorphisme de $X_{\rm adm}(a)$ 
sur son image, contenue dans $\t_a(X)_{\rm adm}$. 
Si ${\boldsymbol b}\in X^{(n-1)}_{\rm adm}(a)$, alors $\t_a({\boldsymbol b})\in \t_a(X)^{(n-1)}_{\rm adm}$.
\epr
\bpf
Soit $\,(a,{\boldsymbol b})$ un $n$-uplet admissible avec ${\boldsymbol b}=(b_2,\ldots,b_n)$.
On peut choisir la cible de la projection $\t_a$, soit :
$$
\t_a: \, \P^N\dasharrow \oplus_{i=2}^{n-1} X_{b_i}(\rho_i).
$$
On note $\t=\t_a$, $\,X'=\t(X)$ et  $q'=q-(\rho_1+1)$.

\sk
Le d\'ebut est une r\'ep\'etition de la d\'emonstration de la Proposition \ref{proj}.
\`A d\'efaut d'\^etre clair, on essaiera d'\^etre bref.

\sk
L'image $X'$ de $X$ engendre un espace de dimension $\pi_{r,n-1}(q')$ et pour 
$i=2,\ldots,n$, l'application $\t: X\dasharrow X'$ induit un diff\'eomorphisme local 
du germe de $X$ en $b_i$ sur un germe lisse et $\rho'$-r\'egulier ${X}'_i\subset X'$ en $b_i$,
de dimension $r+1$. On a choisi la cible de telle sorte que $(X'_i)_{b_i}(\rho_i)=X_{b_i}(\rho_i)$.

\sk
Si un hyperplan $H$ contient $X_a(\rho_1)$, une courbe $C\in \Si_q(X;a)$ coupe $H$ au moins \`a l'ordre $(\rho_1+1)$ en $a$.
Le degr\'e de $C'=\t(C)$ est donc $\leq q'$. D'autre part, pour 
tout $i\in \{2,\ldots,n-1\}$, l'image $C'_i$ du germe de $C$ en $b_i$
est un germe lisse donc $\rho'$-r\'egulier en $b_i$ et les 
osculateurs $(C'_i)_{b_i}(\rho_i)$ sont en somme directe projective.
Il en r\'esulte que $C'$
engendre un espace de dimension $\geq q'$. Finalement $C'\in {\rm CRN}_{q'}(X')$.

\sk
En faisant varier $\boldsymbol b$, on obtient que $X'$ appartient \`a la classe $\cal{X}_{r+1,n-1}(q')$.

\sk
Pour montrer que, si ${\bs b}\in X_{\rm adm}^{(n-1)}(a)$, alors 
$\t({\bs b})={\bs b}$ est un $(n-1)$-uplet admissible de $X'$, 
il suffit, voir la D\'efinition \ref{D3.1}, 
de montrer que $\boldsymbol b$ est contenu dans $X'_{\rm reg}$.
On utilise les r\'esultats obtenus sur les 
projections osculatrices associ\'ees \`a un $(n-2)$-uplet admissible
pour le v\'erifier.
C'est la raison pour laquelle on n'a pas trait\'e d'embl\'ee
les projections osculatrices g\'en\'erales.

\sk
Il suffit de montrer que $b_{n-1}\in X'_{\rm reg}$.
Soit $\phi$ la projection $\oplus_{i=2}^{n-1} X_{b_i}(\rho_i) \dasharrow X_{b_{n-1}}(\rho_{n-1})$
de centre $\oplus_{i=2}^{n-2} X_{b_i}(\rho_i)$. La compos\'ee 
$\wh{\t} = \phi\circ \t$ est une projection osculatrice qui envoie 
$X$ sur une vari\'et\'e 
de Veron\`ese $X''$ d'ordre $\rho_{n-1}$.
Comme $\wh{\t}: X\dasharrow X''$ est birationnelle d'apr\`es la 
Proposition \ref{proj}, il en va de m\^eme
de $\t: X\dasharrow X'$ et de $\phi: X'\dasharrow X''$.
De m\^eme, comme $\wh{\t}$ induit un  diff\'eomorphisme d'un voisinage 
de $b_{n-1}$ dans $X$ sur un voisinage de $b_{n-1}$ dans $X''$, 
le germe $X'_{n-1}$ est le germe de $X'$ en $b_{n-1}$. Donc $b_{n-1}$ appartient \`a $X'_{\rm reg}$.
\epf
Il est clair, compte tenu en particulier de la derni\`ere partie de 
l'\'enonc\'e, qu'on obtient par it\'eration un \'enonc\'e analogue 
pour une projection osculatrice associ\'ee \`a un $p$-uplet 
admissible, quel que soit $p\in\{1,\ldots,n-2\}$.

\sct{Int\'egrabilit\'e}

\ssct{Introduction}

Dans cette section, nous terminons la d\'emonstration du Th\'eo\-r\`eme \ref{Th4}.
Pr\'esentons la strat\'egie suivie.

\sk
Si $X_0\in \cal{X}_{r+1,n}(n-1)$,
une section hyperplane 
g\'en\'erique de $X_0$ dans $\P^{r+n-1}$ est une vari\'et\'e minimale
de m\^eme degr\'e. Le syst\`eme $\cal{Y}(X_0)$ des sections hyperplanes de $X_0$ est donc 
un syst\`eme lin\'eaire de dimension $r+n-1$, dont 
l'\'el\'ement g\'en\'erique $Y$ appartient \`a la classe $\cal{X}_{r,n}(n-1)$.
Compte tenu de la propri\'et\'e fondamentale d'une application $\phi : \; X\dasharrow X_0$
qui associe une vari\'et\'e {\em standard} $X\in \cal{X}_{r+1,n}(q)$  \`a $X_0$ au sens de la D\'efinition \ref{D2},
le syst\`eme lin\'eaire $\phi^\star(\cal{Y}(X_0))$ (comme $X$ n'est pas suppos\'ee lisse, on pr\'ecisera ce qu'on entend par l\`a)
est 
de dimension $r+n-1$ et son \'el\'ement g\'en\'erique appartient \`a la classe $\cal{X}_{r,n}(q)$.
En g\'en\'eral, le probl\`eme est de d\'efinir, dans les cas favorables,
un syst\`eme lin\'eaire $\cal{Y}(X)$ sur $X$ qui v\'erifie ces propri\'et\'es.

\sk
Si $n=2$, $\,X$ est une vari\'et\'e de Veronese d'ordre $q$
et $\cal{Y}(X)$ est le syst\`eme des vari\'et\'es de Veronese $Y\subset X$, de dimension $r$ et d'ordre $q$.
C'est aussi l'image du syst\`eme des hyperplans de $X_0=\P^{r+1}$
par un isomorphisme de Veronese $\,\P^{r+1}\rightarrow X$.

En g\'en\'eral, la solution est donn\'ee par les projections osculatrices. 
En effet, soit $C$ une courbe admissible, ${\boldsymbol a}\subset C$
un $(n-2)$-uplet et $\t_{\boldsymbol a}$ une application birationnelle associ\'ee de $X$ 
sur une vari\'et\'e de Veronese $\t_{\boldsymbol a}(X)$ d'ordre $\rho$.
Si $Y\in \cal{X}_{r,n}(q)$ contient $C$, son image $\t_{\boldsymbol a}(Y)$
appartient \`a $\cal{X}_{r,2}(\rho)$
pour la m\^eme raison que $\t_{\boldsymbol a}(X)$ appartient \`a $\cal{X}_{r+1,2}(\rho)$.

\sk
Ceci sugg\`ere de consid\'erer, pour tout ${\boldsymbol a}\in X^{(n-2)}_{\rm adm}$,
le syst\`eme lin\'eaire $\t_{\boldsymbol a}^\star(\cal{Y}(\t_{\boldsymbol a}(X)))$. 
On introduira une notion provisoire d'\og int\'egrabilit\'e \fg\,
de $X$ dont on montrera qu'elle est \'equivalente 
au fait que les syst\`emes $\t_{\boldsymbol a}^\star(\cal{Y}(\t_{\boldsymbol a}(X)))$
appartiennent \`a un m\^eme syst\`eme lin\'eaire, puis au fait 
que la vari\'et\'e $X$ est standard.

Un argument g\'eom\'etrique tr\`es  simple montrera que $X$ est int\'egrable si 
$\rho\geq m$
dans la division euclidienne $q=\rho(n-1)+m-1$ de $q$ par $n-1$. En particulier, c'est le cas si 
$X\in \cal{X}_{r+1,3}(q)$ et $q\neq 3$. 

\sk
Le cas g\'en\'eral sera r\'esolu gr\^ace \`a un r\'esultat 
d'h\'er\'edit\'e de l'int\'egrabilit\'e, qui a son int\'er\^et propre : 

\sk
{\em si $\,n\geq 4$, la vari\'et\'e $X\in \cal{X}_{r+1,n}(q)$ est standard 
si et seulement si $\t_a(X)$ est un \'el\'ement standard de $\cal{X}_{r+1,n-1}(q-(\rho_a +1))$ 
pour tout point admissible pond\'er\'e $a\in X$.}

\sk
Cet \'enonc\'e est faux si n=3.

\ssct{Int\'egrabilit\'e}

Les hypoth\`eses  g\'en\'erales et les notations sont les m\^emes 
que dans le Chapitre 2. On suppose toujours $n\geq 3$.
On utilisera le Th\'eor\`eme~\ref{Xlisse} et selon les 
cas, par commodit\'e ou par hasard, la Proposition~\ref{proj} ou la Proposititon~\ref{proj-jumelle}.

\sk
On notera $\pi_{\boldsymbol a} : X\dasharrow \P^{r+1}$ une application 
birationnelle associ\'ee \`a un $(n-2)$-uplet admissible $\boldsymbol a$ de $X$
et $\cal{H}$  le syst\`eme lin\'eaire des hyperplans de $\P^{r+1}$.

On notera $\t_{\boldsymbol a} : X\dasharrow \t_{\boldsymbol a}(X)$ une application 
birationnelle sur une vari\'et\'e de Veronese d'ordre $\rho$
associ\'ee \`a ${\boldsymbol a}\in X^{(n-2)}_{\rm adm}$
et souvent, $\cal{H}_{\boldsymbol a}$ au lieu de $\cal{Y}(\t_{\boldsymbol a}(X))$,  le syst\`eme lin\'eaire 
des diviseurs de Veronese d'ordre $\rho$ de $\t_{\boldsymbol a}(X)$.
\bd
\label{D4.1}
Un {\em diviseur admissible} de $X$ est une sous-vari\'et\'e $Y\subset X$ irr\'eductible 
de dimension $r$ qui poss\`ede la propri\'et\'e suivante : il existe une courbe admissible 
$C\subset Y$ et un $(n-2)$-uplet ${\boldsymbol a}$ de $C$ tels que 
$\pi_{\boldsymbol a}(Y)\in \cal{H}$.
\ed
Les diviseurs admissibles {\em g\'en\'eriques} seront, dans les cas favorables, des g\'en\'erateurs 
du syst\`eme lin\'eaire qu'on va construire.
On les obtient tous de la mani\`ere suivante. 

\sk
Soit $C$ une courbe admissible, ${\boldsymbol a}\subset C$ un $(n-2)$-uplet et $x_0$ un point 
de $C\bck {\boldsymbol a}$. Si $h\subset T_{x_0}X$ est un hyperplan qui contient 
$T_{x_0}C$ et si $H$ est l'hyperplan de $\P^{r+1}$ passant par $\pi(x_0)$ et de direction 
l'image de $h$ dans $T_{\pi(x_0)}\P^{r+1}$, il est clair que $\pi_{\boldsymbol a}^{-1}(H)$\footnote{
La notation $\pi_{\bs a}^{-1}(H)$ correspond \`a l'image de $H$
par $\pi_{\bs a}^{-1}: \, \P^{r+1}\dasharrow X$, voir la note en bas de page 
signal\'ee au d\'ebut de la d\'emonstration de la Proposition \ref{proj}.}
est l'unique diviseur admissible $Y$ qui contient $C$ et est tel que $\pi_{\boldsymbol a}(Y)=H$.

Par construction, si $\wh{\boldsymbol x}\in Y^n$ est un $n$-uplet admissible
qui contient $\boldsymbol a$, la courbe $\g(\wh{\boldsymbol x})$ est contenue dans $Y$,
mais rien ne dit que $Y$ contienne d'autres courbes admissibles que 
celles-ci.
\bd
\label{D4.2}
Si $n\geq 3$, on dit que $X\in \cal{X}_{r+1,n}(q)$ est {\em int\'egrable} (en tant qu'\'el\'ement de 
la classe $\cal{X}_{r+1,n}(q)$) si tout diviseur admissible 
de $X$ appartient \`a la classe $\cal{X}_{r,n}(q)$.
On convient que tout $X\in \cal{X}_{r+1,2}(q)$ est int\'egrable.
\ed
Le r\'esultat suivant donne des caract\'erisations de l'int\'egrabilit\'e qu'on utilisera :
\bpr
\label{equivalences}
Soit $Y$ une sous-vari\'et\'e irr\'eductible de $X$, de dimension $r$,
contenant au moins une courbe admissible. Elle engendre un espace de 
dimension plus grande que $\pi_{r-1,n}(q)$ et les propri\'et\'es suivantes sont \'equi\-valentes.
\be

\item La vari\'et\'e $Y$ appartient \`a la classe $\cal{X}_{r,n}(q)$.

\item On a $Y\cdot C=n-1$ pour une ou toute courbe admissible $C\not\subset Y$.

\item On a $Y\cdot C\leq n-1$ pour une ou toute courbe admissible $C\not\subset Y$.

\item Pour toute courbe admissible $C\subset Y$ 
et tout $(n-2)$-uplet ${\bs a}\subset C$, on a $\pi_{\boldsymbol a}(Y)\in \cal{H}$.

\item Le diviseur $Y$ est admissible et, pour toute courbe admissible $C\subset Y$ et pour  tout $(n-2)$-uplet
${\boldsymbol a}=(a_1,\ldots,a_{n-2})\subset C$, si $b_1\in C\bck {\boldsymbol a}$ et $\,{\boldsymbol b}=(b_1,a_2,\ldots,a_{n-2})$,
on a l'implication :
\beq
\label{impli}
\pi_{\boldsymbol a}(Y)\in \cal{H}  \; \Rightarrow \pi_{{\boldsymbol b}}(Y)\in \cal{H}.
\eeq
\ee
\epr 
\bpf
Soit $Y$ une sous-vari\'et\'e irr\'eductible de $X$, de dimension $r$
et contenant au moins une courbe admissible. Comme toute courbe $C\in \Si_q(X)$ est contenue 
dans $X_{\rm reg}$ le nombre d'intersection $Y \cdot C$ est bien d\'efini si $C\not\subset Y$.
D'autre part, si $C$ est contenue dans $Y$
et si $\wh{\boldsymbol a}\subset C$ est un $n$-uplet, un $n$-uplet $\wh{\boldsymbol b}\in Y_{\rm reg}^n$
assez voisin de $\wh{\bs a}$ est admissible. Les osculateurs $Y_{b_1}(\rho),\ldots,Y_{b_m}(\rho),
Y_{b_{m+1}}(\rho-1),\ldots,Y_{b_{n-1}}(\rho-1)$ sont alors de dimension maximale 
et, comme sous-espaces d'espaces en somme directe projective, ils le sont aussi.
Donc $Y$ engendre un espace de dimension $\geq \pi_{r-1,n}(q)$.

\bk
Supposons $Y\in \cal{X}_{r,n}(q)$. Si ${\wh{\boldsymbol x}}$ est un $n$-uplet 
admissible de $X$ contenu dans $Y$,
la courbe $\g({\wh{\boldsymbol x}})$
est le {\em seul} \'el\'ement de ${\rm CR}_q(X)$ qui contient $\wh{\boldsymbol x}$.
Par hypoth\`ese, il existe un \'el\'ement de ${\rm CR}_q(Y)$
qui contient $\wh{\bs x}$. C'est n\'ecessairement $\g(\wh{\bs x})$.
Il en r\'esulte qu'une courbe admissible qui a $n$ points deux-\`a-deux 
distincts dans $Y$ est contenue dans $Y$.

Par hypoth\`ese, $Y$ contient une courbe admissible $C$. 
Soit ${\boldsymbol a}\subset C$ un $(n-2)$-uplet et 
${\boldsymbol b}=(b_1,b_2)\subset C\bck {\boldsymbol a}$ un $2$-uplet.
Soit ${\boldsymbol x}\in X^2$ voisin de $\boldsymbol b$
et $C({\bs x})=\g({\bs a},{\bs x})$.
Si ${\bs x}\subset Y$, la courbe $C({\bs x})$ a au moins
$n$ points distincts dans $Y$ donc est contenue dans $Y$.
Il en r\'esulte que, pour tout $2$-uplet ${\boldsymbol x}'\subset \pi_{\boldsymbol a}(Y)$
voisin de $\pi_{\boldsymbol a}({\boldsymbol b})$, la droite qui contient ${\boldsymbol x}'$
est contenue dans $\pi_{\boldsymbol a}(Y)$. Autrement dit $\pi_{\boldsymbol a}(Y)$
est un hyperplan $H$ de $\P^{r+1}$.  

On choisit maintenant ${\bs x}=(b_1,x)$ avec $x\notin Y$.
La courbe $C(b_1,x)$ n'est pas contenue dans $Y$ donc coupe $Y$ 
en $a_1,\ldots,a_{n-1},b_1$ et nulle part ailleurs.
Son image par $\pi_{\bs a}$ est une droite $l\not\subset H$ :
elle n'est pas tangente \`a $H$ en $\pi_{\bs a}(b_1)$.
Ainsi, la courbe $C(b_1,x)$ coupe $Y$ transversalement en 
$b_1$. Par sym\'etrie en $a_1,\ldots,a_{n-1}$ et $b_1$,
elle coupe $Y$ transversalement en tous ces  points et 
donc $Y\cdot C(b_1,x)=n-1$. 
Ceci montre que (1) implique (2). \'Evidemment (2) implique (3).

\bk
Supposons maintenant qu'on a $Y\cdot C\leq n-1$ pour toute courbe admissible
$C\not\subset Y$. En particulier, une courbe admissible qui a $n$ points deux-\`a-deux 
distincts dans $Y$ est contenue dans $Y$. Comme c'est la seule propri\'et\'e
utilis\'ee, dans le paragraphe pr\'ec\'edent,
pour montrer que $\pi_{\bs a}(Y)$ est un hyperplan pour tout $(n-2)$-uplet
$\bs a$ contenu dans une courbe admissible contenue dans $Y$,
on obtient que (3) implique (4).

\bk
Il est clair que (4) implique (5). 
On suppose enfin que $Y$ v\'erifie (5). 
Notons $\boldsymbol A$ l'ensemble des $(n-2)$-uplets $\boldsymbol a$ contenus
dans au moins une courbe admissible contenue dans $Y$ 
et tels que $\pi_{\boldsymbol a}(Y)\in \cal{H}$.

\sk
Par hypoth\`ese $Y$  est un diviseur admissible, donc $\boldsymbol A$ est non vide.
On fixe ${\boldsymbol a}\in {\boldsymbol A}$ et une courbe admissible $C\subset Y$ 
qui contient $\boldsymbol a$. On sait que $\pi_{\bs a}$ induit un diff\'eomorphisme
d'un voisinage de  $C\bck {\boldsymbol a}$ sur son image. En particulier,
$Y$ est lisse au voisinage de $C\bck {\bs a}$.

\sk
Si ${\boldsymbol b}\subset C$ est un $(n-2)$-uplet disjoint de $\boldsymbol a$,
l'hypoth\`ese (5), appliqu\'ee $(n-2)$ fois donne que ${\boldsymbol b}$
appartient \`a $\boldsymbol A$. En particulier 
$Y$ est lisse au voisinage de $C\bck {\bs b}$ donc au voisinage de $C$.
Si ${\bs b}\subset C$ n'est pas disjoint de $\bs a$,
on obtient encore que ${\bs b}$ appartient \`a ${\bs A}$ 
en passant par l'interm\'ediaire d'un $(n-2)$-uplet ${\bs c}\subset C$,
disjoint \`a la fois de $\bs a$ et de $\bs b$.
En r\'esum\'e, {\em si ${\bs a}\in {\bs A}$, tout 
$(n-2)$-uplet d'une courbe admissible contenue 
dans $Y$ qui contient $\bs a$ appartient \`a $\bs A$ et,
d'autre part, on sait que toute courbe admissible qui 
contient un $n$-uplet admissible $\wh{\bs y}\subset Y$
qui contient $\bs a$ est contenue dans $Y$}.

\sk
Soit $\wh{\bs a}=(a_1,\ldots,a_n)$ 
un $n$-uplet de $C$
et  $\wh{\bs y}=(y_1,\ldots,y_n)\in Y^n$ voisin de $\wh{\bs a}$.
On passe de $\wh{\bs a}$ \`a $\wh{\bs y}$ par l'interm\'ediaire des $n$-uplets
admissibles 
$$
\wh{\bs y}_p = (y_1,\ldots,y_p,a_{p+1},\ldots,a_n)\in Y^n, \qquad p=0,\ldots,n.
$$
La courbe $\g(\wh{\bs y}_0)=C$ est contenue dans 
$Y$ et $(a_1,\ldots,a_{n-2})\in {\bs A}$ donc aussi tout 
$(n-2)$-uplet de $\g(\wh{\bs y}_0)$. En particulier $(a_2,\ldots,a_{n-1})\in {\bs A}$
et donc la courbe $\g(\wh{\bs y}_1)$, qui contient ce $(n-2)$-uplet, est contenue dans 
$Y$. Comme elle contient un \'el\'ement de $\bs A$, on peut it\'erer
le raisonnement. On obtient finalement que $\g(\wh{\bs y})$
est contenu dans $Y$ et donc que $Y$ appartient \`a la 
classe $\cal{X}_{r,n}(q)$. Ceci montre que (5) implique (1)
et termine la d\'emonstration.
\epf

\ssct{Premiers exemples d'int\'egrabilit\'e} 

Le cas $r=1$ est rapidement r\'egl\'e. En effet un diviseur admissible n'est alors rien d'autre 
qu'une courbe admissible :
\bpr
\label{P4.1}
Toute surface appartenant \`a la classe $\cal{X}_{2,n}(q)$ est int\'egrable.
\epr
On suppose maintenant  $r\geq 2$. Le cas particulier suivant, dont la d\'emonstration 
est tr\`es simple,  est crucial. Il a comme corollaire imm\'ediat
qu'une vari\'et\'e $X\in \cal{X}_{r+1,3}(q)$ est int\'egrable si $q\neq 2n-3$
et comme cons\'equence indirecte,  comme on le 
verra dans la prochaine section, 
qu'une vari\'et\'e $X\in \cal{X}_{r+1,n}(q)$ est int\'egrable si $q\neq 2n-3$,
ce qui est une version pr\'eliminaire du Th\'eor\`eme \ref{Th4}.
\bpr
\label{P4.2}
Soit $X$ une vari\'et\'e de la classe $\cal{X}_{r+1,n}(q)$. Si  $\rho\geq m$
dans la division euclidienne $q=\rho(n-1)+m-1$ de $q$ par $n-1$, alors $X$  est int\'egrable.
\epr
\bpf
Soit $Y$ un diviseur admissible de $X$. On peut introduire  une courbe admissible $C\subset Y$,
un $(n-2)$-uplet $\bs a$ de $C$ et $\t_{\bs a}: X \dasharrow \t_{\bs a}(X)$,
une application birationnelle sur une vari\'et\'e de Veronese d'ordre $\rho$,
telle que $\t_{\bs a}(Y)\in \cal{H}_{\bs a}$.

\sk
Fixons un $n$-uplet $\wh{b}=(b_1,\ldots,b_n)\subset C\bck{\bs a}$.
Si ${\bs{\wh{y}}}\in Y^n$ est assez voisin de $\bs{\wh{b}}$, la courbe $\g(\bs{\wh{y}})$ 
n'est pas contenue dans le centre de la projection $\t_{\bs a}$. On peut donc consid\'erer
son image  $\t_{\bs a}(\g({\bs{\wh{y}}}))$. 
C'est une courbe de degr\'e $q'\leq q$, ce qu'on \'ecrit sous la forme 
$$
q' \leq \rho(n-1) + m-1.
$$
D'autre part son degr\'e $q'$ est un multiple de $\rho$ car c'est l'image 
d'une courbe de $\P^{r+1}$ par un plongement de Veronese d'ordre $\rho$.

\sk
On suppose  maintenant $\rho \geq m$. Alors le degr\'e de la courbe 
$\t_{\bs a}(\g({\bs{\wh{y}}}))$ est $\leq \rho(n-1)$. 

\sk
Comme $\wh{\bs b}\subset C$ est disjoint de $\bs a$ et que $\t_{\bs a}$
induit un diff\'eomorphisme de $C\bck {\bs a}$ sur son image,
pour $\wh{\bs y}$ assez voisin de $\wh{\bs b}$, 
la courbe $\t_{\bs a}(\g(\wh{\bs y}))$ a $n$ points diff\'erents dans $\t_{\bs a}(Y)$.
Elle est donc contenue dans $\t_{\bs a}(Y)$,
puisque c'est 
l'image, par un plongement de Veronese $v: \P^{r+1}\rightarrow \t_{\bs a}(X)$ d'ordre $\rho$,
d'une courbe de degr\'e
$\leq n-1$ qui a $n$ points diff\'erents dans l'hyperplan $v^{-1}(\t_{\bs a}(Y))$.
La courbe $\g({\bs{\wh{y}}})$ est contenue dans $Y$.

Ainsi le diviseur admissible $Y$ engendre un espace de dimension $\geq \pi_{r-1,n}(q)$
et, pour tout ${\bs{\wh{y}}}\in Y^n$ voisin de $\bs{\wh{b}}$, il contient
la courbe $\g(\wh{\bs y})$.
Il appartient \`a la classe $\cal{X}_{r,n}(q)$.

Par d\'efinition, la vari\'et\'e $X$ est int\'egrable.
\epf
On en d\'eduit une version pr\'eliminaire du Th\'eor\`eme \ref{Th4} dans le cas $n=3$.
\bpr
\label{P4.3}
Toute vari\'et\'e $X\in \cal{X}_{r+1,3}(q)$ est int\'egrable si $q\neq 3$.
\epr
\bpf
En effet, si $q=2$, $X$ est une vari\'et\'e minimale et si $q\geq 4$,
alors $q=2\rho + m-1$ avec $m\in \{1,2\}$ et $\rho\geq 2$.
\epf

\ssct{H\'er\'edit\'e de l'int\'egrabilit\'e}

Dans ce paragraphe, on se donne une pond\'eration
$$
(\rho_1,\ldots,\rho_{n-1}), \;\;\; \text{avec} \;\; \rho_{n-1}\geq 1.
$$
Rappelons qu'auparavant, on a souvent choisi la pond\'eration (\ref{ponderation})
seulement pour fixer les id\'ees. On utilisera la Proposition \ref{proj-gene}.

\sk
Supposons $n\geq 3$ et $q\geq n$. Soit $X \in \cal{X}_{r+1,n}(q)$ une vari\'et\'e int\'egrable, 
$\,a\in X_{\rm adm}$ un point affect\'e du poids $\rho_1$ et $q'=q-(\rho_1+1)$. 
Il r\'esulte facilement de la Proposition \ref{proj-gene},
qu'un diviseur admissible {\em g\'en\'erique} de la vari\'et\'e 
$\t_a(X)$ appartient \`a la classe 
$\cal{X}_{r,n-1}(q')$. La restriction de g\'en\'ericit\'e emp\^eche 
d'en d\'eduire que  l'int\'egrabilit\'e est stable par projection osculatrice,
 \`a cause du choix qu'on fait 
d'une d\'efinition \og forte \fg\, de l'int\'egrabilit\'e.
On verra que c'est vrai \`a la fin de ce chapitre.
Quoi qu'il en soit, cette propri\'et\'e est banale. Ce n'est pas le cas de la r\'eciproque partielle suivante.
\bpr
\label{heredite1}
Soit $n\geq 4$, $\,q\geq n$  et $X\in \cal{X}_{r+1,n}(q)$. Si pour tout point $c_1\in X_{\rm adm}$
affect\'e du poids $\rho_1$,  la vari\'et\'e $\t_{c_1}(X)$
est int\'egrable, alors $X$ est int\'egrable.
\epr
Ce r\'esultat est important. Il est faux si $n=3$. On le r\'e\'enoncera 
en termes de vari\'et\'es standards \`a la fin de ce chapitre.
Bien s\^ur, la d\'emonstration ci-dessous n'utilise pas la stabilit\'e 
de l'int\'egrabilit\'e par projection osculatrice.
\bpf
On utilise la propri\'et\'e (5) de la Proposition \ref{equivalences} comme 
crit\`ere d'int\'e\-grabilit\'e.
On se donne un diviseur admissible $Y$ de $X$, une courbe admissible $C\subset Y$
et deux $(n-2)$-uplets de $C$ de la forme 
$$
{\bs a} = (c_1,\ldots,c_{n-3},a), \qquad {\bs b}= (c_1,\ldots,c_{n-3},b).
$$
Il s'agit de montrer qu'on a l'implication : 
$\t_{\bs a}(Y)\in \cal{H}_{\bs a} \Rightarrow \t_{\bs b}(Y)\in \cal{H}_{\bs b}$.

\bk
Notons ${\bs c}= (c_1,{\bs c}') = (c_1,\ldots,c_{n-3})$. Pour fixer les id\'ees, on se donne 
$d,e\in C\bck({\bs a}\cup {\bs b})$ avec $d\neq e$ et on \'ecrit les d\'ecompositions en sommes directes projectives 
$$
\P^N = X_{c_1}(\rho_1)\oplus \Q,
\;\;\; 
\Q=(\oplus_{i=2}^{n-3} X_{c_i}(\rho_i))\oplus \Q',
\;\;\;  
\Q' = X_d(\rho_{n-2})\oplus X_e(\rho_{n-1}),
$$
qu'on utilise pour choisir les cibles des projections osculatrices  
$$
\t_{c_1}: \P^N\dasharrow \Q, \;\;\;  \t_{\bs c}: \P^N \dasharrow \Q', \;\;\; 
\t_{\bs a} : \P^N \dasharrow X_e(\rho_{n-1}), \;\;\; \t_{\bs b} : \P^N \dasharrow X_e(\rho_{n-1}).
$$
Soit $\phi_{{\bs c}'} : \Q \dasharrow \Q'$ la projection de centre $\oplus_{i=2}^{n-3} X_{c_i}(\rho_i)$.
(C'est l'identit\'e si $n=4$.) On a  
$$
\t_{\bs c}=\phi_{{\bs c}'}\circ \t_{c_1}.
$$

\sk
Compte tenu de la Proposition \ref{proj-gene},  la vari\'et\'e $\t_{\bs c}(X)$
appartient \`a la classe $\cal{X}_{r+1,3}(q')$ avec $q'=\rho_{n-2}+\rho_{n-1}+1$ et de plus, $(\t_{\bs c}(a),e)$ et $(\t_{\bs c}(b),e)$
sont des paires admissibles de $\t_{\bs c}(X)$.
On leur associe les projections osculatrices $\s_a$ et $\s_b$ 
de $\Q'$ sur $X_e(\rho_{n-1})$,
de centres respectifs $\t_{\bs c}(X)_{\t_{\bs c}(a)}(\rho_{n-2})$ et $\t_{\bs c}(X)_{\t_{\bs c}(b)}(\rho_{n-2})$.
Par construction, on a : 
$$
\t_{\bs a} = \s_a \circ \t_{\bs c}= \s_a\circ \phi_{{\bs c}'} \circ \t_{c_1},
 \qquad
\t_{\bs b} = \s_b \circ \t_{\bs c}= \s_b\circ \phi_{{\bs c}'} \circ \t_{c_1}.
$$

\sk
Finalement, revenons au probl\`eme pos\'e au d\'ebut de la d\'emonstration.
On suppose $\t_{\bs a}(Y)\in \cal{H}_{\bs a}$. La vari\'et\'e $\t_{c_1}(Y)$
contient la courbe admissible $\t_{c_1}(C)$ et 
$$
(\s_a\circ \phi_{{\bs c}'})(\t_{c_1}(Y))= \t_{\bs a}(Y)\in \cal{H}_{\bs a},
$$
donc $\t_{c_1}(Y)$ est un diviseur admissible de $\t_{c_1}(X)$. 
Sous les hypoth\`eses de l'\'enonc\'e, la vari\'et\'e $\t_{c_1}(X)$
est int\'egrable donc $(\s_b\circ \phi_{{\bs c}'})(\t_{c_1}(Y))$
appartient \`a $\cal{H}_{\bs b}$. On a obtenu 
que $\t_{\bs b}(Y)$ appartient \`a $\cal{H}_{\bs b}$, ce qui termine la d\'emonstration.
\epf
\bco
\label{reduction}
Toute vari\'et\'e $X$ de la classe $\cal{X}_{r+1,n}(q)$ est int\'egrable si
$q\neq 2n-3$.
\eco
\bpf
On sait que le r\'esultat est vrai si $n=3$. On suppose $n\geq 4$ et que le r\'esultat
est vrai \`a l'ordre $n-1$. On \'ecrit la division euclidienne 
$q = \rho(n-1)+m-1$ avec,  par hypoth\`ese $(\rho,m)\neq (1,n-1)$. 

Soit $a$ un point admissible de $X$,
pond\'er\'e par $\s\in \{\rho-1,\rho\}$, avec $\s=\rho$ si $m=n-1$.
La vari\'et\'e $\t_a(X)$ appartient \`a $\cal{X}_{r+1,n-1}(q')$,
$$
q' = \rho(n-1) + m - 1- (\s+1) = \rho(n-2) + m - 1 - (\s + 1 - \rho).
$$
Soit aussi $q'  = \rho'(n-2) + m'-1$ la division euclidienne de $q'$ par $(n-2)$.

\sk
Si $m \geq 2$, on choisit $\s = \rho$. Alors   
$(\rho',m')=(\rho,m-1)\neq (1,n-2)$ et $\t_a(X)$ est int\'egrable par
hypoth\`ese de r\'ecurrence.

Si $m = 1$, on choisit $\s=\rho-1$. Alors $(\rho',m')=(\rho,1)$ et on a la m\^eme
conclusion. 

D'apr\`es la Proposition \ref{heredite1}, la vari\'et\'e $X$ 
est int\'egrable.
\epf

\ssct{Construction d'un syst\`eme lin\'eaire sur $X$}

On note $R(X)$ le corps des fonctions rationnelles sur $X$. 
On n'exclut pas {\em a priori}
que $X_{\rm sing}$ ait une  composante irr\'eductible de dimension $r$.
Si une fonction rationnelle sur $X$ s'annule sur une telle composante, la d\'efinition 
de son ordre d'annulation n'est pas \'el\'ementaire, voir Fulton~\cite{Fu}.
C'est une difficult\'e bien connue de la th\'eorie des 
syst\`emes lin\'eaires. La circonstance suivante permet d'\'echapper \`a 
cette difficult\'e.
\ble
\label{fond}
Si la restriction de $f\in R(X)$ \`a 
$X_{\rm reg}$ n'a pas de p\^ole (ou de z\'ero) sur une 
courbe admissible, $f$ est constante.
En particulier, une fonction $f\in R(X)\bck\{0\}$ est d\'etermin\'ee modulo $\C^*$
par le diviseur qu'elle d\'efinit sur $X_{\rm reg}$ par restriction.
\ele 
\bpf
Si $f$ n'a pas de p\^ole sur la courbe admissible $C$ et si $a\in C$,
$f$ est constante, \'egale \`a $f(a)$,  sur toute courbe admissible voisine de $C$
qui passe par $a$.
Comme ces courbes recouvrent un ouvert non vide de $X$,
$f$ est constante.
\epf
Compte tenu de cette remarque, en n\'egligeant les \'eventuelles composantes 
irr\'eductibles de dimension $r$ de $X_{\rm sing}$,
on dispose d'une th\'eorie des syst\`emes lin\'eaires sur $X$
(on dira plut\^ot sur $X_{\rm reg}$) tout \`a fait analogue 
\`a celle qu'on aurait si $X$ \'etait lisse. On fait quelques rappels et on renvoie le lecteur
\`a la pr\'esentation de Mumford \cite{Mu} pour toute pr\'ecision.

\sk
Si $f\in R(X)\bck \{0\}$, le diviseur $(f_{|X_{\rm reg}}) = (f_{|X_{\rm reg}})_0 - (f_{|X_{\rm reg}})_\infty$
des z\'eros et des p\^oles de $f$ dans $X_{\rm reg}$ est bien d\'efini.
Deux diviseurs $D$ et $D'$ de  $X_{\rm reg}$ sont {\em lin\'eairement \'equivalents},
ce qu'on note $D\sim_{X_{\rm reg}} D'$, s'il existe une fonction $f\in R(X)$
telle que $(f_{|X_{\rm reg}}) = D'-D$. 
\`A tout  diviseur {\em effectif} $\,Y_0$ de $X_{\rm reg}$, on  associe :
\be
\item l'ensemble $|Y_0|$ des diviseurs effectifs de $X_{\rm reg}$ lin\'eairement \'equivalents 
\`a $Y_0$, c'est {\em le syst\`eme lin\'eaire complet engendr\'e} par $Y_0$ ;
\item l'espace vectoriel $Q(Y_0)$ des $f\in R(X)$
telles que, ou bien $f=0$, ou bien le diviseur $(f_{|X_{\rm reg}}) + Y_0$ est effectif.
Si $f\neq 0$, il revient au m\^eme de dire que le diviseur $Y_0-(f_{|X_{\rm reg}})_\infty$
est effectif ou encore que $(f_{|X_{\rm reg}})=Y-Y_0$ o\`u $Y\in |Y_0|$.
\ee

\sk
Le point est que, compte tenu du Lemme \ref{fond}, si $Y\sim_{X_{\rm reg}} Y_0$,
il existe une fonction $f\in R(X)\bck \{0\}$, {\em unique modulo $\C^\star$}
telle que $(f_{|X_{\rm reg}}) = Y - Y_0$. On a donc, comme c'est le cas 
pour une vari\'et\'e lisse, une correspondance  biunivoque entre 
les droites vectorielles de $Q(Y_0)$ et les \'el\'ements de $|Y_0|$.

\sk
Par d\'efinition, un syst\`eme lin\'eaire $\cal{Y}$ sur $X_{\rm reg}$
est un ensemble de diviseurs effectifs de $X_{\rm reg}$ associ\'e
par une telle correspondance \`a un sous-espace non nul 
$F$ d'un espace $Q(Y_0)$. Si cet espace est de dimension finie, 
$\dim F - 1$ est la dimension du syst\`eme $\cal{Y}$.

\bk
Soit $\pi: X\dasharrow \P^{r+1}$ une application birationnelle 
associ\'ee \` a un $(n-2)$-uplet admissible de $X$. 
Si $H_0,H\in \cal{H}\,$ sont deux hyperplans distincts et si $u\in R(\P^{r+1})\bck\{0\}$
est une fonction rationnelle de diviseur $(u)=H-H_0$,
la	 fonction rationnelle non constante $u\circ \pi$ 
induit un diviseur non nul sur $X_{\rm reg}$,
toujours en vertu du Lemme \ref{fond}.

\sk
Compte tenu de cette remarque on peut, de la m\^eme fa\c{c}on par exemple 
que dans \cite{Mu},  associer \`a tout $H\in \cal{H}$
un diviseur {\em non nul} $\,(\pi_{|X_{\rm reg}})^\star (H)$ de $X_{\rm reg}$.
Ce diviseur est  aussi la somme (non vide) $\,\sum_{i=1}^m n_iW_i$, 
o\`u $W_i$ d\'ecrit la famille des 
sous-vari\'et\'es  irr\'eductibles de dimension $r$  de $X_{\rm reg}$
d'image $H$ et o\`u $n_i$ est le degr\'e de 
l'application induite $\t: \, W_i\dasharrow H$.

\sk
Ces diviseurs sont les \'el\'ements d'un syst\`eme lin\'eaire 
$\,(\pi_{|X_{\rm reg}})^\star(\cal{H})$. Si $H_0\in \cal{H}$,
c'est le sous-syst\`eme 
du syst\`eme lin\'eaire complet $|(\pi_{|X_{\rm reg}})^\star(H_0)|$ associ\'e 
au sous-espace de $Q((\pi_{|X_{\rm reg}})^\star(H_0))$ des fonctions
de la forme $u\circ \pi$, o\`u $u\in R(\P^{r+1})$ et, ou bien $u=0$,
ou bien $(u)_\infty=\emptyset$, ou bien $(u)_\infty = H_0$. 

\sk
{\em On notera abusivement $\pi^\star$ au lieu de $\pi_{|X_{\rm reg}}^\star$.}

\sk
Soit $\Gamma_\pi \subset X\times \P^{r+1}$ le graphe de l'application $\pi$. 
Si $K\subset X$ est un ensemble alg\'ebrique,
on note $\pi[K]$ le sous-ensemble alg\'ebrique $\{x'\in \P^{r+1}, \, \exists \,x \in K, \; (x,x')\in \Gamma_\pi\}$
de $\P^{r+1}$. 
Si $H\in \cal{H}$, on d\'efinit $\pi^{-1}[H]$ de fa\c{c}on analogue.
D'ailleurs, $\pi^{-1}[H] \cap X_{\rm reg}$ est aussi le support 
du diviseur $\pi^\star(H)$.

\sk
On note $\cal{H}(x)$ la sous-vari\'et\'e des $H\in \cal{H}$
qui passent  par le point $x\in \P^{r+1}$.  
On utilisera le lemme suivant :
\ble
\label{dense}
Soit $\,{\bs a}$ un $(n-2)$-uplet admissible et $b\in X_{\rm adm}({\bs a})$.
L'ensemble des $H\in \cal{H}$ qui contiennent l'image d'une courbe 
$C\in \Si_q(X;{\bs a},b)$ et tels que $\pi_{\bs a}^\star(H)$ soit irr\'eductible,
autrement dit que $\pi_{\bs a}^{-1}(H)=\pi_{\bs a}^\star(H)$, 
contient un ouvert dense de $\cal{H}(\pi_{\bs a}(b))$.
\ele
\bpf
Comme $\bs a$ est fix\'e, on ne le note pas en indice.
Notons $b'=\pi(b)$. On rappelle que $\pi$ induit un diff\'eomorphisme de $X_{\rm adm}({\bs a})$ 
sur son image dans $\P^{r+1}$. 
Comme $({\bs a},b)$ est admissible, l'image de $\Si_q(X;{\bs a},b)$  est 
un ouvert dense de $\Si_1(\P^{r+1};b')$. La premi\`ere propri\'et\'e 
d\'efinit donc un ouvert dense $\Omega$ de $\cal{H}(b')$.

\sk
Si $H\in \Omega$, le diviseur $\pi^\star(H)$
est de la forme $\pi^\star(H) = \pi^{-1}(H) + \sum_{i=1}^m n_iW_i$, 
o\`u les diviseurs irr\'eductibles 
$W_1,\ldots,W_m$ sont contenus dans $X\bck X_{\rm adm}({\bs a})$.
Pour $i=1,\ldots,m$, on peut choisir un point $w_i\in W_i$ tel que $\pi[w_i]$
soit fini : $\pi[w_i]=\{w'_{i,1},\ldots,w'_{i,m_i}\}$. 
D'autre part, $\pi$ est un diff\'eomorphisme d'un  voisinage de $b$
sur un voisinage de $b'$ et, compte tenu d'un cas particulier 
\'el\'ementaire d'un  th\'eor\`eme de Zariski,
la fibre $\pi^{-1}[b']$ est connexe donc r\'eduite au point $b$.
Il en r\'esulte que $b'$ n'est pas un des points $w'_{i,j}$.
Pour que $\pi^\star(H)$ soit irr\'eductible, il suffit 
que $H$ ne passe par  aucun des points $w'_{i,j}$, ce qui d\'efinit 
un ouvert dense de $\Omega$.
\epf 

On peut enfin \'enoncer :
\bpr
\label{systeme}
Soit $X$ une vari\'et\'e int\'egrable de la classe $\cal{X}_{r+1,n}(q)$.
La r\'eunion, pour $\bs a \in X^{(n-2)}_{\rm adm}$, des syst\`emes lin\'eaires 
$\cal{Y}(X;{\bs a})=\pi_{\bs a}^\star(\cal{H})$ est contenue 
dans un syst\`eme lin\'eaire $\cal{Y}(X)$ sur $X_{\rm reg}$,
uniquement d\'etermin\'e. Il est complet, de dimension $r+n-1$.

Un diviseur admissible {\em g\'en\'erique}
appartient \`a $\cal{Y}(X)$ et $\,Y\cdot C=n-1$ pour tout $Y\in \cal{Y}(X)$ et toute courbe admissible
$C\not\subset Y$.
\epr
\bpf
Le point important est de montrer que deux sys\-t\`emes $\cal{Y}(X;{\bs a})$
et $\cal{Y}(X;{\bs b})$ 
sont contenus dans un m\^eme syst\`eme lin\'eaire sur $X_{\rm reg}$, autrement dit qu'on a
\beq
\label{implication}
Y_1\in \cal{Y}(X;{\bs a}), \;\; Y_2\in \cal{Y}(X;{\bs b}), \;\; \Rightarrow \; Y_1\sim_{X_{\rm reg}} Y_2.
\eeq
En utilisant la transitivit\'e de la relation $\sim_{X_{\rm reg}}$ on se ram\`ene,
un peu comme dans la derni\`ere partie de la d\'emonstration de la Proposition \ref{equivalences},
\`a un cas particulier.

\sk
Comme $X_{\rm adm}({\bs a})$ et $X_{\rm adm}({\bs b})$ sont des ouverts denses 
de $X$, on peut choisir $c_1\in X$ tels que les $(n-1)$-uplets $(c_1,a_1,\ldots,a_{n-2})$ et $(c_1,b_1,\ldots,b_{n-2})$
soient admissibles et donc aussi les $(n-2)$-uplets 
$(c_1,a_2,\ldots,a_{n-2})$ et $(c_1,b_2,\ldots,b_{n-2})$.
De proche en proche, on construit un $(n-2)$-uplet admissible ${\bs c}=(c_1,\ldots,c_{n-2})$
tels que tous les $(n-2)$-uplets 
$$
(c_1,\ldots,c_p,a_{p+1},\ldots,a_{n-2}),\;\;\; (c_1,\ldots,c_p,b_{p+1},\ldots,b_{n-2}),\qquad  p=0,\ldots,n-2,
$$
soient admissibles. On passe donc de $\bs a$ \`a $\bs c$ puis de $\bs c$
\`a $\bs b$ par deux suites  
de $(n-2)$-uplets admissibles, tels que deux $(n-2)$-uplets cons\'ecutifs aient $n-3$
\'el\'ements communs, \`a la m\^eme position (ceci est important \`a cause de la 
pond\'eration).

\sk
On peut donc supposer ${\bs a}=(a_1,a_2,\ldots,a_{n-2})$
et ${\bs b}=(b_1,a_2,\ldots,a_{n-2})$. Finalement, soit $c\in X$
tel  que $(c,{\bs a})$ et $(c,{\bs b})$ soient 
admissibles. Deux $(n-2)$-uplets cons\'ecutifs
de la suite $\bs a$, $(c,a_2,\ldots,a_{n-2})$, $\bs b$ sont contenus 
dans une m\^eme courbe admissible.

\sk
On s'est ainsi ramen\'e \`a ne consid\'erer que le cas particulier 
de deux $(n-2)$-uplets admissibles ${\bs a}=(a_1,a_2,\ldots,a_{n-2})$
et ${\bs b}=(b_1,a_2,\ldots,a_{n-2})$ contenus dans une m\^eme courbe admissible.
On d\'emontre alors l'implication (\ref{implication})
en v\'erifiant que $\cal{Y}(X;{\bs a})\cap \cal{Y}(X;{\bs b})$ est non vide.

\sk
Notons $b'_1=\pi_{\bs a}(b_1)$ et $a'_1=\pi_{\bs b}(a_1)$.
Le lemme pr\'ec\'edent, appliqu\'e \`a $\bs a$ et \`a $\cal{H}(b'_1)$,
permet d'introduire un ouvert non vide $\cal{U}$ de $\cal{H}(b'_1)$ tel que
tout $H\in \cal{U}$  contient l'image $\pi_{\bs a}(C)$ d'une 
courbe $C\in \Si_q(X;{\bs a}\cup\{b_1\})$ et tel que $\pi_{\bs a}^{-1}(H)=\pi_{\bs a}^\star(H)$.

\sk
Comme $X$ est int\'egrable, $\pi_{\bs b}(\pi_{\bs a}^{-1}(H))\in \cal{H}(a'_1)$
pour tout $H\in \cal{U}$. On v\'erifie facilement (voir le commentaire 
qui suit la D\'efinition \ref{D4.1}) que l'application qui \`a $H\in \cal{H}(b'_1)$ associe 
$\pi_{\bs b}(\pi_{\bs a}^{-1}(H))\in \cal{H}(a'_1)$
est un hom\'eo\-mor\-phisme local.

\sk
On applique encore le Lemme \ref{dense}, cette fois \`a $\pi_{\bs b}$
et \`a $\cal{H}(a'_1)$.
Il existe un \'el\'ement $H_0$ de $(\pi_{\bs b}\circ \pi_{\bs a}^{-1})(\cal{U})$
tel que $\pi_{\bs b}^{-1}(H_0)=\pi_{\bs b}^\star(H_0)$.
Alors $\pi_{\bs b}^{-1}(H_0)$ est un \'el\'ement de $\cal{Y}(X;{\bs b})$
qui est aussi de la forme $\pi_{\bs a}^{-1}(H)$ avec $H\in \cal{U}$
donc  appartient au syst\`eme $\cal{Y}(X;{\bs a})$.

\bk
On a montr\'e que les syst\`emes lin\'eaires $\cal{Y}(X;{\bs a})$ 
sont contenus dans un m\^eme syst\`eme lin\'eaire $\cal{Y}(X)$
et qu'un diviseur admissible {\em g\'en\'erique} 
appartient au syst\`eme $\cal{Y}(X)$, voir le Lemme~\ref{dense}.
Par invariance lin\'eaire, on a bien $Y\cdot C=n-1$ pour tout $Y\in \cal{Y}(X)$
et toute courbe admissible $C\not\subset Y$.
Il reste \`a montrer que le syst\`eme $\cal{Y}(X)$ est n\'ecessairement de dimension $r+n-1$, c'est-\`a-dire que,
si ${\bs m}\in X^{r+n-1}$ est un $(r+n-1)$-uplet g\'en\'erique, il existe un et un seul 
diviseur $Y\in \cal{Y}(X)$ qui contient $\bs m$. En effet, si l'on d\'emontre cela,
comme $\cal{Y}(X)$ est contenu dans le syst\`eme lin\'eaire complet $\cal{Y}_1(X)$
d\'efini par un \'el\'ement quelconque de l'un des syst\`emes $\cal{Y}(X;{\bs a})$,
on aura n\'ecessairement $\cal{Y}_0(X)=\cal{Y}_1(X)$.

\bk
On \'ecrit ${\bs m}=({\bs a},{\bs d})=({\bs a},{\bs d}',{\bs d}'')$, la juxtaposition d'un  
$(n-2)$-uplet $\bs a$, d'un $2$-uplet ${\bs d}'$ et d'un $(r-1)$-uplet ${\bs d}''$. On peut supposer 
$({\bs a},{\bs d}')$ admissible, que les \'el\'ements de ${\bs d}$
appartiennent \`a $X_{\rm adm}({\bs a})$ et que les $r+1$ \'el\'ements de $\pi_{\bs a}({\bs d})$
sont en position g\'en\'erale dans $\P^{r+1}$, {\em i.e.}
engendrent un hyperplan $H$.

Notons $Y=\pi_{\bs a}^\star(H)$. C'est un \'el\'ement de $\cal{Y}(X;{\bs a})$ qui contient $\bs m$. 

\sk
R\'eciproquement, soit $Y_1\in \cal{Y}(X)$ un diviseur qui contient $\bs m$.
Il contient $({\bs a},{\bs d}')$ donc la courbe admissible $C_0=\g({\bs a},{\bs d}')$
puisque $Y_1\cdot C=n-1$ si $C\not\subset Y_1$ est une courbe admissible.
Soit $Y_0$ une composante irr\'eductible de $Y$ telle que $C_0\subset Y_0$. 
Si $C\not\subset Y_0$ est une courbe admissible, 
on a \'evidemment $Y_0\cdot C\leq n-1$, donc $Y_0\cdot C=n-1$ compte tenu de la Proposition \ref{equivalences}, 
qui donne aussi que $\pi_{\bs a}(Y_0)$ est un hyperplan $H_0\in \cal{H}$.

\sk
Si $Y_1=Y_0+W$, on a $W\cdot C=0$ si $C\not\subset Y_1$ est une courbe admissible, donc 
$W$ ne rencontre pas $X_{\rm adm}$. On en d\'eduit que $Y_0$ contient $\bs m$.
En particulier $H_0=H$.

Finalement, les deux diviseurs $Y=\pi^\star_{\bs a}(H)$ et $Y_1$ 
sont lin\'eairement \'equivalents et induisent  le m\^eme diviseur au voisinage de $C_0$. 
D'apr\`es le Lemme \ref{fond}, ils sont \'egaux. 
\epf

\ssct{Fin de la d\'emonstration du th\'eor\`eme principal}

Rappelons que, compte-tenu en particulier du Corollaire \ref{reduction},
une vari\'et\'e $X\in \cal{X}_{r+1,n}(q)$ est int\'egrable
si $n=2$, ou si $r=1$, ou si $n\geq 3$ et $q\neq 2n-3$.
le Th\'eor\`eme \ref{Th4}  est donc une cons\'equence du 
r\'esultat plus pr\'ecis  suivant.
\bt
\label{final}
Soit $r\geq 1$, $n\geq 2$ et $q\geq n-1$.
Toute vari\'et\'e int\'egrable de la classe $\cal{X}_{r+1,n}(q)$
est standard. Plus pr\'ecis\'ement, si $\phi : X\dasharrow \P^{r+n-1}$
est une application rationnelle d\'efinie par le syst\`eme lin\'eaire $\cal{Y}(X)$
de la Proposition \ref{systeme},
alors :
\be

\item $X_0=\phi(X)\subset \P^{r+n-1}$  est une vari\'et\'e minimale
de dimension r+1 et de degr\'e $n-1$ et $\phi$ induit une application birationnelle 
$\phi: X\dasharrow X_0$ ;

\item celle-ci induit un diff\'eomorphisme de $X_{\rm adm}$ sur son image, contenue 
dans $(X_0)_{\rm adm}$ ;

\item l'image d'une courbe admissible de $X$ est une courbe admissible de $X_0$
et $\phi$ induit un diff\'eomorphisme de $\Si_q(X)$ sur son image dans $\Si_{n-1}(X_0)$.

\ee
\et
\bpf
On associe des applications rationnelles au syst\`eme lin\'eaire $\cal{Y}(X)$
de la m\^eme fa\c{c}on que si $X$ 
\'etait lisse, voir le d\'ebut de la Section~3.5.
On choisit $Y_0\in \cal{Y}(X)$ et une base 
$(f_0,\ldots,f_{r+n-1})$ de l'espace $Q(Y_0)$.
On pose  
\beq
\label{phi}
\phi(x)=[f_0(x): \cdots : f_{r+1}(x) : \cdots : f_{r+n-1}(x)].
\eeq
L'application $\phi: X\dasharrow \P^{r+n-1}$ d\'efinie 
par (\ref{phi}) ne d\'epend des choix qu'\`a 
composition \`a gauche pr\`es par un automorphisme de $\P^{r+n-1}$.
De plus $X_0=\phi(X)$ engendre $\P^{r+n-1}$.

\sk
Soit $C$ une courbe admissible de $X$ et ${\bs a}\subset C$ un $(n-2)$-uplet.
Pour analyser $\phi$ au voisinage de $C\bck {\bs a}$, 
on choisit $Y_0\in \cal{Y}(X;{\bs a})$ et la base de $Q(Y_0)$ tels
que dans (\ref{phi}), $(f_0,\ldots,f_{r+1})$ soit une base
du sous-syst\`eme lin\'eaire $\cal{Y}(X;{\bs a})$.
L'application rationnelle
$$
\pi_{\bs a}(x)=[f_0(x): \cdots : f_{r+1}(x)]
$$
n'est rien d'autre qu'une application birationnelle $X\dasharrow \P^{r+1}$
associ\'ee \`a $\bs a$. On sait qu'elle induit un morphisme injectif 
sur $X_{\rm adm}({\bs a})$,
de rang constant $r+1$. 

\sk
L'application $\phi: X\dasharrow X_0$ est donc  birationnelle et,
en faisant varier  $\bs a$, on obtient qu'elle est d\'efinie comme 
morphisme sur $X_{\rm adm}$, de rang constant $r+1$.
Elle est injective sur $X_{\rm adm}$ car, pour tout $x_1,x_2\in X_{\rm adm}$,
il existe un $(n-2)$-uplet admissible $\bs a$ tel que $x_1,x_2\in X_{\rm adm}({\bs a})$.
Donc $\phi$ induit un diff\'eomorphisme de $X_{\rm adm}$
sur son image, une sou-vari\'et\'e lisse $\phi(X_{\rm adm})$ de $\P^{r+n-1}$, peut-\^etre non ferm\'ee,
contenue dans $X_0$.

\sk
La vari\'et\'e $X_0$ engendre $\P^{r+n-1}$. Compte tenu des propri\'et\'es du syst\`eme
$\cal{Y}(X)$, si $C\in \Si_q(X)$, son image $\phi(C)$ coupe un hyperplan g\'en\'erique en
$n-1$ point. D'autre part, pour $(x_1,\ldots,x_n)\in X_0^n$ g\'en\'erique, il 
existe une courbe $\phi(C)$ qui passe par ces points. 
Comme un $n$-uplet g\'en\'erique de points de $X_0$ engendre un $\P^{n-1}$,
on conclut que $X_0$ est une vari\'et\'e minimale de degr\'e $n-1$.

\sk
La classification des vari\'et\'es minimales, voir le d\'ebut du Chapitre 5,
montre que $(X_0)_{\rm sing}$
n'a pas de composante irr\'eductible de dimension $r$. Il en r\'esulte 
que $\phi(X_{\rm adm})$ est un ouvert $\Omega$ de $X_0$, contenu dans la partie lisse  
de $X_0$. Pour d\'emontrer la derni\`ere partie de l'\'enonc\'e, puisque 
$\phi: X_{\rm adm} \rightarrow \Omega$ est un diff\'eomorphisme 
et que $\Omega$ est la r\'eunion des images des courbes admissibles de $X$,
il suffit de montrer que l'image d'une courbe admissible de $X$
est une courbe admissible de $X_0$.

\sk
Il suffit pour \c{c}a de r\'ep\'eter le raisonnement de la d\'emonstration 
du Th\'eor\`eme \ref{Th1} en l'appliquant \`a $X_0$ et, au lieu
de ${\rm CRN}_{n-1}(X_0)$, \`a la famille des  courbes $\phi(C)$ avec 
$C\in \Si_q(X)$. Si $\Gamma$ est une telle courbe, on sait que $\Gamma$
est lisse de degr\'e $n-1$ et contenue dans $(X_0)_{\rm reg}$.
Si $(a'_1,\ldots,a'_{n-1})$ est un $(n-1)$-uplet 
de la courbe $\Gamma$ on obtient que  les osculateurs $(X_0)_{a'_1},\ldots(X_0)_{a'_{n-1}}$
sont en somme directe projective. Ainsi, tout $n$-uplet de $\Gamma$
est admissible pour $X_0$. Le th\'eor\`eme est d\'emontr\'e.
\epf

\ssct{Vari\'et\'es int\'egrables et vari\'et\'es standards}

Nous avons d\'emontr\'e qu'une vari\'et\'e int\'egrable $X\in \cal{X}_{r+1,n}(q)$
est standard, mais pas la r\'eciproque.
Du coup, la Proposition \ref{heredite1}, qui est importante, est \'enonc\'ee en termes 
d'int\'egrabilit\'e, une notion transitoire. D'autre part,  
le commentaire qui pr\'ec\`ede cette proposition 
laissait une question ouverte. Les r\'esultats suivants 
r\'epondent \`a ces questions. D'abord, on a : 
\ble
\label{int-st}
Une vari\'et\'e $X\in \cal{X}_{r+1,n}(q)$ est int\'egrable si et seulement si
elle est standard.
\ele
On se donne maintenant une pond\'eration $(\rho_1,\ldots,\rho_{n-1})$ avec $\rho_{n-1}\geq 1$.
On a :
\bt
\label{heredite2}
Soit $X\in \cal{X}_{r+1,n}(q)$ avec $n\geq 3$ et $q\geq n$. Si $a\in X_{\rm adm}$,
on affecte $a$ du poids $\rho_1$, on note 
$\t_a$ une projection osculatrice associ\'ee et $q'=q-(\rho_1+1)$.

\be 

\item Si $X$ est un \'el\'ement standard de $\cal{X}_{r+1,n}(q)$, alors 
$\t_a(X)$ est un \'el\'ement standard de  
$\cal{X}_{r,n-1}(q')$, pour tout $a\in X_{\rm adm}$.

\item Si $n\geq 4$ et si $\t_a(X)$ est un \'el\'ement standard de la classe 
$\cal{X}_{r,n-1}(q')$ pour tout $a\in X_{\rm adm}$,
alors $X$ est un \'el\'ement standard de $\cal{X}_{r+1,n}(q)$.

\ee

\et
Nous d\'emontrerons ces deux \'enonc\'es dans la Section 4.4.

\sk
Une autre d\'emonstration que celle qu'on proposera est possible, d'ailleurs plus \'el\'emen\-taire,
sur la base de la classification des vari\'et\'es standards \'etablie 
dans le Chapitre 5.
Indiquons seulement \`a quoi se r\'eduit la d\'emonstration.
Compte tenu de la Proposition~\ref{heredite1},
la seconde 
partie du th\'eor\`eme ci-dessus est une cons\'equence du Lemme \ref{int-st}.
Compte tenu du fait qu'une vari\'et\'e $X\in \cal{X}_{r+1,n}(q)$ 
est toujours standard si $q\neq 2n-3$, il suffit de d\'emontrer ce lemme et la premi\`ere 
partie du th\'eor\`eme dans le cas $q=2n-3$. Enfin, si la premi\`ere 
partie du th\'eor\`eme est d\'emontr\'ee, on obtient le lemme
par r\'ecurrence sur $n$, \`a partir du cas $n=3$, 
\`a nouveau gr\^ace \`a la Proposition \ref{heredite1}.

\sk
En r\'esum\'e, il suffit de montrer qu'une vari\'et\'e standard de la classe $\cal{X}_{r+1,3}(3)$
est int\'egrable et que, si $X$ est une vari\'et\'e standard d'une classe $\cal{X}_{r+1,n}(2n-3)$
avec $n\geq 4$, pour tout $a\in X_{\rm adm}$, son image $\t_a(X)\in \cal{X}_{r+1,n-1}(2n-5)$
par une projection de centre $X_a(1)$ est une vari\'et\'e standard.
La v\'erification, \`a partir de la description de ces vari\'et\'es dans le Chapitre 5, 
ne pr\'esente pas de difficult\'e.

\sct{La structure infinit\'esimale de la vari\'et\'e des courbes admissibles}

\ssct{Structures quasi-grassmanniennes}

Dans cette section, nous rappelons les d\'efinitions et les r\'esultats 
de la th\'eorie des structures quasi-grassmanniennes dont nous aurons 
besoin, voir par exemple Hangan \cite{Han}, Goldberg \cite{Go}.

\sk
Notons $G_s(W)$ la grassmannienne des sous-espaces vectoriels de dimension $s$ d'un espace vectoriel $W$. 
Soit $V$ un espace vectoriel de dimension $rn$ et soit $u: \C^r\otimes \C^n \rightarrow V$
un isomorphisme.
Si $\rho\in \{1,\ldots,r\}$ et $\nu\in \{1,\ldots,n\}$, l'application 
\beq
\label{type}
G_\rho(\C^r)\times G_\nu(\C^n) \rightarrow G_{\rho \nu}(V),
\qquad 
(E,F) \mapsto u(E\otimes F),
\eeq
d\'efinit une famille de sous-espaces de dimension $\rho\nu$ de $V$. Ce
sont les {\em sous-espaces de type $(\rho,\nu)$} de $V$, pour l'isomorphisme $u$.
Un isomorphisme $u': \C^r\otimes \C^n \rightarrow V$
d\'efinit les m\^emes familles de sous-espaces si et seulement 
l'automorphisme $u'\circ u^{-1}$ de $\C^r\otimes \C^n$
est de la forme $\g\otimes \delta$ avec $\g\in {\rm GL}(\C^r)$ et $\delta\in {\rm GL}(\C^n)$.

Une {\em structure tensorielle vectorielle de type $(r,n)$} sur $V$ est donn\'ee 
par une famille maximale d'isomorphismes lin\'eaires $\C^r\otimes \C^n\rightarrow V$ tels
que si $u,u'$ sont de la famille, $u'\circ u^{-1}$ 
est de la forme qu'on vient de d\'ecrire.
Elle induit naturellement  des structures de type $(\rho,\nu)$
sur les sous-espaces d\'efinis par (\ref{type}). 

\sk
Revenons \`a l'application (\ref{type}).
Un sous-espace de type $(r,n-p)$ est l'intersection
de $p$ sous-espaces $u(\C^r\otimes F_\a)$ de type $(r,n-1)$
en position g\'en\'erale. Il revient au m\^eme de dire 
que les hyperplans $F_1,\ldots,F_p$ de $\C^n$ sont en position g\'en\'erale.
Une remarque analogue vaut pour les sous-espaces de type $(r-p,n)$.

\sk
D'autre part, un sous-espace de type $(r-1,n-1)$ s'\'ecrit 
$u(\C^r\otimes F)\cap u(E\otimes \C^n)$, l'intersection 
de deux sous-espaces de type respectif $(r,n-1)$ et $(r-1,n)$. 
Notons que {\em ces espaces ne sont pas en position 
g\'en\'erale dans $V$} :
$$
\dim \C^r \otimes F + \dim E\otimes \C^n - \dim E\otimes F = 
r(n-1) + (r-1)n - (r-1)(n-1) = rn - 1.
$$

\sk
On aura l'occasion d'utiliser le lemme suivant, emprunt\'e \`a Goldberg \cite{Go}.
\ble 
\label{AG}
Soit $F_0,F_1,\ldots,F_n$ des sous-espaces de codimension $r$ de $V$
en position g\'en\'erale, i.e. tels que $n$ quelconques 
des espaces $F_\a^\perp=\{l\in V',\;\; l_{|F_\a}=0\}$ engendrent le dual
$V'$ de $V$. Il existe une et une seule
structure tensorielle vectorielle de type $(r,n)$ sur $V$ pour laquelle
ces espaces sont des sous-espaces de type $(r,n-1)$ de $V$.
\ele
\bpf
Pour l'existence, on se donne une base $\phi_1,\ldots,\phi_r$ 
de $F_0^\perp$ et l'on d\'e\-compose les $\phi_j$ 
suivant la d\'e\-com\-position $V' = \oplus_{\a =1}^n \, F_\a^\perp$, soit 
$\phi_j = \sum_{\a=1}^n  m_{j \a}$, avec $\,m_{j \alpha}\in F_\a^\perp$.
On v\'erifie que $(m_{j \alpha})$ est une base de $V'$ 
et que $F_0,\ldots,F_n$ sont des espaces de type $(r,n-1)$
pour la structure tensorielle vectorielle d\'efinie par la base 
duale de $(m_{j \a})$. 

\sk
R\'eciproquement, soit $u: \C^r\otimes \C^n\rightarrow V$ un isomorphisme 
et $H_0,\ldots,H_n$ des hyperplans de $\C^n$ en position 
g\'en\'erale, tels que 
$u(\C^r\otimes H_\a)=F_\a$, $\,\a=0,\ldots,n$. 
On peut, sans changer la structure, supposer 
qu'on a $H_\a=\{x\in \C^n, \, x_\a=0\}$ si 
$\a=1,\ldots,n$ et qu'on a $H_0=\{x\in \C^n, \; x_1+\cdots+x_n=0\}$.
Si $u': \C^r\otimes \C^n\rightarrow V$ a la m\^eme
propri\'et\'e, on peut faire la m\^eme r\'eduction sans changer 
la structure d\'efinie par $u'$. On obtient alors 
$(u'\circ u^{-1}) (\C^r\otimes H_\a)=\C^r\otimes H_\a$ pour tout $\a$.
On en  d\'eduit que $u'\circ u^{-1}$ est de la forme 
$\g\otimes I_{\C^n}$ donc que $u$ et $u'$ d\'efinissent la 
m\^eme structure.
\epf

\sk
Dans la suite, on utilisera la terminologie des \og $G$-structures \fg.
Identifions $\C^{rn}$ \`a un produit tensoriel $\C^r\otimes \C^n$
et notons $(g_{j \a ,k \b})$ les \'el\'ements du groupe
de matrices ${\rm GL}(rn)$,
par l\'eg\`eret\'e sans \'ecrire les domaines de variation des indices. 
On note ${\rm G}_{r,n}$ le groupe des matrices 
$(g_{j\a,k\b})$ de la forme 
$$
g_{j \a,k \b} = C_{jk}\, A_{\a\b}, \qquad (C_{jk})_{j,k=1}^r\in {\rm GL}(r), \;\; (A_{\a\b})_{\a,\b=1}^n\in {\rm GL}(n).
$$ 
Une {\em ${\rm G}_{r,n}$-structure vectorielle} sur $V$,
ou {\em structure grassmannienne vectorielle de type $(r,n)$},
est d\'efinie par la donn\'ee d'une famille maximale $\cal{E}$ de bases 
de $V$, dites {\em bases distingu\'ees}, telles que la matrice de passage 
entre deux \'el\'ements $(v_{j\a})$ et $(w_{j\a})$
de $\cal{E}$ appartienne au groupe ${\rm G}_{r,n}$. Autrement dit
la formule de passage est de la forme 
\beq
\label{passage}
w_{j \a} = \sum_{k=1}^r \sum_{\b =1}^n C_{j k} A_{\a \b} \, v_{k \b}.
\eeq
Comme le groupe de matrices ${\rm G}_{r,n}$ est invariant par transposition,
il revient au m\^eme de se donner une famille maximale $\cal{E}'$ de 
bases de l'espace dual $V'$ de $V$, telles que les matrices de passage
appartiennent au groupe ${\rm G}_{r,n}$.

\sk
Toute base distingu\'ee $(v_{j \a})$ de $V$ d\'etermine un isomorphisme 
$u : \C^r\otimes \C^n \rightarrow V$ d\'efini par $u(\varepsilon_j\otimes \varphi_\a)=v_{j \a}$,
o\`u $(\varepsilon_j)_{j=1}^r$ est la base canonique de $\C^r$ et $(\varphi_\a)_{\a=1}^n$
la  base canonique de $\C^n$, et la famille des isomorphismes 
ainsi obtenus d\'efinit une structure tensorielle de type $(r,n)$.
Une ${\rm G}_{r,n}$-structure vectorielle est donc la m\^eme chose qu'une structure 
tensorielle vectorielle de type $(r,n)$.

\bk
Soit maintenant $M$ une vari\'et\'e analytique lisse
de dimension $rn$. Une {\em ${\rm G}_{r,n}$-structure} sur $M$,
ou {\em structure 
quasi-grassmannienne de type $(r,n)$} sur $M$, 
est d\'efinie par la donn\'ee d'une ${\rm G}_{r,n}$-structure vectorielle
sur chaque espace tangent $T_xM$, d\'ependant analytiquement
du point $x$ de $M$. En pratique, elle sera d\'efinie 
par une famille de bases locales de $1$-formes 
telle que les matrices de passage sont des fonctions 
analytiques \`a valeurs dans le groupe ${\rm G}_{r,n}$.
Une telle base est dite {\em distingu\'ee}.

\sk
Une sous-vari\'et\'e lisse $N$ de $M$ est une {\em vari\'et\'e int\'egrale de type $(\rho,\nu)$}
si, pour tout $x\in N$, $\,T_xN$ est un sous-espace  
de type $(\rho,\nu)$ de $T_xM$.

\sk
Par exemple, si $V$ est un espace vectoriel, pour $x\in V$ on identifie  
$T_xV$ \`a $V$. Une ${\rm G}_{r,n}$-structure sur $V$ est {\em constante} 
si elle est d\'efinie par une base de
$1$-formes constantes. Deux ${\rm G}_{r,n}$-structures constantes 
sont \'equivalentes par transformation lin\'eaire.

\sk
Une ${\rm G}_{r,n}$-structure est {\em int\'egrable} ou {\em plate}
si elle est localement diff\'eomorphe \`a une ${\rm G}_{r,n}$-structure 
constante. L'exemple le plus important 
est la ${\rm G}_{r,n}$-structure naturelle sur la grassmannienne $\G_{r,n}$ 
(de dimension $rn$) des $\P^{n-1}\subset \P^{r+n-1}$.
Ses vari\'et\'es int\'egrales de type $(r,n-1)$ sont les sous-vari\'et\'es
de codimension $r$ dont les \'el\'ements passent par un point donn\'e 
de $\P^{r+n-1}$ et ses vari\'et\'es int\'egrales de type $(r-1,n)$ sont les sous-vari\'et\'es
de codimension $n$ dont les \'el\'ements sont contenus dans un hyperplan 
donn\'e.

\sk
On se donne 
une ${\rm G}_{r,n}$-structure sur une vari\'et\'e analytique lisse $M$.
\bd
\label{QG1}
La ${\rm G}_{r,n}$-structure est {\em $\g$-int\'egrable}, respectivement {\em $\delta$-int\'egrable},
si, pour tout $x\in M$ et tout sous-espace $N_x$ de $T_xM$  de type $(r,1)$, respectivement
de type $(1,n)$, il existe un germe en $x$ de vari\'et\'e int\'egrale $N$ de m\^eme
type que $N_x$, tel que $T_xN=N_x$.
\ed
Les lettres $\g$ et $\delta$ sont cens\'ees \'evoquer la gauche et la droite.

La $\g$- et la $\delta$-int\'egrabilit\'e de la structure sont caract\'eris\'ees
par l'annulation de certains tenseurs. Il en r\'esulte, on utilisera cette remarque,
qu'il suffit, dans la d\'efinition pr\'ec\'edente, que les conditions 
soient  v\'erifi\'ees  pour $x\in M$ et $N_x\subset T_xM$ g\'en\'eriques.

\sk
On rappelle le r\'esultat suivant, voir Goldberg \cite{Go}  :
\bt
\label{QG2}
Une ${\rm G}_{r,n}$-structure est int\'egrable si et seulement si elle est $\g$- et $\delta$-int\'egrable.
\et

\ssct{La structure quasi-grasssmannienne de $\Si_q(X)$}

L'objet de ce chapitre est de d\'emontrer le r\'esultat suivant.
\bt
\label{tqg}
Si $X\in \cal{X}_{r+1,n}(q)$, la vari\'et\'e lisse $\Si_q(X)$ des courbes admissibles de $X$
admet une ${\rm G}_{r,n}$-structure naturelle, d\'etermin\'ee par la propri\'et\'e
suivante : 
\bi

\item pour tout point admissible $a$ de $X$, la sous-vari\'et\'e $\Si_q(X;a)$
des courbes admissibles de $X$ qui passent par $a$ est une sous-vari\'et\'e int\'egrale 
de type $(r,n-1)$ de $\Si_q(X)$.

\ei
Cette structure est $\g$-int\'egrable. Elle est int\'egrable si et seulement si
$X$ est un \'el\'ement standard de $\cal{X}_{r+1,n}(q)$.
\et
La premi\`ere partie de l'\'enonc\'e est une cons\'equence du Th\'eor\`eme
\ref{structure-Si} et du fait que l'espace $H^0(\P^1,\L)$ 
a une structure tensorielle vectorielle de type $(r,n)$ naturelle 
si $\L$ est la somme directe de $r$ fibr\'es en droites de degr\'e $n-1$.
Cet espace admet d'ailleurs une structure plus fine qui distingue,
parmi tous les sous-espaces de type $(r,n-1)$, la famille 
des sous-espaces de la forme $\{\xi\in H^0(\P^1,\L), \; \xi(t)=0\}$  quand 
$t$ d\'ecrit $\P^1$. Cette structure, qu'on peut aussi d\'ecrire 
en termes de $G$-structure, est int\'eressante mais n'intervient
pas dans la d\'emonstration. Elle ne sera pas discut\'ee ici,
voir aussi Gindikin \cite{Gi}
\`a ce sujet, dans un cadre plus g\'en\'eral. 

\bk
Soit $C$ une courbe admissible et $V$ l'espace tangent $T_C\Si_q(X)$,
qu'on identifie \`a l'espace $H^0(C,N_CX)$ des sections globales de $N_CX$. \'Ecrivons 
$N_CX = \oplus_{j=1}^r  \L_j$, o\`u $\L_1,\ldots,\L_r$ sont des 
sous-fibr\'es en droites de degr\'e $n-1$ de $N_CX$.
Choisissons un point $x_0$ de $C$, un isomorphisme
$\phi: C\rightarrow \P^1$ tel que $\phi(x_0)=\infty$
et pour $j=1,\ldots,r$, 
une section non nulle $e_j$ de $\L_j$ ayant un z\'ero 
d'ordre $n-1$ en $x_0$. Les $rn$ sections analytiques 
\beq
\label{base}
x\in C, \qquad e_{j \a}(x) = \phi(x)^{\a-1} e_j(x), \;\;\;  \a=1,\ldots,n,
\eeq
de $N_CX$ forment une base $(e_{j \a})$ de $V$.
On munit $V$ de l'unique ${\rm G}_{r,n}$-structure vectorielle pour laquelle 
$(e_{j \a})$ est une base distingu\'ee.
Rappelons que cette structure est aussi la structure tensorielle 
vectorielle d\'efinie par l'isomorphisme lin\'eaire $u: \C^r\otimes \C^n\rightarrow V$
donn\'e par $u(\varepsilon_j\otimes \varphi_\a)=e_{j \a}$, o\`u $(\varepsilon_j)_{j=1}^r$
est la base canonique de $\C^r$ et $(\varphi_\a)_{\a=1}^n$ celle de $\C^n$.
Notons $(m_{j \a})$ la base duale de la base $(e_{j \a})$.

\sk
Un sous-espace de type $(r,n-1)$ de $V$ est  d\'efini par un syst\`eme de la forme :
$$
\sum_{\a=1}^n  t_\a m_{j \a}(\xi) = 0, \qquad j=1,\ldots,r,
$$
o\`u $[t_1:\cdots:t_n]\in \P^{n-1}$.
D'autre part, une section $\xi = \sum_{j, \a} m_{j \a}(\xi)\, e_{j \a}$
s'annule en un point donn\'e $a\in C$ si et seulement si
$$
\sum_{\a=1}^n  \phi(a)^{\a-1} m_{j \a}(\xi) = 0, \qquad j=1,\ldots,r.
$$
Ce syst\`eme d\'efinit un sous-espace de type $(r,n-1)$ de $V$,
d'ailleurs d'une forme particuli\`ere. 
Compte tenu du Lemme \ref{AG}, la ${\rm G}_{r,n}$-structure vectorielle
d\'efinie par la base $(e_{j \a})$ est la seule pour laquelle ces espaces sont de type $(r,n-1)$.
En particulier, elle ne d\'epend pas des choix qu'on a faits.

\sk
Admettons provisoirement que la structure 
qu'on vient de d\'efinir sur $T_C\Si_q(X)$ d\'epend analytiquement de $C$.
On a donc construit une ${\rm G}_{r,n}$-structure sur $\Si_q(X)$.
Comme l'espace tangent $T_C\Si_q(X;a)$ s'identifie au sous-espace 
des sections $\xi\in V$ qui s'annulent en $a$, on
obtient la premi\`ere partie de l'\'enonc\'e.

\bk
Soit $C\in \Si_q(X)$ et $(e_{j \a})$
une base de $V$, de la forme (\ref{base}).

Soit $N_C$ un sous-espace  de type $(r,1)$ de $V$. Par d\'efinition, 
$N_C$ admet une pr\'esentation param\'e\-trique de la forme
$$
\xi = \big(\sum_{\a=1}^n  s_\a \phi^{\a-1}\big) \big(\sum_{j=1}^r t_j e_j\big), \qquad (t_1,\ldots,t_r)\in \C^r,
$$
o\`u $(s_1,\ldots,s_n)\in \C^n\bck \{0\}$.
\'Etant donn\'e $a\in C\bck\{x_0\}$, toutes les  sections $\xi\in N_C$ s'annulent
en ce point si et seulement si $u=\phi(a)$ v\'erifie $\sum_{\a=0}^{n-1} s_\a u^{\a-1}= 0$. 
Pour $s\in \C^n$ g\'en\'eral, les solutions d\'efinissent un $(n-1)$-uplet ${\bs a}=(a_1,\ldots,a_{n-1})$
de $C$ et l'on obtient 
$$
N_C=\bigcap_{i=1}^{n-1} T_C\Si_q(X;a_i) = T_C\Si_q(X;{\bs a}).
$$
Comme $\Si_q(X;{\bs a})=\bigcap_{i=1}^{n-1} \Si_q(X;a_i)$ est une vari\'et\'e int\'egrale de type $(r,1)$
de $\Si_q(X)$, compte tenu de la D\'efinition~\ref{QG1} et du commentaire qui la suit,
on obtient que la ${\rm G}_{r,n}$-structure de $\Si_q(X)$ est $\g$-int\'egrable.

\bk
Pour  montrer que la ${\rm G}_{r,n}$-structure vectorielle qu'on a 
construite 
sur chaque espace $T_C\Si_q(X)$ d\'epend analytiquement de $C\in \Si_q(X)$,
fait qu'on a provisoirement admis, on peut proc\'eder de la fa\c{c}on suivante.

Soit $C\in \Si_q(X)$, $\,(a_0,\ldots,a_n)$
un $(n+1)$-uplet de $C$ et, pour $\a=0,\ldots,n$, $\,Y_\a$ un germe d'hypersurface lisse 
transverse \`a $C$ en $a_\a$. Comme on l'a vu \`a la fin de la Section~2.6,
l'application $\g$, voir le Lemme \ref{anal}, induit un diff\'eomorphisme local
$$
Y_1\times \cdots \times Y_n \rightarrow \Si_q(X). 
$$
Il en r\'esulte que, pour $\a\in \{1,\ldots,n\}$ fix\'e, les vari\'et\'es 
$\Si_q(X;y)$ avec $y\in Y_\a$ forment un germe de feuilletage (r\'egulier)
$\cal{F}_\a$ de codimension $r$ au voisinage de $C$ dans $\Si_q(X)$.
De plus les feuilletages $\cal{F}_1,\ldots,\cal{F}_n$ sont en position g\'en\'erale.
On d\'efinit le feuilletage $\cal{F}_0$ de fa\c{c}on analogue.
Par sym\'etrie, on obtient $n+1$ feuilletages de codimension $r$, en position 
g\'en\'erale.  

On choisit un syst\`eme de $1$-formes analytiques $(\phi_j)_{j=1}^r$
tel que le feuilletage $\cal{F}_0$ soit d\'efini par le 
syst\`eme d'\'equations $\{\phi_j=0, \, j=1,\ldots,r\}$,
et on d\'ecompose chaque $1$-forme $\phi_j$ en $\phi_j = \sum_{\alpha=1}^n  \omega_{j \alpha}$,
o\`u, pour $\a=1,\ldots,n$,  le syst\`eme $\{\omega_{j \a}=0, \; j=1,\ldots,r\}$ 
d\'efinit le feuilletage $\cal{F}_\a$.
Par construction et compte tenu du Lemme \ref{AG},  on obtient ainsi un germe $(\omega_{j \a})$ de base 
de $1$-formes analytiques au voisinage de $C$ et pour tout $C'\in \Si_q(X)$
voisin de $C$, la base de $T_{C'}\Si_q(X)$ induite par cette base 
est distingu\'ee pour la ${\rm G}_{r,n}$-structure vectorielle de $T_{C'}\Si_q(X)$. Ceci donne 
le r\'esultat voulu.

\ssct{Fin de la d\'emonstration du Th\'eor\`eme \ref{tqg}}

Si $X_0\subset \P^{r+n-1}$ est une vari\'et\'e minimale de degr\'e $n-1$,
la structure quasi-grassmannienne de $\Si_{n-1}(X_0)$
est induite par l'application qui associe \`a un $\P^{n-1}\subset \P^{r+n-1}$
g\'en\'erique la section $X_0\cap \P^{n-1}$.
Elle est int\'egrable.
Par d\'efinition d'une vari\'et\'e standard de la classe $\cal{X}_{r+1,n}(q)$, \'eventuellement
en invoquant le Th\'eor\`eme \ref{final} pour plus de pr\'ecision,
la ${\rm G}_{r,n}$-structure naturelle de $\Si_q(X)$ est int\'egrable si la vari\'et\'e $X$
est standard.

\sk
On traite la r\'eciproque. 
On suppose que la ${\rm G}_{r,n}$-stucture de $\Si_q(X)$
est int\'egrable. Compte tenu du Th\'eor\`eme \ref{final}, il s'agit de montrer que $X$ est int\'egrable
au sens de la D\'efinition \ref{D4.2}.

\sk
Soit $C$ une courbe admissible de $X$,
$\,{\bs a}\subset C$ un $(n-2)$-uplet et $\pi: X\dasharrow \P^{r+1}$ une application birationnelle 
associ\'ee. Soit $x$ un point 
de $C\bck {\boldsymbol a}$, $\,h\subset T_xX$ un hyperplan qui contient 
$T_xC$ et $H$ l'hyperplan de $\P^{r+1}$ passant par $\pi(x)$ et de direction 
l'image de $h$ dans $T_{\pi(x)}\P^{r+1}$. La vari\'et\'e $Y_0=\pi^{-1}(H)$
est un diviseur admissible qui contient $C$ et est tel que $T_xY_0=h$.
Comme on l'a remarqu\'e apr\`es la
D\'efinition \ref{D4.1},  
on obtient tous les diviseurs admissibles de $X$ par cette construction.
Il s'agit donc de montrer que $Y_0$ appartient \`a la classe $\cal{X}_{r,n}(q)$.

\sk
Soit $V=T_C\Si_q(X)$, qu'on identifie \`a $H^0(C,N_CX)$, 
et $(e_{j \a})$ une base distingu\'ee de $V$ de la forme (\ref{base}). 
Par d\'efinition, un sous-espace $L_C$ de type $(r-1,n)$ de $V$
admet une pr\'esen\-tation param\'e\-trique de la forme 
$$
\xi = \big( \sum_{\a=1}^n s_\a\phi^{\a-1}\big)\big( \sum_{j=1}^r t_j e_j\big),
\qquad (s_1,\ldots,s_n)\in \C^n, \;\; (t_1,\ldots,t_r)\in K,
$$
o\`u $K$ est un hyperplan de $\C^r$. La fibre $L_C(p) = \{\xi(p)\in (N_CX)_p, \;\; \xi\in L_C\}$
de $L_C$ en tout point $p\in C$ est donc de rang $r-1$.
On choisit le sous-espace $L_C$, de type $(r-1,n)$,
dont la fibre $L_C(x)$ au point donn\'e $x\in C\bck {\bs a}$ est $h/T_xC$. 

\sk
Comme la ${\rm G}_{r,n}$-structure de $\Si_q(X)$
est suppos\'ee int\'egrable,
il existe en $C$ un germe $L\subset \Si_q(X)$ de vari\'et\'e int\'egrale de type $(r-1,n)$,
tel que $T_CL=L_C$. 
Consid\'erons la projection canonique 
$$
\kappa : \, \Big\{(C', x')\in L\times X, \;\; x'\in C'\Big\} \rightarrow X.
$$
Son image est aussi l'image d'un rel\`evement $V: L\times \P^1 \rightarrow X$
de l'inclusion $v: L \rightarrow \Si_q(X)$. Les d\'eriv\'ees de $v$ et de $V$
en $C'\in L$ et $(C',t)\in L\times \P^1$ sont reli\'ees 
par la formule~(\ref{diff}) du Chapitre 2. Si $x'=V(C',t)$, 
l'image de la d\'eriv\'ee de $V$ est le sous-espace de $T_{x'}X$ qui contient $T_{x'}C'$
et dont le quotient par $T_{x'}C'$ est la fibre en $x'$ du sous-espace 
$T_{C'}L$ de $H^0(C',N_{C'}X)$. En particulier, comme cette fibre est de dimension
constante $r-1$, l'application $V$ est de rang constant $r$ et donc, si le germe 
$L$ est choisi assez petit,  son image $Y=\kappa (L)$
est une hypersurface analytique (non ferm\'ee) lisse qui contient $C$. Par construction
$T_xY=h$. 

\sk
Consid\'erons \`a pr\'esent la projection canonique 
$$
\kappa^{(n)} : \, \Big\{(C', {\bs x}')\in L\times X^{(n)}_{\rm adm}, \;\; {\bs x}'\subset C'\Big\} \rightarrow X^{(n)}_{\rm adm}.
$$
La vari\'et\'e de d\'epart est de dimension $(r-1)n + n=rn$ et si 
${\bs x}'$ est dans l'image de $\kappa^{(n)}$,  $\,(\g({\bs x}'),{\bs x}')$
est son seul ant\'ec\'edent. L'application $\kappa^{(n)}$,
\`a valeurs dans $Y^n$, est donc de rang $rn = \dim Y^n$ au point g\'en\'erique 
de la vari\'et\'e source.
Rappelons que $Y$ est lisse le long de $C$ donc engendre 
un espace de dimension $\geq \pi_{r-1,n}(q)$. On a obtenu :

\sk
{\em la vari\'et\'e irr\'eductible lisse (non  ferm\'ee) $Y=\kappa(L)$ engendre un espace 
de dimension $\geq \pi_{r-1,n}(q)$ et il existe un $n$-uplet ${\bs y}_0\in Y^n$
tel que, pour tout ${\bs y}\in Y^n$ voisin de ${\bs y}_0$,
il existe une courbe $C\in {\rm CRN}_q(Y)$ qui contient~$\bs y$.}

\sk
Bien que la vari\'et\'e $Y$ ne soit pas ferm\'ee, le  m\^eme argument qu'au d\'ebut de 
la d\'emonstration de la Proposition~\ref{equivalences} montre
que $\pi(Y)$ est contenu dans un hyperplan de $\P^{r+1}$, n\'ecessairement l'hyperplan 
$H$. Il en r\'esulte que $Y$ est un ouvert de $Y_0=\pi^{-1}(H)$
puis que $Y_0$ appartient \`a la classe $\cal{X}_{r,n}(q)$. 
Le Th\'eor\`eme \ref{tqg} est d\'emontr\'e.

\ssct{Application : retour sur la Section 3.7}

On vient de montrer que si la vari\'et\'e $X\in \cal{X}_{r+1,n}(q)$
est standard, la ${\rm G}_{r,n}$-structure naturelle
de $\Si_q(X)$ est int\'egrable, puis que, si cette structure 
est int\'egrable, alors $X$ est int\'egrable.
Pour $X$, les deux propri\'et\'es, d'\^etre int\'egrable et d'\^etre standard, sont 
donc \'equivalentes. Ceci d\'emontre le Lemme~\ref{int-st} et la seconde partie 
du Th\'eor\`eme \ref{heredite2}.

\bk
Reste \`a d\'emontrer la premi\`ere partie de ce th\'eor\`eme. Soit $X$ une vari\'et\'e 
standard de la classe $\cal{X}_{r+1,n}(q)$ avec $n\geq 3$ et $q\geq n$.
On se donne un point $a\in X_{\rm adm}$ qu'on affecte du poids $\rho_1$.
Il s'agit de montrer que la vari\'et\'e $X'=\t_a(X)\in \cal{X}_{r+1,n-1}$
est standard, o\`u $\t_a$ est une projection osculatrice associ\'ee \`a $a$
et $q'=q-(\rho_1+1)$.

\sk
On sait que $\t_a$ induit un diff\'eomorphisme de
$X_{\rm adm}(a)$ sur son image, un ouvert dense de $X'_{\rm adm}$,
ainsi qu'une application $\wh{\t}_a : \Si_q(X;a)\rightarrow \Si_{q'}(X')$.
En utilisant les rappels faits au d\'ebut de la Section 2.6,
on calcule facilement la d\'eriv\'ee de cette application en $C\in \Si_q(X)$
et l'on v\'erifie que c'est un isomorphisme.
On en d\'eduit que $\wh{\t}_a$ est un diff\'eomorphisme de $\Si_q(X;a)$ sur son image, 
un ouvert dense de $\Si_{q'}(X')$. 

\sk
Finalement, comme la ${\rm G}_{r,n}$-structure de $\Si_q(X)$ est int\'egrable,
la ${\rm G}_{r,n-1}$-structure induite sur $\Si_q(X;a)$, une vari\'et\'e
int\'egrale de type $(r,n-1)$,
est int\'egrable. Compte tenu du Th\'eor\`eme \ref{tqg},
le diff\'eomorphisme $\wh{\t}_a$ est compatible avec les ${\rm G}_{r,n-1}$
structures de la source et du but. 
Il en r\'esulte que la ${\rm G}_{r,n-1}$-structure de $\Si_{q'}(X')$
est int\'egrable sur un ouvert dense de $\Si_{q'}(X')$. Elle est
donc int\'egrable et, d'apr\`es le m\^eme th\'eor\`eme,
la vari\'et\'e $X'$ est standard.

\ssct{Une autre d\'emonstration du Th\'eor\`eme \ref{Th4} ?\,}

Dans la Section 4.3, pour montrer que $X$ est une vari\'et\'e standard
si la ${\rm G}_{r,n}$-structure de $\Si_q(X)$
est int\'egrable, nous avons utilis\'e 
la notion interm\'ediaire de \og vari\'et\'e int\'egrable \fg\,
et le Th\'eor\`eme~\ref{final}.
D'un autre c\^ot\'e, les m\'ethodes d\'ecrites dans ce chapitre 
sugg\`erent une autre d\'emonstration du Th\'eor\`eme \ref{Th4}.
L'esquisse ci-dessous est seulement plausible.
Faute de temps et aussi faute de motivation, nous 
n'avons pas v\'erifi\'e les d\'etails.

\sk
Un probl\`eme g\'en\'eral int\'eressant est de caract\'eriser les paires
$(X,\Si(X))$, o\`u $X$ est une vari\'et\'e alg\'ebrique irr\'eductible 
de dimension $r+1$ et $\Si(X)$ une vari\'et\'e alg\'ebrique irr\'eductible 
de dimension $rn$ de courbes de $X$,
qui sont \'equivalentes par  une transformation birationnelle 
\`a une paire $(X_0,\Si(X_0))$ o\`u la vari\'et\'e $X_0$ engendre un espace $\P^{r+n-1}$
et $\Si(X_0)$ est la vari\'et\'e d\'efinie \`a partir des 
sections $X_0\cap \P^{n-1}$. Disons que la paire $(X,\Si(X))$
est {\em standard} si c'est le cas. 
C'est l'analogue pour les courbes du probl\`eme de la 
caract\'erisation  des syst\`emes de diviseurs qui sont 
lin\'eaires. Dans ce cas-ci, un th\'eor\`eme d'Enriques \cite{En}
donne une r\'eponse importante.

\sk
On consid\`ere la situation qu'on vient de d\'ecrire,
sous l'hypoth\`ese que la relation d'incidence \og $x\in C$ \fg\, 
entre points $x\in X$ et courbes $C\in \Si(X)$
d\'etermine, au sens de la premi\`ere partie du Th\'eor\`eme~\ref{tqg},
une structure quasi-grassmannienne sur un ouvert 
dense de $\Si(X)$ et que celle-ci est int\'egrable.
On suppose de plus que, pour ${\bs x}\in X^n$
g\'en\'erique, il existe une {\em et une seule} courbe $C\in \Si(X)$
qui contient $\bs x$.
Il est tentant de penser (est-ce connu ?) que, sous ces hypoth\`eses,
la paire $(X,\Si(X))$ est standard.

La Section 4.3 nous dit au moins comment essayer de construire le syst\`eme de diviseurs qui,
si l'on montre qu'il est lin\'eaire,
r\'esout le probl\`eme pos\'e. Un \'el\'ement g\'en\'erique $Y$
du syst\`eme devra contenir la r\'eunion des courbes 
qui sont \'el\'ements d'une vari\'et\'e (ou d'un  germe de vari\'et\'e) int\'egrale 
de type $(r-1,n)$ de $\Si(X)$, d\'ependant de $Y$.
 
\sk
Si les conjectures pr\'ec\'edentes sont correctes et d\'emontr\'ees,
une d\'emons\-tration alternative du Th\'eor\`eme \ref{Th4} consiste \`a ne retenir du Chapitre 3
que la Proposition \ref{P4.2}, dont la d\'emonstration 
est tr\`es facile, \`a en d\'eduire le r\'esultat voulu 
si $n=3$ et \`a conclure par l'\'enonc\'e suivant 
de la th\'eorie des structures quasi-grassmanniennes.
\bt
\label{Ge}
Soit $r\geq 2$, $\,n\geq 3$ et $M$ une vari\'et\'e munie 
d'une ${\rm G}_{r,n}$-structure $\g$-int\'egrable.
Soit $\cal{F}_0,\cal{F}_1,\ldots,\cal{F}_n$ des feuilletages de codimension $r$ au voisinage de $x_0\in M$,
en position g\'en\'erale, dont les feuilles sont des vari\'et\'es int\'egrales de type $(r,n-1)$.
Si la ${\rm G}_{r,n-1}$-structure induite sur chaque feuille de chaque 
feuilletage est int\'egrable et si $n\geq 4$, alors la ${\rm G}_{r,n}$-structure 
de $M$ est  int\'egrable au voisinage de $x_0$.
\et
On remarquera l'analogie entre l'\'enonc\'e pr\'ec\'edent et le Th\'eor\`eme  \ref{heredite2}
sur l'h\'er\'edit\'e de l'int\'egrabilit\'e. En fait nous avons d'abord d\'emontr\'e
un analogue de ce dernier en utilisant le Th\'eor\`eme \ref{Ge}.
La d\'emonstration directe qu'on a donn\'ee dans le Chapitre 3 n'a \'et\'e d\'ecouverte 
que plus tard.

\bk
Le Th\'eor\`eme \ref{Ge} appara\^{\i}t dans Goldberg \cite{Go}, o\`u le cas $n=3$ est permis,
mais tout exemple d'une vari\'et\'e sp\'eciale d'une classe $\cal{X}_{r+1,3}(3)$,
on en construira dans le Chapitre~6, 
montre que le r\'esultat est faux si $n=3$.

En effet, soit $X\in \cal{X}_{r+1,3}(3)$ une vari\'et\'e sp\'eciale. La vari\'et\'e
$\Si_3(X)$ est munie d'une ${\rm G}_{r,3}$-structure $\gamma$-int\'egrable naturelle,
qui n'est pas int\'egrable puisqu'on suppose que $X$ est sp\'eciale. 
D'autre part, si $a$ est un point admissible de $X$, on a vu que 
la ${\rm G}_{r,2}$-structure induite sur la vari\'et\'e $\Si_3(X;a)$, une
vari\'et\'e int\'egrale de type $(r,2)$, est localement \'equivalente 
\`a la ${\rm G}_{r,2}$-structure naturelle sur la vari\'et\'e $\Si_2(X')$,
o\`u $X'$ est une vari\'et\'e de Veronese d'ordre $2$, l'image de $X$
par une projection osculatrice de centre $X_a(1)$. 
Cette structure est donc toujours int\'egrable.

D'autre part, on a vu \`a la fin de la Section 4.2 comment construire 
des feuilletages $\cal{F}_0,\cal{F}_2,\ldots,\cal{F}_n$ de codimension $r$ au voisinage de $C\in \Si_3(X)$,
en position g\'en\'erale, dont les feuilles sont des vari\'et\'es de la forme $\Si_3(X;a)$
avec $a\in X_{\rm adm}$. Ceci montre que la condition $n\geq 4$ est 
n\'ecessaire dans l'\'enonc\'e ci-dessus.

\sk
En revanche, la $\delta$-int\'egrabilit\'e des ${\rm G}_{r,n}$- et ${\rm G}_{r,n-1}$-structures
de l'\'enonc\'e se lit sur des  tenseurs dont des formules explicites 
sont donn\'ees dans Goldberg \cite{Go} et {\em qui sont faciles \`a calculer sous l'hypoth\`ese $n\geq 4$}.
Bien que ce ne soit pas la d\'emonstration choisie dans \cite{Go} (elle est incorrecte),
peut-\^etre dans l'espoir de couvrir aussi le cas $n=3$, la d\'emonstration du Th\'eor\`eme \ref{Ge}
est tr\`es simple \`a partir de ces formules.

\sct{La classification des vari\'et\'es standards}

\ssct{Pr\'eliminaires}

Deux vari\'et\'es projectives $X$ et $X'$ sont {\em \'equivalentes} 
s'il existe un isomorphisme $\phi: \lan X\ran \rightarrow \lan X'\ran$
tel que $\phi(X)=X'$.
Dans ce chapitre, nous donnons la classification des vari\'et\'es standards 
\`a \'equivalence pr\`es.

\sk
Le r\'esultat pr\'ecis est \'enonc\'e dans la prochaine section.
La d\'emonstration occupe le reste du chapitre et suit la classification 
des vari\'et\'es minimales, qu'on rappellera.

\sk
Soit $X_0\in \cal{X}_{r+1,n}(n-1)$ une vari\'et\'e minimale
et $X\in \cal{X}_{r+1,n}(q)$ une vari\'et\'e standard, voir 
la D\'efinition \ref{D2}.
On dit qu'une application 
birationnelle $\phi : \; X_0\dasharrow X$, telle que l'image 
d'une section de $X_0$ par un $\P^{n-1}\subset \P^{r+n-1}$ g\'en\'erique
 appartient \`a $\text{\rm CRN}_q(X)$,
{\em associe $X$ \`a $X_0$}  et, si une telle 
application existe, que {\em $X$ et  $X_0$ sont associ\'ees}.

\sk
Pour chaque vari\'et\'e minimale $X_0$, on construira
explicitement des vari\'et\'es standards qui lui sont
associ\'ees. On montrera que la liste donn\'ee est exhaustive
en d\'eterminant tous les syst\`emes lin\'eaires sur $X_0$
tels que les applications rationnelles qu'ils d\'efinissent associent
$X_0$ \`a des vari\'et\'es standards. On utilisera sans d\'emonstration
les propri\'et\'es connues du groupe de Picard de $X_0$\footnote{
Rappelons en particulier que si $X_0$ est un c\^one au-dessus de $X'_0$,
le groupe de Picard de $X_0$ est naturellement isomorphe \`a celui de $X'_0$ ;
voir R. Hartshorne : Algebraic Geometry, {\em Graduate Texts in Mathematics~52},
Springer-Verlag, Chapitre II, Exercice 6.3. On utilisera cette propri\'et\'e.} ainsi que le lemme suivant. 
\ble
Soit $X_0\in \cal{X}_{r+1,n}(n-1)$ une vari\'et\'e minimale et $\phi : X_0\dasharrow X$
une application birationnelle qui associe une vari\'et\'e standard $X\in \cal{X}_{r+1,n}(q)$
\`a $X_0$. L'application $\phi$ est d\'efinie par un syst\`eme lin\'eaire 
complet $|W|$ de dimension $\pi_{r,n}(q)$, tel que $C\cdot W=q$
pour toute section de $X_0$ par un $\P^{n-1}\subset \P^{r+n-1}$ g\'en\'erique.
\ele
\bpf
On se place sous les hypoth\`eses de l'\'enonc\'e. D'abord,
la partie singuli\`ere d'une vari\'et\'e minimale
est vide ou de codimension $\geq 2$, ce qui permet d'utiliser la th\'eorie des 
syst\`emes lin\'eaires. On note $\text{Pic}\, (X_0)$ le groupe des classes de diviseurs 
de $X_0$ modulo \'equivalence lin\'eaire : deux diviseurs $W,W'\in \text{Div}(X_0)$
sont lin\'eairement \'equivalents si $W-W'$ est le diviseur d'une 
fonction rationnelle $f\in R(X_0)\bck\{0\}$.

\sk
L'application $\phi$ est d\'efinie par un syst\`eme lin\'eaire $\cal{L}$
sans composante fixe et de dimension $\pi_{r,n}(q)$, uniquement d\'etermin\'e. 
C'est un sous-syst\`eme lin\'eaire 
d'un syst\`eme lin\'eaire complet $|W|$ sans composante fixe,
d\'efini par une classe $W\in \text{Pic}\,(X_0)$,

\sk
Comme $\phi: \, X_0 \dasharrow X$ est birationnelle et qu'une section 
g\'en\'erique $C=X_0\cap \P^{n-1}$  
est contenue dans $(X_0)_{\rm reg}$ et ne rencontre pas 
le lieu d'ind\'etermination de $\phi$, son image 
$\phi(C)$ est de degr\'e $C\cdot W$, ce qui donne $C\cdot W = q$. 

\sk
L'image $X'$ de $X_0$ par une application rationnelle 
d\'efinie par le syst\`eme lin\'eaire {\em complet} $|W|$
a aussi la propri\'et\'e que, pour ${\bs x}\in X'^n$
g\'en\'erique, il existe un \'el\'ement de $\text{CRN}_q(X')$
qui contient $\bs x$. Compte tenu du Th\'eor\`eme \ref{Th1},
la vari\'et\'e $X'$ engendre un espace de dimension $\leq \pi_{r,n}(q)$. 
Comme le sous-syst\`eme $\cal{L}$ du syst\`eme $|W|$ est d\'ej\`a de dimension $\pi_{r,n}(q)$, on obtient $\cal{L}=|W|$.
\epf

\ssct{La classification}

Rappelons d'abord la classification des vari\'et\'es minimales $X_0\subset \P^{r+n-1}$, de dimension $r+1$
et de degr\'e $n-1$.

Si $n=2$, $X_0$ est un espace projectif $\P^{r+1}$. 

Si $n=3$ et $r\geq 2$, $X_0$ peut \^etre une hyperquadrique de $\P^{r+2}$, de rang 
$\geq 5$. 

Si $n=5$, $X_0$ peut \^etre un c\^one au-dessus d'une surface de Veronese d'ordre $2$.

Si $n\geq 3$ et si $X_0$ n'est pas de l'une des deux formes pr\'ec\'edentes,
$X_0$ est un scroll rationnel normal. Les scrolls rationnels normaux de dimension 
$r+1$ et de degr\'e $n-1$ sont caract\'eris\'es \`a \'equivalence pr\`es 
par une famille d'entiers $a_0,\ldots,a_r$ tels que :
\beq
\label{S-a}
a_0\geq  \cdots \geq a_r \geq 0, \qquad a_0 +  \cdots + a_r = n-1.
\eeq
Le mod\`ele $S_{a_0,\ldots,a_r}$ d'un tel scroll est donn\'e dans
le paragraphe ci-dessous.

\sk
Comme on va le voir, une vari\'et\'e standard est le plus souvent 
\'equivalente \`a la compl\'et\'ee projective d'une sous-vari\'et\'e d'un 
espace $\C^N$ param\'etr\'ee par une famille de mon\^omes.
On introduit une d\'efinition purement utilitaire.
\bd
\label{var-mon}
Soit $A\subset \N^{r+1}\bck\{0\}$ un ensemble fini de
cardinal $N_A\geq 1$. On note $x=(x_\a)_{\a\in A}$ le point courant de $\C^{N_A}$
et, pour tout $s=(s_1,\ldots,s_{r+1})\in \C^{r+1}$ et tout $\a=(\a_1,\ldots,\a_{r+1})\in A$,
on note $s^\a= s_1^{\a_1}\cdots s_{r+1}^{\a_{r+1}}$. 
{\em La vari\'et\'e mon\^omiale d\'efinie par l'ensemble $A$}
est la vari\'et\'e projective $X\subset \P^{N_A}$ 
telle que 
\beq
X\cap \C^{N_A} = \big\{ (s^\a)_{\a\in A},  \;\; s \in \C^{r+1}\big\}.
\eeq
\ed
On distinguera souvent certaines des composantes de $s\in \C^{r+1}$ en 
changeant de notation.
Par exemple, \'etant donn\'ee une famille d'entiers $(a_0,\ldots,a_r)$ 
qui v\'erifie (\ref{S-a}), la vari\'et\'e mon\^omiale 
d\'efinie par l'ensemble 
$$
A =\Big\{(k,\a)\in (\N\times \N^r)\bck \{0\}, \;\;\;  |\a|\leq 1, \;\; k \leq (\rho-1)a_0 + \sum_{j=1}^r \a_j a_j\Big\}.
$$
est un scroll rationnel normal de dimension $r+1$ et de degr\'e $n-1$, qu'on note $S_{a_0,\ldots,a_r}$.

\bk
Le th\'eor\`eme suivant donne la liste, sans omission ni r\'ep\'etition,
des vari\'et\'es standards.
\bt
\label{class}
Soit $X\in \cal{X}_{r+1,n}(q)$ une vari\'et\'e standard avec $r\geq 1$, 
$\,n\geq 2$ et $q\geq n-1$.  Elle est associ\'ee 
\`a une seule vari\'et\'e minimale $X_0$, \`a \'equivalence pr\`es.
\be
\item  Si $n=2$ donc $X_0=\P^{r+1}$, $\,X$ est une vari\'et\'e de Veronese d'ordre $q$.
\item  Si $n=3$ et si $X_0$ est une  hyperquadrique de rang $\geq 5$, $q$ est pair et $X$ est l'image de $X_0$
par un plongement de Veronese d'ordre~$q/2$.
\item  Si $n=5$ et si $X_0$ est un c\^one au-dessus d'une surface de Veronese d'ordre~$2$,
$q$ est pair et $X$ est \'equivalente \`a la vari\'et\'e mon\^omiale 
d\'efinie par l'ensemble 
$$
A(q) = \Big\{(i,j,\a)\in \N^2\times \N^{r-1}, \;\;\; 1\leq 2(i+j)+4|\a| \leq q\Big\}.
$$
\item Si $n\geq 3$ et si $X_0$ est le scroll rationnel normal $S_{a_0,\ldots,a_r}$, de degr\'e $n-1$,
$X$ est \'equivalente \`a une vari\'et\'e mon\^omiale d\'efinie par un
ensemble 
$$
A(\rho,\chi)=\Big\{(k,\a)\in (\N\times \N^r)\bck \{0\}, \;\;\;  |\a|\leq \rho, \;\; k \leq (\rho-|\a|)a_0 + \sum_{j=1}^r \a_j a_j + \chi\Big\},
$$
o\`u $q=\rho(n-1)+\chi$, $\,\rho\geq 1$ est un entier et $\,\chi\in \{-1,\ldots,n-2\}$.

\nk
Si $q\not\equiv -1$ modulo $n-1$, $\, q = \rho (n-1) + \chi$ est la division euclidienne de $q$ par $n-1$.

\nk
Si  $q\equiv -1$ modulo $n-1$, les ensembles $A(\rho,n-2)$ et $A(\rho+1,-1)$ d\'efinissent des vari\'et\'es 
standards ; elles sont \'equivalentes si et seulement si $n=3$ ou $a_0=n-1$. 
\ee
\et
Le premi\`ere partie de l'\'enonc\'e est une cons\'equence de la seconde. En effet, si 
une vari\'et\'e $X\in \cal{X}_{r+1,n}(q)$ est 
associ\'ee \`a deux vari\'et\'es $X_0, X_1\in \cal{X}_{r+1,n}(n-1)$, $X_1$  est associ\'ee \`a $X_0$ et 
donc \'equivalente \`a $X_0$ d'apr\`es la classification.

\sk
Remarquons qu'une m\^eme vari\'et\'e peut appartenir \`a plusieurs classes $\cal{X}_{r+1,n}(q)$,
comme le montre la derni\`ere partie de la classification.
Si $1\leq \chi\leq n-2$, l'ensemble $A(1,\chi)$ d\'efinit une 
vari\'et\'e standard de la classe $\cal{X}_{r+1,n-1}(n-1+\chi)$,
qui est aussi le scroll rationnel normal $S_{a_0+\chi,\ldots,a_r+\chi}\in \cal{X}_{r+1,n'}(n'-1)$ 
avec $n'=n+(r+1)\chi$.

\sk
Avant de passer \`a la d\'emonstration du th\'eor\`eme, faisons quelques 
remarques sur la D\'efinition \ref{var-mon}.
Les ensembles de multi-indices $A\subset \N^{r+1}\bck\{0\}$ 
qui interviendront auront tous les propri\'et\'es suivantes :
\bi
\item si $\a\in \N^{r+1}$ est de longueur $|\a|=1$, alors $\a\in A$ ;
\item si $\a\in A$ et si $\a' \in \N^{r+1}\bck\{0\}$ est tel que 
$\a'_i\leq \a_i$ pour $i=1,\ldots,r+1$, alors $\a'\in A$.
\ei

\sk
Soit $X$ une vari\'et\'e mon\^omiale d\'efinie par un ensemble $A$
qui a ces propri\'et\'es. 
Il r\'esulte de la premi\`ere propri\'et\'e que $X\cap \C^{N_A}$ est une vari\'et\'e 
isomorphe \`a $\C^{r+1}$.

Il r\'esulte de la seconde propri\'et\'e que $X\cap \C^{N_A}$
est homog\`ene. Plus pr\'ecis\'ement, pour $s_\star\in \C^{r+1}$ donn\'e, $X\cap \C^{N_A}$
est aussi param\'etr\'e par $x_\a = (s_\star + s)^\a$, $\,\a\in A$.
Compte tenu de la seconde condition, les polyn\^omes $(s_\star + s)^\a - s_\star^\a$, $\,\a\in A$, 
engendrent le m\^eme espace vectoriel que les mon\^omes $s^\a$, $\,\a\in A$. Il existe 
donc une transformation affine de $\C^{N_A}$ qui conserve $X$ et envoie $0\in X$
sur un point donn\'e de $X\cap \C^{N_A}$.

Si $X,X'$ sont des vari\'et\'es mon\^omiales d\'efinies 
par des ensembles $A,A'\subset \N^{r+1}\bck \{0\}$ qui v\'erifient les conditions pr\'ec\'edentes
et si $A'\subset A$,
la projection $(x_\a)_{\a\in A}\mapsto (x_\a)_{\a\in A'}$
induit un isomorphisme de $X\cap \C^{N_A}$ sur $X'\cap \C^{N_{A'}}$.

\ssct{Les trois premiers points de la classification}

On les d\'emontre au cas par cas.

\sk
Si $n=2$, le groupe de Picard de $X_0=\P^{r+1}$ est engendr\'e 
par la classe $H$ d'un hyperplan. Si $q \geq 1$,
les applications rationnelles induites par le syst\`eme $|qH|$
sont les plongements de Veronese d'ordre $q$. Ceci donne 
le premier point du th\'eor\`eme. Bien s\^ur,
le Th\'eor\`eme \ref{Th2} est bien plus fort, puiqu'on n'y suppose 
pas que $X$ soit standard.

\sk
Soit $n=3$ et $X_0\subset \P^{r+2}$ une hyperquadrique 
de rang $\mu+1 \geq 5$, ce qui impose en particulier $r\geq 3$.
Comme $X_0$ est un c\^one au-dessus d'une hyperquadrique lisse
de rang $\geq 5$, son groupe 
de Picard est le groupe libre engendr\'e par 
la classe $H$ d'une section hyperplane de $X_0$. Si $\rho\geq 1$,
l'application rationnelle d\'efinie par le syst\`eme lin\'eaire $|\rho H|$
est la restriction \`a $X_0$ d'un plongement de Veronese d'ordre $\rho$
de $\P^{r+2}$ et l'image d'une section de $X_0$
par un $\P^2\subset \P^{r+2}$ g\'en\'erique 
est une courbe rationnelle 
normale de degr\'e~$2\rho$.

\sk
D'autre part, dans un syst\`eme convenable de coordonn\'ees 
homog\`enes $[U_0 : \cdots : U_{r+2}]$, la quadrique $X_0$ est 
donn\'ee par l'\'equation $\sum_{j=0}^\mu U_j^2 = 0$.
En notant $U=(U_0,U')$, comme un polyn\^ome homog\`ene de degr\'e $\rho$ s'\'ecrit 
$$
F(U) = G_0(U') + U_0G_1(U') + (\sum_{j=0}^\mu U_j^2)G_2(U),
$$
o\`u $G_0(U')$ est homog\`ene de degr\'e $\rho$ et $G_1(U')$ homog\`ene de degr\'e $\rho-1$,
on voit que  
$$
\dim |\rho H|  = { r+1+\rho \choose r+1} + { r+\rho \choose r+1} - 1  = \pi_{r,3}(2\rho).
$$
On a obtenu le  deuxi\`eme point du th\'eor\`eme.

\sk
Soit $n=5$ et $X_0\subset \P^{r+4}$ un c\^one au-dessus d'une surface de Veronese $S\subset \P^5$.
Le groupe de Picard de $S$ est le groupe libre  engendr\'e 
par une conique de $S$, donc celui de $X_0$ est le groupe libre engendr\'e par la classe
$R$ d'un c\^one au-dessus d'une conique de $S$ et $2R$ est la classe d'une section 
hyperplane de $X_0$. 

Si  $C\in \text{CRN}_4(X_0)$ est la section de $X_0$
par un $\P^4\subset \P^{r+4}$ g\'en\'erique,
alors $C\cdot \s R = 2\s$. Il en r\'esulte que, si $X\in \cal{X}_{r+1,5}(q)$
est associ\'ee \`a $X_0$, $q$ est pair et $X$ d\'etermin\'ee
par $q$, \`a \'equivalence pr\`es. En particulier, si $q$ est un multiple de $4$,
$X$ est l'image de $X_0$ par un plongement de Veronese d'ordre $q/4$ de $X_0$.

\sk
R\'eciproquement, $r\geq 1$ \'etant fix\'e, soit $q\geq 4$ un entier pair et  
$$
A(q) = \Big\{(i,j,\a)\in \N^2\times \N^{r-1}, \qquad 1\leq 2(i+j)+4|\a| \leq q\Big\}.
$$
Notons $X(q)$ la vari\'et\'e mon\^omiale d\'efinie  par l'ensemble $A(q)$.
On v\'erifie que $X(4)$ est un c\^one au-dessus d'une surface de Veronese.
On peut supposer que c'est $X_0$.

L'inclusion $A(4)\subset A(q)$ 
induit une application birationnelle $\phi: X_0\dasharrow X(q)$.
Une section g\'en\'erique $X_0\cap \P^4$ de $X_0$ est param\'etr\'ee par des \'equations de la forme :
$$
t_j = \fr{T_j(\t)}{T_0(\t)}, \;\; j=1,2; \qquad s_j= \fr{S_j(\t)}{T_0(\t)^2}, \;\; j=1,\ldots,r-1,
$$
o\`u $\t\in \P^1$, les $T_j(\t)$ sont des polyn\^omes de degr\'e $\leq 2$ et les $S_j(\t)$
des polyn\^omes  de degr\'e $\leq 4$. Son image par $\phi$ est param\'etr\'ee par  
$$
x_{(i,j,\a)}(\t) = \fr{T_1(\t)^iT_2(\t)^jS_1(\t)^{\a_1}\ldots S_{r-1}(\t)^{\a_{r-1}}}{T_0(\t)^{i+j+2|\a|}},
\;\;\;
(i,j,\a)\in A(q). 
$$
C'est en g\'en\'eral une courbe rationnelle normale de degr\'e $q$.
La vari\'et\'e $X(q)$ v\'erifie donc la premi\`ere propri\'et\'e dans la D\'efinition \ref{D1}.
Pour montrer qu'elle appartient \`a $\cal{X}_{r+1,5}(q)$, il suffit de v\'erifier qu'elle engendre 
un espace de dimension $\pi_{r,5}(q)$.

\sk
On \'ecrit $q=2\s$ et on d\'ecompose $A(q) = \{(i,j,\a)\in \N^{r+1}, \;\; 1\leq i+j+2|\a| \leq \s \}$
selon les parit\'es de $i$ et de $j$ : le cardinal de $A(q)$ est la somme des cardinaux 
des quatre ensembles de $(i,j,\a)\in \N^{r+1}$ respectivement d\'efinis par :
$$
0<2i+2j+2|\a| \leq \s, \;\;\;  2i+2j+1+2|\a| \leq \s,\;\;\;   2i+2j+1+2|\a| \leq \s,\;\;\;   2i+2j+2+2|\a| \leq \s. 
$$

Si $\s=2\rho$ est pair, on obtient :  $\text{card}(A(q)) + 1 = {r+\rho+1 \choose r+1} + 3\, {r+\rho\choose r+1} = \pi_{r,5}(q) + 1$.

Si $\s=2\rho+1$ est impair, on a :  $\text{card}(A(q)) + 1 = 3 \, {r+\rho+1 \choose r+1} +  {r+\rho\choose r+1} = \pi_{r,5}(q) + 1$.

\nk
Ceci d\'emontre le troisi\`eme point du th\'eor\`eme.

\ssct{Si $X_0$ est un scroll rationnel normal}

On suppose enfin $n\geq 3$ et que $X_0$ est un scroll rationnel normal 
de dimension $r+1$ et de degr\'e $n-1$, d\'efini par des entiers $a_0,\ldots,a_r$ 
qui v\'erifient (\ref{S-a}).
Le scroll $X_0$ \'etant fix\'e \`a \'equivalence pr\`es, on associe \`a toute paire d'entiers 
$(\rho,\chi)$ tels que 
$$
\rho(n-1)+\chi \geq n-1, \qquad  \rho\geq 1, \qquad  -1\leq \chi\leq n-2,
$$
l'ensemble 
\beq
\label{AA}
A(\rho,\chi)=\Big\{(k,\a)\in (\N\times \N^r)\bck \{0\}, \;\;\;  |\a|\leq \rho, \;\; k \leq (\rho-|\a|)a_0 + \sum_{j=1}^r \a_j a_j + \chi\Big\}
\eeq
et la vari\'et\'e mon\^omiale $X(\rho,\chi)$ d\'efinie par cet ensemble.

\sk
Remarquons que la premi\`ere condition sur la paire $(\rho,\chi)$
exclut la paire $(1,-1)$. D'autre part, l'\'ecriture $q=\rho(n-1)+\chi$ est la division 
euclidienne de $q$ par $n-1$ sauf si $\chi=-1$, auquel cas $q$ est congru 
\`a $-1$ modulo $n-1$.

\sk
La vari\'et\'e $X(1,0)$ est le scroll rationnel normal qu'on a not\'e  $S_{a_0,\ldots,a_r}$
dans la Section~5.2. On peut supposer $X_0 = X(1,0)$. La vari\'et\'e 
$X_0\cap \C^{r+n-1}$ est param\'etr\'ee par 
\beq
\label{param}
(t,\ldots,t^{a_0}, s_1,s_1t, \ldots, s_1t^{a_1}, \ldots, s_r,s_rt,\ldots,s_rt^{a_r}),
\qquad (t,s)\in \C\times \C^r,
\eeq
et une section de $X_0$ par un $\P^{n-1}\subset \P^{r+n-1}$ g\'en\'erique est donn\'ee par les relations 
\beq
\label{TTT}
s_k = \fr{P_k(t)}{P_0(t)}, \qquad k=1,\ldots,r,
\eeq
entre les param\`etres $t$ et $s$, 
o\`u les $P_k(t)$ sont des polyn\^omes de degr\'e respectif $\leq n-1-a_k$.

\sk
Consid\'erons une vari\'et\'e $X(\rho,\chi)$. Comme $A(1,0)\subset A(\rho,\chi)$,
on a un isomorphisme naturel de $X_0\cap \C^{r+n-1}$ sur 
$X(\rho,\chi)\cap \C^{{\rm card}\, A(\rho,\chi)}$,
l'inverse d'une projection.
L'image d'une section de $X_0$ par un $\P^{n-1}\subset \P^{r+n-1}$ g\'en\'erique 
donn\'ee par (\ref{TTT}) est une courbe rationnelle normale de degr\'e $q = \rho(n-1) + \chi$.
En effet, elle est donn\'ee par 
$$
x_{(k,\a)}(t) = \fr{t^k P_1(t)^{\a_1}\cdots P_r(t)^{\a_r}P_0(t)^{\rho - |\a|}}{ P_0(t)^\rho },
\qquad 
(k,\a)\in A(\rho,\chi).
$$
Le degr\'e du d\'enominateur est major\'e par $\rho(n-1-a_0)\leq q$ et celui du num\'erateur  
est major\'e par $k + \sum_{j=1}^r \a_j(n-1-a_j) + (n-1-a_0) (\rho - |\a|)\leq q$,
par construction.

\sk
La vari\'et\'e $X(\rho,\chi)$ a donc la premi\`ere propri\'et\'e dans la D\'efinition \ref{D1}.
Pour montrer qu'elle a aussi la deuxi\`eme, on calcule la dimension  $N(\rho,\chi)={\rm card}\,A(\rho,\chi)$
de l'espace qu'elle engendre. 
On pose $\a_0=\rho-|\a|$, $\wh{\a}=(\a_0,\a)\in \N^{r+1}$ et $\wh{\a}\cdot a=\sum_{\mu=0}^r \a_\mu a_\mu$.
Soit $\tau$ une permutation circulaire de $\{0,\ldots,r\}$. Le nombre $N(\rho,\chi)+1$ est \'egal \`a 
\beqn
\sum_{|\wh{\a}|=\rho} \, ( \wh{\a}\cdot a + \chi + 1 ) 
 & = &
\fr{1}{r+1} \sum_{|\wh{\a}|=\rho} \, \sum_{j=0}^r ( \a_{\tau^j(0)} a_0 + \cdots + \a_{\tau^j(r)} a_r + \chi + 1)
\\
& = &
\fr{1}{r+1} \sum_{|\wh{\a}|=\rho} \, \rho(n-1) 
+
(\chi+1) \sum_{|\wh{\a}|=\rho} 1   \\
& = & 
(n-1) {r+\rho \choose r+1} + (\chi+1) {r+\rho \choose r},
\eeqn
ou encore 
$$
N(\rho,\chi)+1 = (\chi + 1){r+\rho + 1 \choose r+1}  +  (n-2-\chi){r+\rho \choose r+1}  =  \pi_{r,n}(q) + 1.
$$
La vari\'et\'e $X(\rho,\chi)$ appartient donc \`a la classe $\cal{X}_{r+1,n}(q)$ avec $q=\rho(n-1)+\chi$.
On a obtenu la premi\`ere partie du lemme suivant.
\ble
\label{scroll-1}
On suppose que $X_0$ est le scroll rationnel normal $S_{a_0,\ldots,a_r}$, de degr\'e $n-1$.
Soit $q \geq n-1$ un entier et $q=\rho(n-1) + \chi$
la division euclidienne de $q$ par $n-1$.

La vari\'et\'e mon\^omiale $X(\rho,\chi)$ d\'efinie par l'ensemble (\ref{AA})
est un \'el\'ement standard de la classe $\cal{X}_{r+1,n}(q)$, associ\'e \`a $X_0$. 

Si $\chi=n-2$, la vari\'et\'e $X(\rho+1,-1)$ a aussi cette propri\'et\'e ;
elle n'est pas \'equivalente \`a la vari\'et\'e $X(\rho,n-2)$ sauf si $n=3$ ou si $a_0=n-1$.
\ele
\bpf
Il reste \`a montrer la derni\`ere partie de l'\'enonc\'e. 
Commen\c{c}ons par les cas particuliers.

\sk
Si $a_0=n-1$, la description de $A(\rho,-1)$ dans (\ref{AA}) s'\'ecrit :
$$
|\a|\leq \rho, \qquad k\leq (\rho-|\a|)(n-1) - 1.
$$
Comme ce syst\`eme n'a pas de solution $(k,\a)\in \N^{r+1}$ avec $|\a|=\rho$, on a en fait 
$|\a|=\rho-1$ et $k\leq (\rho-1-|\a|)(n-1)+n-2$. Autrement dit,
$A(\rho,-1)=A(\rho-1,n-2)$ et  donc $X(\rho,-1)=X(\rho-1,n-2)$.

\sk
Si $n=3$ et $a_0=a_1=1$, la description de $A(\rho,-1)$ dans (\ref{AA})
s'\'ecrit :
$$
|\a|\leq \rho, \qquad k \leq \rho-(|\a|-\a_1) - 1.
$$
\'Echangeons les param\`etres $s_1$
et $t$ et les indices $k$ et $\a_1$. On obtient que la vari\'et\'e $X(\rho,-1)$ est \'equivalente 
\`a	 la vari\'et\'e mon\^omiale d\'efinie 
par l'ensemble des $(k,\a)\in \N^{r+1}\bck\{0\}$  qui v\'erifient 
$k \leq (\rho - 1 - |\a|) + \a_1 + 1$ et $|\a| \leq \rho - 1$,
c'est-\`a-dire \`a la vari\'et\'e $X(\rho-1,1)$.

\sk
Avant de traiter le cas g\'en\'eral,  consid\'erons une vari\'et\'e mon\^omiale $X(\rho,\chi)$.
Rappelons que son intersection avec $\C^{N(\rho,\chi)}$ est homog\`ene. 
Il est clair qu'elle n'est pas $(\rho+1)$-r\'eguli\`ere. Elle est $\rho$-r\'eguli\`ere si 
$k+|\a|\leq \rho$ implique la deuxi\`eme in\'egalit\'e dans (\ref{AA}),
c'est-\`a-dire si 
$$
|\a|\leq \rho \; \Rightarrow \; (\rho-|\a|)(a_0-1) + \sum_{j=1}^r \a_j a_j + \chi \geq 0.
$$
C'est \'evidemment le cas si $\chi\geq 0$. Alors, en faisant $t=0$
dans les mon\^omes $t^ks^\a$, $\, (k,\a)\in A(\rho,\chi)$, on voit 
que l'intersection de $X(\rho,\chi)$ avec son espace osculateur \`a l'ordre $\rho$
en $0$ contient une vari\'et\'e 
de Veronese de dimension $r$ et d'ordre $\rho$. 

\sk
On exclut maintenant les cas particuliers d\'ej\`a trait\'es. 
On suppose 
donc $a_2\geq 1$ ou $a_2=0$, $a_1\geq 1$ et $a_0\geq 2$.

Soit $\rho\geq 2$. Compte tenu de ce qu'on vient de voir, pour
montrer que les vari\'et\'es
standards $X(\rho,-1)$ et $X(\rho-1,n-2)$ ne sont pas \'equivalentes,
il suffit de montrer que l'intersection $Z$ de $X(\rho,-1)$ avec son osculateur 
\`a l'ordre $\rho-1$ en $0$ n'a pas de composante de dimension $r$.
On note $x(t,s)\in X(\rho,-1)$ le point de param\`etre $(t,s)\in \C^{r+1}$.

\sk
Si $x(t,s)\in Z$, tous les mon\^omes $t^ks^\a$ de degr\'e $\geq\rho$
avec $(k,\a)\in A(\rho,-1)$ sont nuls.
En particulier $s^\a=0$ si $|\a|=\rho\,$ et $\,\sum_{j=1}^r \a_ja_j-1\geq 0$.
Si $a_2\geq 1$, on obtient $s_1=s_2=0$, ce qui implique que $Z$ est de codimension $\geq 2$.

Si $a_2=0$, on a $a_1\geq 1$ et $a_0\geq 2$ par hypoth\`ese et on obtient encore $s_1=0$. On doit aussi annuler le mon\^ome $t^{\rho a_0 -1}$ 
qui est de degr\'e $\geq\rho$, ce qui donne $t=0$. On obtient encore que 
$Z$ est de codimension $\geq 2$. Le lemme est d\'emontr\'e.
\epf 

\ssct{Si $X_0$ est un scroll rationnel normal, suite et fin}

Il reste \`a montrer que les
vari\'et\'es standards de la classe $\cal{X}_{r+1,n}(q)$
qui apparaissent dans le Lemme \ref{scroll-1} sont les seules 
associ\'ees au scroll rationnel normal $X_0=S_{a_0,\ldots,a_r}$.

\sk
Si $a_0=n-1$, c'est facile. Dans ce cas, $X_0$ est un c\^one au-dessus d'une courbe rationnelle normale
de degr\'e $n-1$, donc le groupe $\text{Pic}\,(X_0)$ est le groupe libre engendr\'e par la classe $R$
du c\^one au-dessus d'un point de cette courbe et $(n-1)R$ est 
la classe d'une section hyperplane de $X_0$. 
Pour $q\geq 1$, on a $C\cdot qR = q$ si $C$ est une section de $X_0$ par un $\P^{n-1}\subset \P^{r+n-1}$
g\'en\'erique. 
Il en r\'esulte que pour tout $q\geq n-1$, il existe {\em au plus} 
une vari\'et\'e $X\in \cal{X}_{r+1,n}(q)$ associ\'ee \`a $X_0$, \`a \'equivalence pr\`es. 
C'est la vari\'et\'e donn\'ee par le Lemme \ref{scroll-1}.

\bk
On suppose maintenant $a_0\neq n-1$. On dit alors que $X_0$ est {\em un scroll g\'en\'eral}.
Le calcul sera plus laborieux. 
Le groupe $\text{Pic}\,(X_0)$ est le $\Z$-module libre $\Z H \oplus \Z R$ engendr\'e 
par la classe $H$ d'une section hyperplane de $X_0$ et la classe $R$ d'un \'el\'ement
du r\'eglage  de $X_0$ par des $\P^r$, dont on obtient un \'el\'ement en fixant 
le param\`etre $t$ dans (\ref{param}).

Si $a\in \Z$, on note $a^+ = \max (a,0)$. Pour tout  $(\rho,\chi)\in \N^\star \times \Z$,
on d\'efinit l'entier $I(\rho,\chi)$ par la formule 
\beq
\label{I}
I(\rho,\chi) = \sum_{|\a|=\rho} (\a_0 a_0 + \cdots + \a_r a_r + \chi + 1)^+.
\eeq
La dimension du syst\`eme lin\'eaire $|\rho H + \chi R|$ est donn\'ee 
par la formule suivante si $\rho \geq 1$, voir par exemple Harris \cite{Ha} :
$$
\dim \, |\rho H + \chi R| + 1 = I(\rho,\chi), \qquad \rho\geq 1.
$$

\sk
Soit $\phi: X_0\dasharrow X$ une application birationnelle
qui associe une vari\'et\'e standard $X\in \cal{X}_{r+1,n}(q)$ \`a $X_0$.
Soir $|\rho H + \chi R|$ le syst\`eme lin\'eaire complet qui d\'efinit $\phi$.

Si $C$ est une section de $X_0$ par un $\P^{n-1}\subset \P^{r+n-1}$ 
g\'en\'erique, on a $C\cdot H = n-1$ et $C\cdot R = 1$, ce qui donne $\rho (n-1)+\chi=q$.

Si $l$ est une droite contenue dans un 
repr\'esentant de $R$, on a $l\cdot H=1$ et $l\cdot R=0$ donc $l \cdot (\rho H+\chi R) = \rho(n-1)$,
ce qui donne $\rho \geq 1$. On a donc :
\ble
\label{CNS}
Si le syst\`eme lin\'eaire $|\rho H  + \chi R|$
d\'efinit une application rationnelle qui associe au scroll g\'en\'eral $X_0$ une vari\'et\'e 
$X\in \cal{X}_{r+1,n}(q)$, la paire $(\rho,\chi)\in \Z^2$ v\'erifie : 
\beq
\label{cond}
\rho\geq 1, \qquad q = \rho (n-1) + \chi, \qquad I(\rho,\chi) = \pi_{r,n}(q)+1.
\eeq
\ele
On a calcul\'e $I(\rho,\chi)$ pour $\rho \geq 1$ et  $\chi\geq -1$
au d\'ebut de la Section 5.4 et obtenu que 
$$
I(\rho,\chi) = (\chi + 1){r+\rho + 1 \choose r+1}  +  (n-2-\chi){r+\rho \choose r+1}, \qquad \chi\geq -1,
$$
ne d\'epend pas des entiers $a_0,\ldots,a_r$ mais seulement de leur somme $n-1$.
En particulier, si $q=\rho(n-1)+\chi$ est la division euclidienne de 
$q$ par $(n-1)$, alors $(\rho,\chi)$ est une solution de (\ref{cond}) et,
si $q=\rho(n-1)+n-2$, $\,(\rho+1,-1)$ est une autre solution de (\ref{cond}).

\sk
Pour montrer que la liste du Lemme \ref{scroll-1} est exhaustive,
ce qui ach\`evera la d\'emons\-tration du Th\'eor\`eme \ref{class},
il reste deux choses \`a v\'erifier.
D'une part, on doit v\'erifier que les solutions pr\'ec\'edentes du syst\`eme (\ref{cond}) sont les seules, 
autrement dit  que (\ref{cond}) implique $\chi \in \{-1,\ldots,n-2\}$. C'est l'objet 
du lemme suivant.

D'autre part on doit v\'erifier, si $n=3$ et $a_0=a_1=1$, que 
les deux syst\`emes lin\'eaires $|\rho H + R|$ et 
$|(\rho+1)H - R|$ d\'efinissent des applications 
rationnelles qui envoient $X_0$ sur des vari\'et\'es \'equivalentes,
ce qu'on fait maintenant.

Dans ce cas, $X_0$ est un c\^one au-dessus d'une quadrique lisse $Q\subset \P^3$
et le groupe ${\rm Pic}(X_0)$ est aussi le groupe libre 
engendr\'e par les classes $R$ et $R'$ des deux r\'eglages de $X_0$ par des 
$\P^r$. De plus $R+R'=H$. On a donc 
$\rho H + R  = (\rho + 1) R + \rho R'$ et $(\rho +1) H - R = \rho R + (\rho + 1)R'$.
Comme les deux r\'eglages $R$ et $R'$ sont \'echang\'es par un automorphisme 
de $X_0$, on obtient le r\'esultat cherch\'e.

\sk
Le lemme suivant termine la discussion.
\ble
\label{dernier point}
Si $(\rho,\chi)$ est une solution du syst\`eme (\ref{cond}), alors 
$\chi\in \{-1,\ldots,n-2\}$.
\ele
\bpf
Posons :
$$
I_0(\rho,\chi) = \sum_{|\a|=\rho} (\a_0 (n-1) + \chi + 1)^+.
$$
Si $\rho\geq 2$, on \'ecrit $\a=(\a_0,\a')$, on distingue selon
que $\a_0$ est nul ou pas et on calcule 
\beqn
I_0(\rho,\chi) 
& = & \sum_{|\a'|=\rho} (\chi + 1)^+ + \sum_{|\a|=\rho,\; \a_0\geq 1} (\a_0 (n-1) + \chi + 1)^+
\\
& = & \sum_{|\a'|=\rho} (\chi + 1)^+ + \sum_{|\a|=\rho-1} ((\a_0+1)(n-1) + \chi + 1)^+
\\
& = & \sum_{|\a'|=\rho} (\chi + 1)^+ + I_0(\rho-1,\chi+n-1).
\eeqn
Ceci montre que, pour $\rho(n-1)+\chi=q$ fix\'e, $I_0(\rho,\chi)$ atteint son maximum
$\pi_{r,n}(q)+1$ en $(\rho,\chi)$ si et seulement si $\chi\leq n-2$.

\sk
Comme $I(\rho,\chi) = I_0(\rho,\chi)$ si $\chi\geq -1$, on obtient d\'ej\`a $\chi\leq n-2$
si $(\rho,\chi)$ est solution de (\ref{cond}). Finalement pour montrer qu'on a 
$\chi\geq -1$, il suffit de montrer que :
\beq
\label{lafin}
\chi<-1 \; \Rightarrow \; I(\rho,\chi) < I_0(\rho,\chi).
\eeq
On  montre que $I(\rho,\chi)$ diminue si, \'etant donn\'e deux indices
distincts $j,k\in \{0,\ldots,r\}$ tels que $a_j\geq a_k\geq 1$, on substitue 
la paire $(a_j+1,a_k-1)$ \`a la paire $(a_j,a_k)$ dans $(a_0,\ldots,a_r)$.
Par sym\'etrie, il suffit de traiter le cas $j=0$, $k=1$.

\sk
On \'ecrit que $I(\rho,\chi)$ est une somme de termes de la forme :     
$$
\sum_{\a_0+\a_1=\mu} (\a_0a_0 + \a_1a_1  + h)^+,\qquad \mu\in \N,\;\; h\in \Z.
$$
Si $\a_0=\a_1$ le terme correspondant de la somme ci-dessus ne d\'epend que de $a_0+a_1$. 
On regroupe les autres termes par paires, soit avec $i>j$ et $i+j=\mu$ :
$$
(ia_0 + ja_1  + h)^+  + (ja_0 + ia_1  + h)^+ = A^+ + B^+, \qquad A\geq B.
$$
D'autre part :
$$
(i(a_0+1) + j(a_1-1)  + h)^+  + (j(a_0+1) + i(a_1-1)  + h)^+ = (A+(i-j))^+ + (B-(i-j))^+.
$$
Il suffit de remarquer que, si $A,B,C\in \Z$,
$$
A\geq B, \;\; C\geq 0 \; \Rightarrow (A+C)^+ + (B-C)^+ \geq A^+ + B^+,
$$
pour obtenir le r\'esultat en vue.

\sk
De proche en proche, on est ramen\'e \`a d\'emontrer (\ref{lafin}) 
pour $a_0=n-2$, $a_1=1$. Alors :
$$
I(\rho,\chi) = \sum_{|\a|=\rho} (\a_0 (n-2) + \a_1 + \chi + 1)^+.
$$
Compte tenu de ce qui pr\'ec\`ede, il suffit de montrer, en consid\'erant seulement les 
paires $(\a_0,\a_1)\in \{(\rho,0),(0,\rho)\}$, qu'on a :
$$
(\rho(n-2) + \chi + 1)^+ + (\rho + \chi + 1)^+  < (\rho(n-1)  + \chi + 1)^+ + (\chi + 1)^+
$$
si $\chi<-1$, ce qui est \'evident, compte tenu du fait que $\rho(n-1)  + \chi = q >0$.
\epf

\sct{Exemples de vari\'et\'es sp\'eciales}

\ssct{Introduction}

D'apr\`es le Th\'eor\`eme \ref{Th4}, seules les classes  $\cal{X}_{r+1,n}(2n-3)$ sont 
susceptibles de contenir des vari\'et\'es sp\'eciales, et cela seulement
si $r\geq 2$ et $n\geq 3$, ce qu'on suppose maintenant. 

\sk
Commen\c{c}ons par rappeler la classification des vari\'et\'es standards 
$X\in \cal{X}_{r+1,n}(2n-3)$. Compte tenu du Th\'eor\`eme \ref{class}, une telle vari\'et\'e 
est associ\'ee \`a un scroll rationnel normal 
$S_{a_0,\ldots,a_r}$ de degr\'e $n-1=a_0+\cdots+a_r$. 

Pour un scroll donn\'e, 
il y a en g\'en\'eral deux vari\'et\'es standards qui lui sont associ\'ees, \`a \'equivalence pr\`es.
L'une est le scroll $S_{a_0 + n-2, \ldots, a_r + n-2}$, dont l'intersection avec $\C^{(r+2)(n-1)-1}$
est param\'etr\'ee par 
\beq
\label{ST1}
(t,\ldots,t^{a_0+n-2},s_1,ts_1,\ldots,t^{a_1+n-2}s_1,\ldots,s_r,ts_r,\ldots,t^{a_r+n-2}s_r), 
\eeq
o\`u $(t,s)\in \C\times \C^r$.
C'est la seule solution si $a_0=n-1$ ou si $n=3$. Sinon, l'autre solution est \'equivalente \`a la vari\'et\'e 
dont l'intersection avec $\C^{(r+2)(n-1)-1}$ est param\'etr\'ee par 
$$
(t,\ldots,t^{2a_0-1}, \ldots , s_i,ts_i,\ldots,t^{a_0+a_i-1}s_i, \ldots, s_js_k, ts_js_k,\ldots,t^{a_j+a_k-1}s_js_k,\ldots),
$$
o\`u $(t,s)\in \C\times \C^r$ et $i,j,k\in \{1,\ldots,r\}$,
avec $j\leq k$ et $a_j+a_k\geq 1$.

\sk
Dans l'attente de r\'esultats plus g\'en\'eraux, on pr\'esente dans ce court 
chapitre quelques exemples de vari\'et\'es sp\'eciales.

\ssct{Quelques vari\'et\'es sp\'eciales dans les classes $\cal{X}_{r+1,3}(3)$ et $\cal{X}_{r+1,4}(5)$}

On a le r\'esultat suivant :
\bpr
\label{3-3}
Soit  $\cal{Q}$ une quadrique de $\P^{r+1}$, de rang $\mu\geq 3$.
L'image de $\P^1\times \cal{Q}$ par le plongement de Segre de type $(1,1)$
est une vari\'et\'e de la classe $\cal{X}_{r+1,3}(3)$. Elle est sp\'eciale si $r\geq 2$.
Ces vari\'et\'es sont class\'ees \`a \'equivalence pr\`es par leur dimension et le rang $\mu\in \{3,\ldots,r+2\}$
de la quadrique $\cal{Q}$.
\epr
Dans le cas $r=2$ et $\mu=4$, la vari\'et\'e obtenue est \'equivalente \`a 
l'image de $\P^1\times \P^1\times \P^1$ par 
le plongement de Segre de type $(1,1,1)$. Nous devons cet exemple,
le premier d'une vari\'et\'e 
sp\'eciale dont nous ayons eu connaissance, \`a F. Russo \cite{Ru}. 
\bpf
Le plongement de Segre $\s: \P^1\times \P^{r+1}\rightarrow \P^{2r+3}$ de type $(1,1)$ est d\'efini par 
$$
\s: \; ([S_0:S_1],[T_0:\cdots:T_{r+1}]) \mapsto [S_0T_0:\cdots:S_0T_{r+1}:S_1T_0:\cdots:S_1T_{r+1}].
$$
Soit $\cal{Q}\subset \P^{r+1}$ une quadrique de rang $\mu \geq 3$,
autrement dit irr\'eductible, et $X\subset \P^{2r+3}$
l'image de $\P^1\times \cal{Q}$ par le plongement $\s$.

\sk
Soit $(\tau_1,q_1),(\tau_2,q_2),(\tau_3,q_3)$ trois points 
de $\P^1\times \cal{Q}$ tels que $\tau_1,\tau_2,\tau_3\in \P^1$
soient deux-\`a-deux distincts et que $q_1,q_2,q_3\in \cal{Q}$ engendrent un $\P^2$
qui coupe $\cal{Q}$ suivant une conique propre $\Gamma$.
Si $\phi: \, \P^1 \rightarrow \Gamma$ est l'isomorphisme d\'etermin\'e par 
$\phi(\tau_i)=q_i$, $\,i=1,2,3$, l'application $\P^1\rightarrow \P^{2r+3}$
d\'efinie par $\tau\mapsto \s(\tau,\phi(\tau))$ 
param\`etre une courbe rationnelle normale de degr\'e $3$ contenue dans 
$X$. Comme $X$ engendre l'espace $\P^{2r+3}$, on obtient que $X$
est une vari\'et\'e de la classe $\cal{X}_{r+1,3}(3)$.

\sk
Pour montrer que $X$ est une vari\'et\'e sp\'eciale, il est commode de travailler dans $\C^{2r+3}$. 
On peut supposer que 
$X\cap \C^{2r+3}$ est param\'etr\'e par :
$$
x(s,t) = (t,s,ts,q(s),tq(s)), \qquad  t\in\C, \;\; s\in \C^r,
$$
o\`u $q$ est une forme quadratique de rang $\mu-2\geq 1$. 
Si la vari\'et\'e $X$ n'est pas une vari\'et\'e mon\^omiale au sens 
de la D\'efinition \ref{var-mon}, elle a toutefois la propri\'et\'e que $X\cap \C^{2r+3}$
est homog\`ene. En effet, si $(t^\star,s^\star)\in \C\times \C^r$,
les  composantes de $x(t+t^\star,s+s^\star) - x(t^\star,s^\star)$
engendrent le m\^eme espace vectoriel de polyn\^omes que celles de $x(t,s)$.

\sk
Consid\'erons alors l'intersection de $X$ avec son espace osculateur $X_0(1)$,
dont la trace dans $\C^{2r+3}$ est donn\'ee par les relations suivantes entre les param\`etres $t\in \C$ et $s\in \C^r$ :
$$
ts = 0, \;\; q(s)=0.
$$
C'est la r\'eunion d'un $\P^1$ et d'une quadrique de dimension $r-1$ 
et de rang $\mu-2$. D'autre part, (\ref{ST1}) rappelle que si $\tilde{X}\in \cal{X}_{r+1,3}(3)$
est une vari\'et\'e sandard, $\tilde{X}_a(1)\cap \tilde{X}$ est de dimension $r$ pour $a\in \tilde{X}$ g\'en\'erique.
Ceci suffit pour conclure que la vari\'et\'e $X$ est sp\'eciale si $r>1$ et que deux quadriques 
de rangs diff\'erents d\'efinissent des vari\'et\'es qui ne sont pas \'equivalentes.
\epf

\sk
Les vari\'et\'es pr\'ec\'edentes admettent des param\'etrages homog\`enes de la forme :
$$
X(T_0,T_1,S) = [T_0^3,T_0^2T_1,T_0^2S,T_0T_1S,T_0q(S),T_1q(S)],
$$
o\`u $S=[S_1:\ldots:S_r]$ et $q(S)$ est une forme quadratique de rang $\mu' = \mu - 2\geq 1$.
Les composantes de $X(T_0,T_1,S)$ s'annulent \`a l'ordre $2$ pour 
$$
T_0=0, \;\; T_1S=0, \;\; q(S)=0,
$$
c'est-\`a-dire au point $p=[0:1:0:\cdots:0]$ et le long de la quadrique $\cal{Q}'$
du $\P^{r-1}$ d'\'equations $T_0=0,\, T_1=0$, d\'efinie par $T_0=T_1=0$ et $q(S)=0$.

\sk
Cette remarque sugg\`ere de consid\'erer l'espace des polyn\^omes homog\`enes de degr\'e $3$
qui s'annulent \`a l'ordre $2$ le long de $\cal{Q}'$.
Pour $\mu'\geq 2$, on obtient le r\'esultat suivant.
\bpr
\label{4-5}
On suppose $r\geq 2$ et l'on se donne un  $\P^{r-1}\subset \P^{r+1}$ et une quadrique $\cal{Q}'$ de $\P^{r-1}$, de rang $\mu'\geq 2$.
Soit $\phi : \P^{r+1}\dasharrow \P^{3r+5}$ une application rationnelle 
d\'efinie par le syst\`eme lin\'eaire des hypersurfaces de degr\'e $3$ de $\P^{r+1}$ qui ont des 
points doubles le long de $\cal{Q}'$. Son image $X=\phi(\P^{r+1})$  est une vari\'et\'e
sp\'eciale de la classe $\cal{X}_{r+1,4}(5)$.
Les vari\'et\'es ainsi obtenues sont class\'ees \`a \'equivalence pr\`es par leur dimension et le rang $\mu'\in \{2,\ldots,r\}$
de la quadrique $\cal{Q}'$.
\epr
\bpf
On reprend les notations des consid\'erations qui pr\'ec\`edent l'\'enonc\'e.
Si $\mu'\geq 2$, un polyn\^ome homog\`ene de degr\'e $3$, qu'on \'ecrit 
\sk
$$
F(T_0,T_1,S) = P(T_0,T_1) + \sum_{j=1}^r P_j(T_0,T_1)S_j + T_0 R_0(S) + T_1R_1(S) + R(S),
$$
s'annule \`a l'ordre $2$ (au moins) le long de $\cal{Q}'$ si et seulement s'il est de la forme 
\beq
\label{hypo-spe1}
F(T_0,T_1,S) = P(T_0,T_1) + \sum_{j=1}^r P_j(T_0,T_1)S_j + (c_0T_0+c_1T_1)q(S).
\eeq
Le syst\`eme lin\'eaire introduit dans l'\'enonc\'e est donc de dimension 
$3r+5=\pi_{r,4}(5)$. 

\sk
D'autre part, \'etant donn\'e un $4$-uplet {\em g\'en\'erique} $(p_i,p_2,p_3,p_4)$
de points de $\P^{r+1}$, ces points engendrent un $\P^3$ qui coupe la quadrique $\cal{Q}'$
en deux points $q_1,q_2$, tels que $p_1,p_2,p_3,p_4,q_1,q_2$ soient six points en position g\'en\'erale 
dans ce $\P^3$. 
Il existe une et une seule cubique gauche qui passe par ces 
six points. On v\'erifie ais\'ement que son image par une application 
rationnelle $\phi : \P^{r+1}\dasharrow \P^{3r+5}$, d\'efinie par le syst\`eme lin\'eaire 
consid\'er\'e, est une courbe rationnelle normale de degr\'e $5$.

Il en r\'esulte que par les points $\phi(p_1),\ldots,\phi(p_4)$ de la vari\'et\'e $X=\phi(\P^{r+1})$
passe une courbe rationnelle normale de degr\'e $5$, puis que $X$ est une vari\'et\'e de la classe $\cal{X}_{r+1,4}(5)$.

\sk
On peut supposer que $X\cap \C^{3r+5}$ est param\'etr\'e par :
\beq
\label{hypo-spe2}
(t,t^2,t^3,s,ts,t^2s,q(s),tq(s)), \qquad t\in \C, \;\; s\in \C^r.
\eeq
C'est une sous-vari\'et\'e homog\`ene de $\C^{3r+5}$ et la trace sur $\C^{3r+5}$
de l'intersection $X\cap X_0(1)$ est obtenue quand les param\`etres $t\in \C$ et $s\in \C^r$ 
v\'erifient $t=0$ et $q(s)=0$. C'est une quadrique de dimension $r-1$ 
et de rang $\mu'$. Ceci montre que deux quadriques $\cal{Q}'$
de rangs diff\'erents d\'efinissent des vari\'et\'es qui ne sont pas \'equivalentes.

\sk
Enfin, pour montrer que les vari\'et\'es obtenues sont sp\'eciales, on peut remarquer 
que, {\em par construction}, si $a\in X$ est g\'en\'erique, l'image de $X$
par la projection de centre $X_a(1)$ est une des vari\'et\'es consid\'er\'ees 
dans l'\'enonc\'e pr\'ec\'ecent, associ\'ee \`a une quadrique $\cal{Q}$ de rang $\mu'+2\geq 4$.
On obtient ainsi toutes ces vari\'et\'es sauf celles qui sont associ\'ees \`a une quadrique de 
rang $3$.
\epf
\bre
Si $\mu'=1$, le param\'etrage (\ref{hypo-spe2}) d\'efinit encore une vari\'et\'e sp\'eciale
de la classe $\cal{X}_{r+1,4}(5)$. Dans ce cas, on modifie la d\'efinition 
du syst\`eme lin\'eaire introduit dans l'\'enonc\'e. On consid\`ere \`a la place 
le syst\`eme des hypersurfaces de degr\'e $3$ d\'efinies par une \'equation $F(T_0,T_1,S)=0$,
o\`u $F(T_0,T_1,S)$ est de la forme (\ref{hypo-spe1}).
On reprend la d\'emonstration pr\'ec\'edente en associant \`a un $4$-uplet {\em g\'en\'erique} 
de $\P^{r+1}$ le point d'intersection $q$ du $\P^3$ qu'il engendre avec le $\P^{r-2}$
\og sous-jacent \fg\, \`a la quadrique $\cal{Q}'$ (qui est de rang $1$)
et la droite intersection 
de ce $\P^3$ avec le $\P^{r-1}$ d'\'equations $T_0=T_1=0$. Une application rationnelle
d\'efinie par le syst\`eme lin\'eaire consid\'er\'e envoie la cubique gauche 
qui passe par les points $p_1,p_2,p_3$ et $p_4$ et qui est tangente \`a cette droite 
au point $q$ sur une courbe rationnelle normale de degr\'e $5$ et $\P^{r+1}$.
On montre ainsi que (\ref{hypo-spe2}) d\'efinit une vari\'et\'e sp\'eciale.
\ere

\ssct{Une vari\'et\'e sp\'eciale de la classe $X_{3,6}(9)$ et une de la classe $\cal{X}_{3,5}(7)$}

L'exemple suivant est curieux :
\bpr
\label{V33}
La vari\'et\'e de Veronese de dimension $3$ et d'ordre $3$ est une vari\'et\'e standard 
de la classe $\cal{X}_{3,2}(3)$ et une vari\'et\'e sp\'eciale de la classe $\cal{X}_{3,6}(9)$.

En projetant cette vari\'et\'e $V$ depuis un espace osculateur $V_a(1)$, on obtient une vari\'et\'e sp\'eciale 
$X$ de la classe $X_{3,5}(7)$, ind\'ependante, \`a \'equivalence pr\`es, du choix de $a\in V$. 
\epr
\bpf
On sait d\'ej\`a que $V$ est un \'el\'ement standard de $\cal{X}_{3,2}(3)$. Cette vari\'et\'e $V$ engendre un espace 
de dimension $\pi_{2,2}(3)= 19$ et on a aussi $\pi_{2,6}(9)=19$.

Soit $v: \P^3\rightarrow \P^{19}$ un plongement de Veronese d'ordre $3$, d'image $V$.
Par six points en position g\'en\'erale dans $\P^3$ passe une 
unique cubique gauche et son image par $v$ est une courbe rationnelle de degr\'e $9$. La vari\'et\'e $V$
appartient donc \`a la classe $\cal{X}_{3,6}(9)$.
Elle est sp\'eciale dans cette classe, puisqu'elle est $3$-r\'eguli\`ere 
et que les vari\'et\'es standards de cette classe ne le sont pas.

La deuxi\`eme partie de l'\'enonc\'e r\'esulte de la premi\`ere et du Th\'eor\`eme \ref{heredite2}
\epf

{\small

\sct{Appendice : propri\'et\'es des espaces osculateurs}

On rappelle dans cet appendice quelques propri\'et\'es des espaces osculateurs 
qu'on utilise dans le corps de l'article, en g\'en\'eral sans r\'ef\'erence. 
Elles sont toutes \'el\'ementaires mais elles ne sont pas n\'ecessairement famili\`eres 
au lecteur. On esquisse quelques d\'emonstrations.

\sk
Soit $X$ un germe en $x_0\in \C^N$ de vari\'et\'e analytique  {\em lisse} 
de dimension $d\geq 1$ et $k\geq 0$ un entier.
Soit $v: (\C^d,p)\rightarrow (\C^N,x_0)$ 
un germe de param\'etrage (r\'egulier, c'est-\`a-dire de rang $d$) de $X$, avec $v(p)=x_0$.
On v\'erifie facilement que l'espace vectoriel engendr\'e par les vecteurs 
\beq
\label{D-osc1}
\fr{\pl^\a v}{\pl t^\a}(p), \qquad 1\leq |\a|\leq k,
\eeq
ne d\'epend pas du choix de $v$. 
{\em L'espace osculateur} (pour l'instant affine), 
ou plus simplement {\em l'osculateur}  $X_{x_0}(k)$ de $X$ \`a l'ordre $k$ en $x_0$  
est l'espace affine qui passe par le point $x_0$ et de direction cet espace vectoriel.

La dimension de $X_{x_0}(k)$  est \'evidem\-ment au plus \'egale
au cardinal de l'ensemble des multi-indices $\a=(\a_1,\ldots,\a_d)$
dont la longueur est comprise entre $1$ et $k$, ce qui donne :
\beq
\label{D-osc2}
\dim\, X_{x_0}(k) + 1 \leq {d+k \choose d}.
\eeq
Le second membre est aussi la dimension de l'espace des polyn\^omes de degr\'e $\leq k$ en $d$
variables ou celle de l'espace des polyn\^omes homog\`enes de degr\'e $k$
en $d+1$ variables.

On dit que le germe $X$ est {\em  $k$-r\'egulier} en $x_0$ si les deux membres de (\ref{D-osc2})
sont \'egaux, autrement dit si la famille (\ref{D-osc1}) est une famille libre.

\sk
L'osculateur $X_{x_0}(k)$ est engendr\'e par les osculateurs $C_{x_0}(k)$
des germes de courbes lisses $C\subset X$ en $x_0$.
Pour le montrer, on peut supposer $x_0=0$ et que le param\'etrage $v$ 
est d\'efini au voisinage de $p=0$, soit 
$v(t) = \sum_{|\a|=1}^{+\infty} t^\a v_\a$.
L'osculateur $X_0(k)$ est le sous-espace engendr\'e
par les vecteurs $v_\a$ avec $1\leq |\a|\leq k$.

Si $c\in \C^d\bck \{0\}$, l'image par $v$ du germe de droite param\'etr\'ee par
$\t \mapsto \t c$ est param\'etr\'ee par $\t\mapsto \sum_{|\a|=1}^{+\infty} \t^{|\a|} c^\a v_\a$.
L'espace engendr\'e par les osculateurs de ces courbes \`a l'ordre $k$ en $0$ contient
donc les vecteurs $\sum_{|\a|=j} c^\a v_\a$ avec $j\in \{1,\ldots,k\}$
et $c\in \C^r$. Il contient aussi leurs d\'eriv\'ees partielles par rapport \`a $c$,
donc la famille des vecteurs $v_\a$ avec $1\leq |\a| \leq k$, ce qui donne 
le r\'esultat.

\bk
La notion d'osculateur est une notion projective :  
si une homographie $\phi: \C^N \dasharrow \C^N$ est d\'efinie au voisinage de $x_0$,
l'osculateur \`a l'ordre $k$ du germe $\phi(X)$ en $\phi(x_0)$
est l'image par $\phi$ de celui de $X$ en $x_0$. 
La v\'erification directe pour une homographie g\'en\'erale peut-\^etre 
un peu compliqu\'ee. Elle est tr\`es simple dans le cas o\`u 
$\phi$ est un isomorphisme affine. On se ram\`ene ainsi 
au cas o\`u $\phi$ est tangente 
\`a l'identit\'e en $0\in \C^N$, {\em i.e.} est de la forme $\phi(x) = (1 + u(x))^{-1} \, x$,
o\`u $u$ est une forme lin\'eaire. La v\'erification est \`a nouveau facile.

Si $X\subset \P^N$ est un germe de vari\'et\'e lisse en $x_0\in \P^N$, on 
note maintenant $X_{x_0}(k)$ son espace projectif osculateur
\`a l'ordre $k\geq 0$ en $x_0$.

\sk
Soit $X$ une sous-vari\'et\'e alg\'ebrique irr\'eductible de $\P^N$.
Si $x\in X_{\rm reg}$, il est clair, \`a partir de la d\'efinition, qu'on a 
$\dim \,X_{x'}(k)\geq \dim \, X_{x}(k)$
pour $x'$ assez voisin de $x$. Soit $m$ le maximum de la dimension 
de $X_x(k)$ quand $x$ varie dans $X_{\rm reg}$.
L'ensemble $\{x\in X_{\rm reg}, \; \dim X_x(k)=m\}$
est un ouvert dense de $X$ et l'application $x\mapsto X_x(k)$
est analytique sur cet ouvert, \`a valeurs dans la grassmannienne des $m$-plans de $\P^N$.
La d\'emonstration est analogue \`a toute d\'emonstration d'un r\'esultat de ce type, quand
il est facile \`a d\'emontrer, ce qui est le cas ici. 

\bk
On rappelle maintenant quelques propri\'et\'es qui jouent un r\^ole 
important dans cet article et qui, pour cette raison, 
m\'eritent qu'on les \'enoncent.
On a d'abord la propri\'et\'e  suivante.
\ble
Soit $X$ un germe de vari\'et\'e lisse et $k$-r\'egulier en $x_0\in \P^N$.
Tout germe $Y\subset X$ de vari\'et\'e lisse de dimension $\geq 1$ en $x_0$ est $k$-r\'egulier.
\ele
\bpf
Soit $d\geq 1$ la dimension de $X$ et $1\leq d'\leq d$ celle de $Y$. On peut choisir la param\'etrisation 
$v: (\C^d,0) \rightarrow (X,x_0)$ de $X$ de 
telle fa\c{c}on que $Y$ soit l'image du germe en $0$ du sous-espace d'\'equations 
$t_j=0$, $\,j=d'+1,\ldots,d$. 
La  sous-famille des vecteurs  (\ref{D-osc1})
obtenue en se restreignant aux multi-indices $\a$ tels que $\a_j=0$
pour $j>s$ est libre, ce qui donne le r\'esultat.
\epf
Les espaces osculateurs se comportent comme on s'y attend par projection r\'eguli\`ere :
\ble
Soit $X$ un germe de vari\'et\'e lisse en $x_0\in \P^N$
et $\pi:\P^N \dasharrow \P^M$ une projection dont le centre $\Q$ est 
en somme directe projective avec $X_{x_0}(k)$.
Alors $\pi(X)_{\pi(x_0)}(k)=\pi(X_{x_0}(k))$.
En particulier, si le germe $X$ est $k$-r\'egulier en $x_0$, son image est $k$-r\'eguli\`ere
en $\pi(x_0)$.
\ele
\bpf
Une r\'ecurrence sur la dimension de $\Q$ permet de se ramener au cas 
o\`u $\Q$ est un point. On peut aussi supposer que 
$X_{x_0}(k)$ est contenu dans l'espace cible de la projection.

On se ram\`ene ainsi \`a  la situation suivante. La projection 
$\pi: \C^{N-1}\times \C  \dasharrow \C^{N-1}\times \{0\}$ est 
d\'efinie par $(x',x_N)\mapsto (x'/(1-x_N),0)$  et au voisinage de $0\in \C^{N-1}\times \C$,
le germe $X$ est donn\'e par un param\'etrage $v(t)=(v'(t),v_N(t))$
avec $v_N(t)=O(|t|^{k+1})$. Son image est alors donn\'ee 
par un param\'etrage $w(t)=(w'(t),0)$ avec $w'(t)=v'(t)+O(|t|^{k+1})$,
ce qui donne le r\'esultat.
\epf
Dans le cas o\`u $X$ est une courbe,  on \'etend le r\'esultat pr\'ec\'edent 
aux projections dont le centre rencontre cette courbe :
\ble
Soit $k\in \N^\star$ et $C$ un germe de courbe lisse $(k+1)$-r\'egulier en $a\in \P^N$.
Soit $\pi: \P^N \dasharrow H$ une projection de centre $a$
et de cible un hyperplan $H$ de $\P^N$ tel que $a\notin H$.
L'image par $\pi$ de $C\bck\{a\}$ se prolonge en un germe de courbe lisse $k$-r\'egulier en $a' = C_a(1)\cap H$
et l'osculateur \`a l'ordre $k$ de $C'$ en $a'$ est l'image 
par $\pi$ de l'osculateur \`a l'ordre $k+1$ de $C$ en~$a$.
\ele
\bpf
On peut supposer $a=0\in \C^N$, que $H$ est donn\'e par $x_1=1$ et $C$ par 
$$
x_j = x_1^j + O(x_1^{k+2}) \;\; \text{si $j=2,\ldots,k+1$}, \;\;\; x_j = O(x_1^{k+2}) \;\; \text{si $j=k+2,\ldots,N$}.
$$
L'image $C'$ de $C$ est alors donn\'ee par $x_1=1$ et pour $t\in \C$ voisin de $0$ :
$$
x_j = t^{j-1} + O(t^{k+1}) \;\; \text{si $j=2,\ldots,k+1$}, \;\;\; x_j = O(t^{k+1}) \;\; \text{si $j=k+2,\ldots,N$}.
$$
\epf
Encore quelques remarques \`a propos des courbes.
Soit $C\subset \P^N$ une courbe alg\'ebrique irr\'eductible et 
$k\in \N^\star$ sa {\em r\'egularit\'e osculatrice}.
On entend par l\`a que 
$C$ est $k$-r\'eguli\`ere en au moins un point $x_0$, donc au point g\'en\'erique,
et n'est $(k+1)$-r\'eguli\`ere en aucun point.
Soit $t\mapsto x(t)$ une param\'etrisation r\'eguli\`ere 
locale de $C$ avec $x(0)=x_0$. Par hypoth\`ese, pout $t\in \C$ voisin de $0$,
l'espace vectoriel $E(t)$ engendr\'e par 
$x'(t),\ldots,x^{(k)}(t)$ est de dimension $k$ et on peut \'ecrire 
$$
x^{(k+1)}(t) = \sum_{j=1}^k \l_j(t)x^{(j)}(t).
$$
En d\'erivant et par r\'ecurrence, on obtient que,
pour tout $j\in \N^\star$, $x^{(j)}(t)\in E(t)$
si $t$ est assez petit. En particulier $x^{(j)}(0)\in E(0)$ pour tout $j\geq 1$ et,
par analyticit\'e, $C\subset E(0)$. On a donc :

\sk
{\em La dimension de l'espace engendr\'e par une courbe irr\'eductible $C\subset \P^N$
est \'egale \`a sa r\'egularit\'e osculatrice.}

\sk
Plus g\'en\'eralement on a :
\ble
\label{osc1}
Soit $k\geq 1$ la r\'egularit\'e osculatrice d'une courbe alg\'ebrique 
irr\'eductible 
$C\subset \P^N$ et $k_1,\ldots,k_m\in \N$
des entiers tels que $k + 1 = \sum_{i=1}^m (k_i+1)$.
Pour tout $m$-uplet $(x_1,\ldots,x_m)\in C^m$ g\'en\'erique,
on a $\lan C\ran = \oplus_{i=1}^m C_{x_i}(k_i)$.

Si $C$ est une courbe rationnelle normale de degr\'e $k$, on a la m\^eme 
conclusion d\`es que les points $x_1,\ldots,x_m$ sont deux-\`a-deux distincts.
\ele
\bpf
On a d\'ej\`a \'etabli le r\'esultat pour $m=1$.
On obtient le cas g\'en\'eral par r\'ecurrence, en projetant 
$C$ sur un hyperplan depuis un point $a\in C$, tel
que $C$ est $k$-r\'eguli\`ere en $a$, et en appliquant les 
lemmes pr\'ec\'edents et l'hypoth\`ese de r\'ecurrence.
La d\'emonstration de la seconde partie est analogue, compte tenu du fait
qu'une courbe rationnelle normale $C$ de degr\'e $k$ est $k$-r\'eguli\`ere 
en chacun de ses points et que sa projection depuis l'un de ses points
est une courbe rationnelle normale de degr\'e $k-1$.
\epf

\ble
Soit $X$ un germe de vari\'et\'e lisse et $k$-r\'egulier en $x\in \P^N$
et $Y\subset X$ un germe en $x$ de vari\'et\'e de dimension $s$.
L'espace $\lan Y\ran \cap X_x(k)$
est de dimension $\geq {s+k \choose s}-1$.
\ele
\bpf
Soit $(x_\nu)_{\nu\in \N}$ une suite de points de 
$Y_{\rm reg}$ qui tend vers $x$. Pour $\nu$ assez grand, $X$ est $k$-r\'egulier en $x_\nu$,
donc $Y$ l'est aussi. L'espace $\lan Y\ran$ contient, pour tout $\nu$,
l'espace $Y_{x_\nu}(k)$, qui est de dimension ${s+k \choose s}-1$.
Comme toutes les valeurs d'adh\'erence de la suite $(Y_{x_\nu}(k))_{\nu\in \N}$
dans la grassmannienne ad\'equate sont contenues dans $\lan Y\ran \cap X_x(k)$,
on obtient le r\'esultat.
\epf

}


\begin{thebibliography}{aa}



\bibitem{Bo} E. Bompiani, {\em Propriet\`a differenziali carracteristiche di enti algebrici},
Rom. Acc. L. Mem. {\bf 26} (1921), 452--474.

\bibitem{CG} S.-S. Chern, P. A.  Griffiths, {\em An inequality for the rank of a web and 
webs of maximum rank}, Ann. Scuola Norm. Sup. Pisa {\bf 5} (1978), 539--557.

\bibitem{En} F. Enriques, {\em Una questione sulla linearit\`a dei sistemi di curve 
appartenenti ad una superficie algebrica}, Rom. Acc. L. Rend. {\bf 5} (1893), 3--8.

\bibitem{Fu} W. Fulton, Intersection Theory (Second Edition), {\em Springer}, 1998.


\bibitem{Gi} S.G. Gindikin, {\em Integral geometry, twistors and generalized conformal 
structures}, J. Geom. Phys. {\bf 5} (1988), 19-35.

\bibitem{Go}  V.V. Goldberg, Theory of multicodimensional 
$(n+1)$-web, {\em Kluwer Academic Publishers, Dordrecht}, 1988.

\bibitem{Gr} P. A.  Griffiths, {\em Variations on a theorem of Abel}, Invent. Math. {\bf 35} (1976), 321--390.

\bibitem{Han} T. Hangan, {\em Sur l'int\'egrabilit\'e des structures tangentes 
produits tensoriels r\'eels}, Ann. Mat. Pura Appl. {\bf 126} (1980), 149--185.


\bibitem{Ha} J. Harris, {\em A bound on the geometric genus of projective 
varieties}, Ann. Scuola Norm. Sup. Pisa {\bf 8} (1981), 35--68.

\bibitem{Ko} K. Kodaira, {\em A theorem of completeness of characteristic systems
for analytic families of compact submanifolds of complex manifolds},
Annals of Math. {\bf 75} (1962), 146--162.

\bibitem{Mu} D. Mumford, Algebraic Geometry I, Complex Analytic Varieties,
{\em die Grundlehren der mat. Wiss., Springer-Verlag}, 1976. 


\bibitem{PT} L. Pirio, J.-M. Tr\'epreau, {\em Sur l'alg\'ebrisation des tissus 
de rang maximal}, en pr\'eparation.

\bibitem{Ru} F. Russo, communication \`a L. Pirio, 2008.


\bibitem{Tr} J.-M. Tr\'epreau, manuscrit 2006 et 
{\em Une nouvelle caract\'erisation des vari\'et\'es de Veronese}, 13 pages (2010),
arXiv 1012.1008.  



\end{thebibliography}
\end{document}